\documentclass[11pt]{article}
\usepackage{times}
\usepackage{setspace}
\usepackage{amssymb}

\usepackage{fontenc,fontaxes}
\usepackage{amsthm}
\usepackage{authblk}
\usepackage{amsmath,amsfonts,fdsymbol,enumerate}
\usepackage[inline,shortlabels]{enumitem}
\usepackage{natbib}
\usepackage{comment}
\usepackage{mathtools}
\usepackage{arydshln}
\usepackage{bbm}
\usepackage{xr}
\usepackage{tikz}
\usepackage{xcolor}
\usepackage[colorlinks = true,
linkcolor = red,
urlcolor  = blue,
citecolor = blue,
anchorcolor = blue]{hyperref}
\usepackage[greek,english]{babel}
\usepackage[margin=1in]{geometry}
\usepackage{sectsty}
\usepackage{titlesec}
\titlespacing*{\section}
{0pt}{2ex plus 1ex minus .2ex}{0.7ex plus .2ex}
\titlespacing*{\subsection}
{0pt}{2ex plus 1ex minus .2ex}{0.5ex plus .2ex}
\titlespacing*{\subsubsection}
{0pt}{1.4ex plus 1ex minus .2ex}{0.4ex plus .2ex}

\usepackage{enumitem}
\setlist[enumerate,1]{itemsep=1pt, topsep=4pt, partopsep=0pt}
\setlist[enumerate,2]{nosep}
\setlist[itemize,1]{itemsep=1pt, topsep=4pt, partopsep=0pt}
\setlist[itemize,2]{nosep}

\makeatletter
\newcommand{\oset}[3][0ex]{%
	\mathrel{\mathop{#3}\limits^{
			\vbox to#1{\kern-2\ex@
				\hbox{$\scriptstyle#2$}\vss}}}}
\makeatother

\expandafter
\def \expandafter \normalsize \expandafter{\normalsize \setlength \abovedisplayskip{5pt plus 2pt minus 3pt}}
\expandafter
\def \expandafter \normalsize \expandafter{\normalsize \setlength \abovedisplayshortskip{0pt plus 2pt}}
\expandafter
\def \expandafter \normalsize \expandafter{\normalsize \setlength \belowdisplayskip{5pt plus 2pt minus 3pt}}
\expandafter
\def \expandafter \normalsize \expandafter{\normalsize \setlength \belowdisplayshortskip{2pt plus 2pt}}

\newcommand*\circled[1]{\tikz[baseline=(char.base)]{
		\node[shape=circle,draw,inner sep=0.3pt] (char) {#1};}}

\newcommand{\RN}[1]{%
	\textup{\uppercase\expandafter{\romannumeral#1}}%
}

\newcommand{\dt}[1] {{#1}^{{\hspace{-0.15em}\smwhitestar}}}
\newcommand{\ddt}[1] {{#1}^{{\hspace{-0.15em}\smwhitestar\hspace{-0.25em}\smwhitestar}}}
\newcommand{\sbot}[1]{{#1}_{\perp \!\!\!\! \perp}}

\newtheoremstyle{slplain}
{.4\baselineskip\@plus.2\baselineskip\@minus.2\baselineskip}
{.4\baselineskip\@plus.2\baselineskip\@minus.2\baselineskip}
{\slshape}
{}
{\bfseries}
{.}
{ }
{}

\theoremstyle{slplain}

\newtheorem{cor}{Corollary}[section]
\newtheorem{lemma}{Lemma}[section]
\newtheorem{proposition}{Proposition}[section]
\theoremstyle{definition}
\newtheorem{remark}{Remark}[section]
\newtheorem{assumption}{Assumption}[section]
\newtheorem*{assumption*}{Assumption}
\newtheorem{definition}{Definition}[section]
\newtheorem{example}{Example}[section]

\mathtoolsset{mathic=true}

\DeclareMathOperator{\ran}{ran}

\DeclareMathOperator{\cl}{cl}

\DeclareMathOperator{\op}{op}

\DeclareMathOperator{\idb}{\it{I}}
\DeclareMathOperator{\idbp}{\it{I_p}}
\DeclareMathOperator{\Ann}{Ann}
\DeclareMathOperator{\rank}{rank}
\newcommand{\PPhi}{\widetilde{\Phi}}
\newcommand{\PP}{\widetilde{\mathrm{P}}}
\newcommand{\NN}{\widetilde{\mathrm{N}}}
\newcommand{\NR}{{\mathrm{N}}}
\newcommand{\PR}{{\mathrm{P}}}

\newcommand{\zero}{$0$}
\numberwithin{equation}{section}

\begin{document}
	
	
	
	\title{ Cointegration and Representation of  \\ Cointegrated Autoregressive Processes in Banach Spaces
	}
	
	
	\author{Won-Ki Seo}
	\affil{School of Economics, University of Sydney} 
	\maketitle	
	\vspace{-2em}
	\begin{abstract}
		We extend the notion of cointegration for time series taking values in a potentially infinite dimensional Banach space. Examples of such time series include stochastic processes in $C[0,1]$ equipped with the supremum distance and those in a finite dimensional vector space equipped with a non-Euclidean distance. We then develop versions of the Granger-Johansen representation theorems for I(1) and I(2) autoregressive (AR)  processes taking values in such a space. To achieve our goal, we first note that an AR($p$) law of motion can be characterized by a linear operator pencil via the companion form representation, and  then study the  spectral properties of a linear operator pencil to obtain a necessary and sufficient condition for a given AR($p$) law of motion to admit I(1) or I(2) solutions. These operator-theoretic results form a fundamental basis for our representation theorems. Furthermore, it is shown that our operator-theoretic approach is in fact a closely related extension of the conventional approach taken in a Euclidean space setting.  Our theoretical results may be especially relevant in a recently growing literature on functional time series analysis in Banach spaces. 	
	\end{abstract} \vspace{0.0em}
	

	
	
	\section{Introduction}\label{sintro}
	Conventionally, the subject of time series analysis is time series taking values in finite dimensional Euclidean space. On the other hand,  a recent literature on functional time series analysis  deals with time series taking values in a possibly infinite dimensional Banach or Hilbert space, for instance, those in $C[0,1]$ equipped with the supremum norm. Examples of such time series are not restricted to function-valued stochastic processes: those in a finite dimensional vector space equipped with a non-Euclidean metric, such as e.g.\ Chebyshev distance or taxicab distance, are also included. 
	
	The property of cointegration, which was introduced by \cite{granger1981} and has been studied in Euclidean space, was recently extended to a more general setting.  A recent paper by \cite{Chang2016152} appears to be the first to consider the possibility of cointegration in an infinite dimensional Hilbert space. More recently, \cite{BSS2017} adopted the notion of cointegration from \cite{Chang2016152} and provided a rigorous treatment of cointegrated linear processes taking values in Hilbert spaces.   
	
	The Granger-Johansen representation theorem is the result on the existence and representation of I(1) (and I(2)) solutions to a given autoregressive (AR) law of motion. Due to crucial contributions by e.g. \cite{engle1987},  \cite{Johansen1991,johansen1992representation,Johansen1996,johansen2008}, \cite{Schumacher1991}, \cite{Hansen2005}, \cite{faliva2002partitioned,Faliva2010,Faliva2011,faliva2021cointegrated}, and \cite{franchi2016,Franchi2017a},  much on this subject is already well known in a Euclidean space setting. For a brief historical overview of this topic, see the introduction of \cite{BS2018}. More recently, \cite{chang2016},  \cite{Hu2016} and \cite{BSS2017} extended the Granger-Johansen representation theorem for the I(1) case to a more general Hilbert space setting. Furthermore, \cite{BS2018} provided representation theorems for I(1) and I(2) AR processes in such a setting based on analytic operator-valued function theory, and \cite{Franchi2017b} developed a more general result for I(d) AR processes for $d\geq1$. 
	

	This paper provides a suitable notion of cointegration and extends the Granger-Johansen representation theorem for Banach-valued, not necessarily Hilbert-valued, AR processes that are I(1) or I(2); that is,  our theory can be applied to more general AR processes, for instance, those taking values in $C[0,1]$,   $L^q[0,1]$ for $1\leq q <\infty$, or any finite dimensional vector space equipped with an arbitrary norm.   Viewed in the light of our purpose,  our representation theorems need to be developed without relying on the following two preconditions commonly required in the literature: (i) a Hilbert space structure and (ii) a special restriction on the  AR polynomial. To see this in detail, we briefly  review the relevant literature. For a given AR($p$) law of motion in a Hilbert space, which is characterized by the AR polynomial $\Phi(z)=I-z\Phi_1 - \ldots -z^p \Phi_p$, \cite{BSS2017} assume that $\Phi_1,\ldots,\Phi_p$ are compact operators when $p>1$ (this compactness assumption is not required if $p=1$), and provide  a sufficient condition for the existence of I(1) solutions and a  characterization of such solutions. In their representation theory, compactness of $\Phi_1,\ldots,\Phi_p$ makes $\Phi(z)$ belong to a special subclass of linear operators, called Fredholm operators, and the mathematical properties of such operators play an important role.  More representation theorems in a Hilbert space setting are provided by  \cite{Hu2016}, \cite{BS2018} and \cite{Franchi2017b}, among which the latter paper  more generally deals with I(d) AR($p$) processes for $d\geq1$ and $p \geq 1$. The representation theorems in those papers are closely related to that provided by \cite{BSS2017} in the sense that Fredholmness of $\Phi(z)$ has a crucial role in their developments. The Fredholm assumption explicitly or implicitly employed in the foregoing papers turns out to  place nontrivial restrictions on solutions to the AR($p$) law of motion:  nonstationarity of such a solution is  driven by a necessarily finite dimensional unit root process even in an infinite dimensional setting.  From another standpoint, \cite{chang2016} employ a different assumption that $\Phi(1)$ is a compact operator and provide an I(1) representation result. As opposed to the results under Fredholmness of  $\Phi(z)$, it turns out that their compactness assumption always leads to I(1) solutions associated with an infinite dimensional unit root process unless the considered Hilbert space is finite dimensional.  To briefly sum up, all of these existing versions are developed in a Hilbert space setting, and each of those relies on a special requirement about $\Phi(z)$ which restricts solutions to the AR($p$) law of motion in a specific way. We thus need a novel approach to overcome these limitations in our more general setting.
	
	To accomplish our goal, we first introduce a suitable notion of cointegration in Banach spaces by defining a  cointegrating functional that properly generalizes the conventional notion of a cointegrating vector. We then characterize the cointegrating space (to be defined as the collection of cointegrating functionals) based on the Phillips-Solo device \citep{phillips1992asymptotics} applied to operator-valued functions. After doing that, representation theorems for I(1) and I(2) AR processes taking values in a Banach space are provided. Our representation theory   is derived under more primitive and weaker mathematical conditions in a general Banach space setting where we do not even have the notion of an angle (inner product)  between two vectors. From \cite{Johansen1991,johansen1992representation} to the foregoing recent papers, it seems that some geometrical properties induced by an inner product, such as orthogonality,  have been thought to be essential for the representation theory. However, it will be clarified in this paper that such a richer geometry is not necessarily required: in our representation theory, geometrical properties induced by an inner product have no essential role. 	
	
	To obtain our representation theorems, we first note that an AR($p$) law of motion in a Banach space allows the companion form AR(1) representation  in a properly defined product Banach space, and it is thus characterized by a linear operator pencil (to be introduced in detail later) denoted by $\widetilde{\Phi}(z)$. By studying the spectral properties of a linear operator pencil, we find necessary and sufficient conditions for  $\widetilde{\Phi}(z)^{-1}$ to have a pole of order 1 and 2, and also obtain a local  characterization of $\PPhi(z)^{-1}$ near $z=1$. These operator-theoretic results not only determine the integration order of solutions to the AR($p$) law of motion, but also lead us to a representation of such solutions in terms of the behavior of $\widetilde{\Phi}(z)$ around $z=1$; that is, our versions of the Granger-Johansen representation theorems for I(1) and I(2)  AR processes are obtained. The fact that solutions to the AR($p$) law of motion are characterized in terms of a local behavior of $\widetilde{\Phi}(z)$,  rather than the original AR polynomial $\Phi(z)$, makes it difficult to compare our results to those developed in a Hilbert/Euclidean space setting.   We thus provide  further representation results so that important characteristics of I(1) or I(2) solutions are expressed in terms of linear operators associated with $\Phi(z)$. To this end, we first show that our necessary and sufficient condition for $\PPhi(z)^{-1}$ to have a pole of order 1 (resp.\ 2)  is equivalent to a natural generalization of the well known Johansen I(1) (resp.\ I(2)) condition, and then we recharacterize solutions to the AR($p$) law of motion in a desired way using those equivalent conditions. By this further effort, not just we can obtain a more detailed characterization of such solutions but also we can better clarify the connection between our representation results and those in the existing literature. 
	
	We structure the remainder of the paper as follows.  In Section \ref{sfcointeg}, we develop a suitable notion of cointegration in Banach spaces and provide some related results. Our representation theory for I(1) and I(2) AR processes are contained in Section \ref{srep}, and concluding remarks follow in Section \ref{sconclude}. Appendix \ref{app:prelim} reviews background material for our study, and Appendix \ref{appenA} collects the proofs of our main results.

	\section{Cointegration in Banach spaces} \label{sfcointeg}
	Let $\{X_t\}_{t\geq 0}$ be an I(1) time series  taking values in Euclidean space of dimension $n$, denoted by $\mathbb{R}^n$. If there exists a nonzero vector $\beta \in \mathbb{R}^n$ such that $\{\beta^\intercal X_t\}_{t\geq0}$ is stationary under a suitable choice of $X_0$, we then say that 	$\{X_t\}_{t\geq 0}$ is cointegrated with respect to $\beta$,  and call $\beta$ a cointegrating vector; see e.g.\ \citet[Definition 3.4]{Johansen1996}. In this conventional definition of cointegration, $\beta$ itself  acts as a scalar-valued map defined on $\mathbb{R}^n$, and what makes $\beta$ a cointegrating vector is \textit{stationarity} of scalar-valued time series $\{\beta^\intercal X_t\}_{t\geq 0}$. We thus may understand cointegration as a property of scalar-valued maps defined on $\mathbb{R}^n$, which leads to the following alternative definition.
	\begin{definition} \label{defcointeg}
		For an I(1) time series $\{X_t\}_{t\geq 0}$ and any scalar-valued linear map $f$ defined on $\mathbb{R}^n$ (i.e., functional on $\mathbb{R}^n$), if $\{f(X_t)\}_{t\geq 0}$ can be stationary under  a suitable choice of $X_0$, then we say that $\{X_t\}_{t\geq 0}$ is cointegrated with respect to $f$, and call $f$ a cointegrating functional.  
	\end{definition}
	\noindent In fact, the above definition is equivalent to the conventional one due to the Riesz representation theorem \citetext{see e.g.,\ \citeauthor{Conway1994},  \citeyear{Conway1994}, p.\ 13}, implying that any functional $f$ on $\mathbb{R}^n$ is uniquely identified as a vector $\beta$ in the following sense: $f(x) = \beta^\intercal x$ for all $x \in \mathbb{R}^n$. Nevertheless, defining cointegration as in Definition \ref{defcointeg} is advantageous especially when we consider a more general vector space; as will be shown, we may replace $\mathbb{R}^n$ with a separable complex Banach space $\mathcal B$ without a serious theoretical complication, and then may obtain a suitable notion of cointegration in $\mathcal B$.
	
	Throughout this section, we formally introduce cointegrated I($1$) and I($2$) processes taking values in  a Banach space and characterize the collection of cointegrating functionals. Prior to a detailed mathematical treatment of those, it may be helpful to see an example of functional time series of economic or statistical interest that can motivate our more general setting.

	\subsection{Example : Banach-valued time series and cointegration}  \label{example}
	Our more general setting is of central relevance for applications involving functional time series. As mentioned and analyzed in \cite{horvath2010testing} and \cite{Hormann2013}, possibly one of the most natural functional time series is a sequence of intraday price curves of a financial asset. Let $X_t(s)$ be the price of a financial asset at time $s \in [s_{\min}, s_{\max}]$ on day $t\in \{1,2,\ldots\}$. By reparametrizing $s$ into $u=(s-s_{\min})/(s_{\max}-s_{\min})$, $X_t \coloneqq \{X_t(u), u \in [0,1]\}$ may be viewed as a random element in $C[0,1]$, the Banach space of continuous functions on $[0,1]$ equipped with the usual sup norm. Linear functionals defined on $C[0,1]$ reveal various characteristics of $X_t$. For example, consider $f_1$, $f_2$, and $f_3$ defined by 
	\begin{equation} \label{eqfns}
		f_1(x) = x(1), \quad\quad f_2(x) = \int_{0}^{1}x(u)du,  \quad\quad f_3(x) = x(1) -  \int_{0}^{1}x(u)du,
	\end{equation}
	where $x \in C[0,1]$. Then, $f_1(X_t)$ (resp.\ $f_2(X_t)$) computes the closing  (resp.\ average) price on day $t$, and  $f_3(X_t)$ gives their difference.
	
	We may assume that price curves observed on two adjacent days,  $X_{t-1}$ and $X_{t}$, are tightly connected in the sense that $X_{t}(0) = X_{t-1}(1)$ or loosely connected in the sense that an overnight jump $\nu_t(0) \coloneqq X_{t}(0) - X_{t-1}(1)$ is allowed. In either case, $\{X_t(u)-X_{t-1}(1)\}_{u \in [0,1]}$ denotes cumulative intraday returns on day $t$. For illustrative purposes, we may model such a sequence returns as follows: for some stationary process $\{\nu_t\}_{t\geq 0}$,\footnote{An empirical evidence about stationarity of cumulative stock return curves of a financial asset was provided in \cite{horvath2010testing}; however, we need to be careful in interpreting such evidence in our context since their results are obtained by viewing intraday price curves as random elements of the usual  Hilbert space $L^2[0,1]$ rather than $C[0,1]$. } 
	\begin{equation*} 
		X_t(u) -  X_{t-1}(1) =  \nu_t (u), \quad u \in [0,1], 
	\end{equation*}
	where $\{\nu_t(1)\}_{t\geq 0}$ is assumed to have a positive variance for reasons to become apparent. By introducing a linear operator $\Phi_1$ defined by $\Phi_1(x)(u) = x(1)$ for $u \in [0,1]$ and then  suppressing dependence on $u$ for convenience, \eqref{eqmodel1} can be written as a curve-valued AR(1) process as follows,
	\begin{equation}\label{eqmodel1}
		X_t = \Phi_1 X_{t-1} + \nu_t.
	\end{equation} 
	If we take the functional $f_1$ to both sides of \eqref{eqmodel1}, we obtain  
	\begin{equation*}
		X_{t}(1) =  X_{t-1}(1) + \nu_t(1). 
	\end{equation*}
	That is, the time series of closing prices $\{f_1(X_t)\}_{t\geq 0}$ is a random walk with stationary increments. Therefore, $\{X_t\}_{t\geq 0}$ is a nonstationary curve-valued process; if it were stationary, $\{f_1(X_t)\}_{t\geq 0}$ would be stationary since $f_1$ is a continuous linear transformation. On the other hand,  note that 
	\begin{equation*}
		f_3(X_t) =  f_3(\nu_t),
	\end{equation*}
	from which we find that $\{f_3(X_t)\}_{t\geq 0}$ is stationary since $f_3$ is a continuous linear transformation and $\{\nu_t\}_{t\geq 0}$ is stationary. That is, $f_3$ transforms the curve-valued nonstationary process $\{X_t\}_{t \geq 0}$ into a scalar-valued stationary process. In view of Definition \ref{defcointeg}, $f_3$ may be called a cointegrating functional, which is,  of course, informal at this point since we have not yet provided our formal definition of a cointegrating functional.
	
	One may be interested in describing the above characteristics of the model \eqref{eqmodel1} using the existing theory of cointegration, assuming that the intraday price curves are random elements of the usual Hilbert space  $L^2[0,1]$, the space of square integrable functions on $[0,1]$ equipped with inner product $\langle f,g \rangle = \int_{0}^{1} f(u)g(u)du$. In this case, however, the AR(1) operator $\Phi_1$ given in \eqref{eqmodel1} and the functionals $f_1$ and $f_3$ given in \eqref{eqfns} lack an essential continuity property required in the existing theory, and dealing with those linear maps in this context is far beyond what has been covered in the literature.\footnote{$\Phi_1$, $f_1$ and $f_3$  are not continuous with respect to the topology of $L^2[0,1]$. Such a linear map is equivalently said to be unbounded. As far as this author knows, unbounded linear operators or functionals have not been considered in the existing theory of cointegration and the Granger-Johansen representation due to serious technical difficulties in dealing with those.} This example, therefore, shows that our general Banach space setting is useful to  accommodate more various functional time series as subjects of the theory of cointegration.
	
	Moreover, Banach space methodology is sometimes naturally in demand when researchers want to adopt a different notion of distance for functional time series analysis. For example, \cite{dette2017} noted that two curves with rather different visual shapes may still have a small $L^2$-distance and thus be identified as similar in the usual $L^2[0,1]$ setting. They therefore employed the sup-distance, which is expected to better reflect the visualization of curve-valued observations in statistical analysis, by assuming that such observations are random elements of $C[0,1]$.

	\subsection{Notation} \label{sec:notation}
	We here review our notation for the subsequent discussions.  The setting for our analysis is a separable complex Banach space $\mathcal B$ equipped with norm $\|\cdot\|_{\mathcal B}$. To conveniently introduce our notation, we let $\widetilde{\mathcal B}$ denote another such space, equipped with norm $\|\cdot\|_{\widetilde{\mathcal B}}$,  which can be set to, for example, $\mathcal B$ or the complex plane $\mathbb{C}$, or the topological dual  $\mathcal B'$  to be defined. 
	
	A linear operator $A:\mathcal B\mapsto\widetilde{\mathcal B}$ is said to be bounded if $\|Ax\|_{\widetilde{\mathcal B}} \leq M\|x\|_{\mathcal B}$ for some $M<\infty$. Such an operator is obviously continuous on $\mathcal B$. Unless otherwise noted, every linear operator considered in this paper is bounded. 
	Let $\mathcal L(\mathcal B,\widetilde{\mathcal B})$ denote the space of bounded linear operators from  $\mathcal B$ to $\widetilde{\mathcal B}$ equipped with the  operator norm $\|A\|_{\op} = \sup_{\|x\|\leq 1} \|Ax\|_{\widetilde{\mathcal B}}$. We are mostly concerned with the case $\mathcal B = \widetilde{\mathcal B}$, so let $\mathcal L(\mathcal B)$ denote $\mathcal L(\mathcal B,\mathcal B)$, and let $I \in \mathcal L({\mathcal B})$ denote the identity operator acting on $\mathcal B$.  For any $A \in \mathcal L(\mathcal B,\widetilde{\mathcal B})$, we let $\ran A$ (resp.\ $\ker A$) denote the set $\{Ax : x \in \mathcal B\}$ (resp.\ $\{x \in \mathcal B : Ax=0\}$). Commonly, $\ran A$ (resp.\ $\ker A$) is called the range (resp.\ the kernel) of $A$.	If $\dim(\ran A)<\infty$, then $A$ is said to be a  finite rank operator. For any subspace $V \subset {\mathcal B}$, let $V'$ denote the space of bounded linear functionals from $V$ to $\mathbb{C}$ equipped with the operator norm, i.e., $V' = \mathcal L(V,\mathbb{C})$, which is commonly called the topological dual of $V$.

	For any subset $V$ of  $\widetilde{\mathcal B}$, let $\cl (V)$ denote the closure of $V$, i.e.,  union of $V$ and its limit points. For subspaces $V_1$ and $V_2$ of $\widetilde{\mathcal B}$, we let $V_1 + V_2$ denote the set $\{v_1+v_2: v_1\in V_1, v_2\in V_2\}$, which is called the algebraic sum of $V_1$ and $V_2$. If $V_1+V_2 =\widetilde{\mathcal B}$ and $V_1 \cap V_2 = \{0\}$, we then say that $\widetilde{\mathcal B}$ is the direct sum of $V_1$ and $V_2$, and write $\widetilde{\mathcal B} = V_1 \oplus V_2$. In this case, $V_1$ (resp.\ $V_2$) is said to be complemented by $V_2$ (resp.\  $V_1$), and $V_2$ (resp.\ $V_1$) is called a complementary subspace of $V_1$ (resp.\  $V_2$). These definitions can be extended for a finite collection of subspaces $V_1,V_2,\ldots,V_k$ in an obvious way. 	For any set $V \subset \widetilde{\mathcal B}$, we let  $\Ann(V)$ denote the annihilator of $V$, defined by the set  $\{f \in \widetilde{\mathcal B}\,' : f(x) = 0, \forall x\in V \}$ which turns out to be a closed subspace of $\widetilde{\mathcal B}\,'$ \citetext{\citeauthor{Fabian2010},  \citeyear{Fabian2010}, p.\ 56}. For any closed subspace $V$ of $\widetilde{\mathcal B}$,  we let $\widetilde{\mathcal B}/V$ denote the quotient space equipped with the quotient norm $\|\cdot\|_{\widetilde{\mathcal B}/V}$, which are briefly reviewed in Appendix \ref{apppre1}.
	
	The definitions of a $\mathcal B$-random variable $X$, its expectation $EX$, covariance $C_X$, and cross-covariance $C_{X,Y}$ with another $\mathcal B$-random variable $Y$ are given in Appendix \ref{apppre3}. We are mostly concerned with the collection of $\mathcal B$-valued random variables $X$ satisfying  $EX = 0$ and $E\|X\|_{\mathcal B}^2 < \infty$, which is denoted by $\mathfrak L^2(\mathcal B)$. For $X \in \mathfrak L^2(\mathcal B)$, we say that $X$ has a positive definite covariance if $fC_X (f) = 0$ implies that $f = 0$.

	\subsection{Cointegrated I(d) processes in Banach spaces}  \label{sec:id}
	Throughout this paper, we will need to consider I($d$) processes in $\mathcal{B}$ and in  $\mathbb{C}$ for $d \in \{1,2\}$, with innovations in $\mathcal B$, so it is convenient to define the I($d$) property with another separable complex Banach space $\widetilde{\mathcal B}$ as in Section \ref{sec:notation}. Our definition of the I($d$) property is adapted from \cite{BS2018} and \cite{Franchi2017b} for our more general setting. As a key building block for the I($d$) property, we first define the I(0) property. 
	\begin{definition} \label{defi0}
		A sequence $X=\{X_t\}_{t \geq t_0}$ in  $\mathfrak L^2(\widetilde{\mathcal B})$ is said to be I(0) if we may write 
		\begin{equation}\label{linearpro}
			X_t - E(X_t) = \sum_{j=0}^\infty \Theta_j \varepsilon_{t-j}, \quad t \geq t_0,
		\end{equation}	
		where $\{\varepsilon_t\}_{t\in \mathbb{Z}}$ is an iid sequence in $\mathfrak L^2({\mathcal B})$ with positive definite covariance $C_{\varepsilon_0}$ and $\{\Theta_j\}_{j\geq 0}$ is a sequence in $\mathcal L (\mathcal B, \widetilde{\mathcal B})$ satisfying $\sum_{j=0}^\infty \|\Theta_j\|_{\op} < \infty$ and $\sum_{j=0}^\infty \Theta_j \neq 0$.
	\end{definition}
	\begin{remark} \label{reminno}
		$\{\varepsilon_t\}_{t\in \mathbb{Z}}$ in Definition \ref{defi0} is an iid sequence in  $\mathfrak L^2({\mathcal B})$, which is called a strong $\mathcal B$-white noise \citep[p.\ 148]{Bosq2000}. The iid condition is imposed for simplicity, and the results to be developed  remain valid even if it is replaced by   stationarity in the covariance structure, i.e., $C_{\varepsilon_t}$ does not depend on $t$ and $C_{\varepsilon_{t},\varepsilon_{s}} =0$ for all $t$ and $s \neq t$. Such a sequence $\{\varepsilon_t\}_{t\in \mathbb{Z}}$ is called a weak $\mathcal B$-white noise \citep[p.\ 161]{Bosq2000}.
	\end{remark}
	As in a Euclidean space setting, \eqref{linearpro} may be conveniently expressed as 
	\begin{equation*}
		X_t - E(X_t) = \Theta(L) \varepsilon_t,
	\end{equation*}
	where $\Theta(z) = \sum_{j=0}^\infty \Theta_j z^j$ and $L$ denotes the lag operator. Note that $\Theta(\cdot)$ is an operator-valued function   on $\mathbb{C}$, which is called an operator pencil (see Appendix \ref{apppre4}); if  $\Theta(\cdot)$ is matrix-valued, it is specially called a matrix pencil. Based on the I(0) property given by Definition \ref{defi0}, we define the I($1$) and the I($2$) properties  as follows.

	\begin{definition} \label{defid} For $d \in  \{1,2\}$, a sequence in $\mathfrak L^2({\widetilde{\mathcal B}})$ is said to be I($d$) if its $d$-th difference is an I(0) process admitting a representation \eqref{linearpro} with $\{\Theta_j\}_{j \geq 0}$ satisfying $\sum_{j=1}^\infty j^{\,d}\, \|\Theta_j\|_{\op} < \infty$.
	\end{definition}
	Note that we require some summability conditions for I(1) and I(2) sequences, which are introduced for mathematical convenience in order to facilitate the use of the Phillips-Solo device in Section \ref{ssinteg}. 
	
	Cointegration of an I($d$) process in $\mathcal B$ may be defined by extending Definition \ref{defcointeg} in an obvious way. Let $X$ be an I($d$) process in $\mathcal B$, $f$ be an element of $\mathcal B'$, and $\Delta \coloneqq I-L$ be the difference operator. If the scalar-valued time series $\{f(\Delta^{d-1}X_t)\}_{t \geq 0}$ is stationary ($\Delta^0$ is understood as the identity operator)  in $\mathbb{C}$ for a suitable choice of $X_0$, we then say that $X$ is cointegrated  and call $f$ a cointegrating functional. Obviously, the collection of cointegrating functionals constitutes a subspace of $\mathcal B'$, so we call it the cointegrating space.

	\subsection{Characterization of the cointegrating space} \label{ssinteg} 
	
	In this section, we characterize the cointegrating space associated with I(1) or I(2) processes in $\mathcal B$. A key input to our results is the Phillips-Solo device  (\citeauthor{phillips1992asymptotics}, \citeyear{phillips1992asymptotics}, Lemma 2.1 and Section 4) for obtaining an algebraic decomposition of a linear filter into long-run and transitory components. Even if  the Phillips-Solo device was presented in \cite{phillips1992asymptotics} as a way to decompose matrix pencils when the usual matrix norm is considered,  it can be directly extended to our Banach space setting by just replacing matrix pencils (resp.\ the usual matrix norm) with operator pencils  (resp.\ the operator norm); no further changes are required  from their proofs. 
	
	For $d \in \{1,2\}$, let $X = \{X_t\}_{t \geq -d+1}$ be an I($d$) sequence in $\mathcal B$, admitting the following representation,
	\begin{equation} 
		\Delta^d X_t = \Theta(L) \varepsilon_{t}, \quad t \geq 1. \label{eqid0}
	\end{equation}
	Under the summability conditions given in Definition \ref{defid}, we may apply the Phillips-Solo device to obtain
	\begin{equation} \label{eqid2}
		\Theta(L) = \Theta(1) + \Delta \dt{\Theta}(L),
	\end{equation}
	where $\dt{\Theta}(L) = - \sum_{j=0}^\infty \dt{\Theta}_j L^j$ and $\dt{\Theta}_j = \sum_{k=j+1}^\infty \Theta_k$.  In   \eqref{eqid2}, $\Theta(1)$ (resp.\ $\Delta \dt{\Theta}(L)$) is called the long-run (resp.\ transitory) component of $\Theta(L)$; see  \cite{phillips1992asymptotics}. We may deduce from \eqref{eqid2} that \eqref{eqid0} allows the following representation, called the Beveridge-Nelson decomposition: for $d \in \{1,2\}$ and for some $\tau_0$, 
	\begin{equation} \label{eqi1bn}
		\Delta^{d-1} X_t = \tau_0 + \Theta(1)\sum_{s=1}^t \varepsilon_s + \nu_t, \quad t \geq 0,
	\end{equation}
	where $\{\nu_t\}_{t \geq 0}$ is stationary and  $\nu_t = \dt{\Theta}(L) \varepsilon_t$ for each $t$. Given \eqref{eqi1bn}, the cointegrating space $\mathfrak C(X)$ associated with \(X\) is formally defined as follows: for $d \in \{1,2\}$,
	\begin{equation*}
		\mathfrak C(X) = \{ f \in \mathcal B' : \{f(\Delta^{d-1}X_t)\}_{t\geq 0} \text{ is stationary for some $\tau_0 \in \mathfrak L^2({\mathcal B})$} \}.
	\end{equation*}
	We also define 
	\begin{equation*}
		\mathfrak A(X) = \ran \Theta(1),
	\end{equation*}
	which is called the attractor space of $X$.  We then provide useful results to characterize $\mathfrak C(X)$ when (i) there is no restriction on $\mathfrak A(X)$ and (ii) the closure of $\mathfrak A(X)$ is complemented in $\mathcal B$,  i.e.,\  for some $V \subset \mathcal B$, \begin{equation} \label{eqdirecsum}
		\mathcal B = \cl (\mathfrak A(X)) \oplus V.
	\end{equation} 
	If $\mathcal B$ is a Hilbert/Euclidean space, then  $V = \mathfrak A(X)^\perp$ always satisfies \eqref{eqdirecsum}; however, there may not exist a subspace $V$ satisfying \eqref{eqdirecsum} if $\mathcal B$ is an infinite dimensional Banach space (Remark \ref{rem33}). Note also  that, even in a very simple case, a subspace $V$ satisfying \eqref{eqdirecsum} is not unique; for example, in the case where $\mathcal {B} = \mathbb{R}^2$ and $\cl (\mathfrak A(X))=\text{span}\{(1,0)\}$,  the span of any arbitrary vector that is not included in $\text{span}\{(1,0)\}$ can be $V$. However, the subsequent results do not depend on any specific choice of $V$, and hence $V$ can be arbitrarily chosen among possible candidates. For such  $V$,  we can define a unique projection $\PR_V \in \mathcal L(\mathcal B)$ onto $V$ along $\cl(\mathfrak A(X))$, which is defined by the following three properties: 
	\begin{equation}
		\PR_V  =  \PR_V^2, \quad\quad \ran \PR_V = V,  \quad \quad 	\ker \PR_V = \cl(\mathfrak A(X)); \label{projec} 
	\end{equation}
	see \citet[Theorem 3.2.11]{megginson1998}. Our first result given below characterizes $\mathfrak C(X)$ in terms of $\mathfrak A(X)$ and $\PR_V$. 
	\begin{proposition} \label{propbasic} If $X = \{X_t\}_{t \geq -d+1}$ is I($d$) for $d \in \{1,2\}$, the following hold.
		\begin{enumerate}[\normalfont(i)]
			\item 	$\mathfrak C(X)= \Ann(\mathfrak A(X))$.  \label{propbasic1}
			\item	If \eqref{eqdirecsum} is satisfied for some $V\subset \mathcal B$, then $\mathfrak C(X) = \{f \circ \PR_V : f \in  \mathcal B'\}$, 		where $\PR_V$ satisfies \eqref{projec}.  \label{propbasic2}
		\end{enumerate}
	\end{proposition}
	Proposition \ref{propbasic}-\ref{propbasic1} shows that $\mathfrak C(X)$ is given by the annihilator of $\mathfrak A(X)$, and thus a closed subspace of $\mathcal B'$ regardless of whether $\mathfrak A(X)$ is closed or not. Moreover, under the direct sum condition \eqref{eqdirecsum}, Proposition \ref{propbasic}-\ref{propbasic2} says that $\mathfrak C(X)$, as a subspace of $\mathcal B'$, is fully characterized by the projection $\PR_V \in \mathcal L(\mathcal B)$, which in turn leads to a natural decomposition of $\{\Delta^{d-1}X_t\}_{t\geq 0}$ into two components with different kinds of cointegrating behavior; see Remark \ref{rem3} below. 
	
	\begin{remark} \label{rem3} 
		In our Banach space setting, the cointegrating space $\mathfrak C(X)$ is, by definition, a subspace of $\mathcal B'$ which is in general different from $\mathcal B$. However, the result given in Proposition \ref{propbasic}-\ref{propbasic2} makes it possible to understand  $\mathfrak C(X)$ as a subspace of $\mathcal B$.  Consider the Beveridge-Nelson decomposition \eqref{eqi1bn}.  Using $\PR_V$ defined under the direct sum  \eqref{eqdirecsum}, we may decompose $\Delta^{d-1}X_t$ into $(\idb-\PR_V) \Delta^{d-1} X_t$ and $\PR_V \Delta^{d-1} X_t$. The former is the unit root component in the sense that $\{f(\idb-\PR_V)\Delta^{d-1}X_t\}_{t \geq 0}$ cannot be stationary for all $f \in \mathcal B'$ as long as $f(\idb-\PR_V) \neq 0$ while the latter is the stationary component in the sense that $\{f\PR_V \Delta^{d-1}X_t\}_{ t \geq 0}$ can be made stationary under a suitable choice of $\tau_0$ for all $f \in \mathcal B'$. This projection-based decomposition of a cointegrated time series is what has been done in a Euclidean space setting \citetext{see e.g.\ \citeauthor{Johansen1996},  \citeyear{Johansen1996}, pp.\ 40--41}, and as discussed, it is also possible in our setting without a richer geometry of a Hilbert space.
	\end{remark} 
	
	\begin{remark}	 \label{rem33}
		If \(\mathcal B\) is an infinite dimensional Banach space, a closed subspace may not be complemented \citep[pp.\ 301--302]{megginson1998}, hence \eqref{eqdirecsum} is  not  generally true. It, however, turns out that either of the following is a sufficient (but not necessary) condition for the existence of $V$ satisfying \eqref{eqdirecsum},
		\begin{equation*}
			\text{(i) } \dim(\mathfrak A(X)) < \infty, \quad \quad \text{(ii) }  \dim(\mathcal B/\mathfrak A(X)) < \infty;
		\end{equation*}
		see Theorem 3.2.18 in \cite{megginson1998}. The two conditions lead to different dimensionalities of the cointegrating space of $\mathcal B'$.  
		In case (i), $V$ is necessarily infinite dimensional, and we deduce from Proposition \ref{propbasic} that the cointegrating space is also infinite dimensional. In case (ii), on the other hand, $V$ is finite dimensional, hence the cointegrating space is finite dimensional as well.
	\end{remark}

	One may be interested in how the general results given by Proposition \ref{propbasic} reduce to what we have known about cointegration in a Hilbert/Euclidean space. Let $\mathcal H$ be  a separable complex Hilbert space with inner product \(\langle \cdot, \cdot \rangle\). If $\mathcal B = \mathcal H$, the Riesz representation theorem \citetext{see e.g.\ \citeauthor{Conway1994},  \citeyear{Conway1994}, p.\ 13} implies that every $f \in \mathcal H'$ is given by the map $\langle \cdot , y\rangle : \mathcal H \mapsto \mathbb{C}$ for a unique element $y \in \mathcal H$. Therefore, we may alternatively define the cointegrating space as follows: for $d \in \{1,2\}$, \begin{equation} \label{eqalterdef}
		\mathfrak C_{\mathcal H}(X) = \{ y \in \mathcal H : \{\langle \Delta^{d-1}X_t, y \rangle \}_{t\geq 0} \text{ is stationary for some $\tau_0 \in \mathfrak L^2({\mathcal H})$} \}.
	\end{equation}
	Moreover, in this case, the closure of any subspace is complemented by its orthogonal complement, i.e.,\  the direct sum \eqref{eqdirecsum} holds for $V= \mathfrak A(X)^\perp$ \citep[pp.\ 35--36]{Conway1994}. For this choice of $V$, $\PR_V$ becomes an orthogonal projection. Under all these simplifications, Proposition  \ref{propbasic} reduces to the following characterization, which is identical to the description of $\mathfrak C_{\mathcal H}(X)$ given by \cite{BSS2017}. 	
	\begin{cor} \label{propbasic3} If $X = \{X_t\}_{t \geq -d+1}$ is I($d$) for $d \in \{1,2\}$ and $\mathcal B = \mathcal H$, then  $\mathfrak C_{\mathcal H}(X) = \mathfrak A(X)^\perp$.
	\end{cor} 
	
	We close this section with some remarks on Proposition \ref{propbasic} and Corollary \ref{propbasic3}.
	\begin{remark}
		Suppose that $\mathcal B = \mathbb{R}^n$ or $\mathbb{C}^n$ equipped with the usual inner product $\langle x,y \rangle = x^\intercal y$. This is of course a special case of a Hilbert space, so we may consider the alternative definition of the cointegrating space, $\mathfrak C_{\mathcal H}(X)$, given in \eqref{eqalterdef}. If there exists a nonzero cointegrating vector,  the long-run component $\Theta(1)$ from the Phillips-Solo decomposition in this setting is a reduced rank matrix, i.e., $\rank \Theta(1) = s < n$. If so, there are two full column rank $n\times s$ matrices $\theta_1$ and $\theta_2$ satisfying $\Theta(1) = \theta_1 \theta_2^\intercal$ \citetext{see e.g. \citeauthor{Faliva2010}, \citeyear{Faliva2010}, Theorem 3 in Chapter 1}.	As a result, $\mathfrak C_{\mathcal H}(X)$ is given by the collection of vectors that are orthogonal to the columns of $\theta_1$, which is an $(n-s)$-dimensional subspace of  $\mathbb{R}^n$ or $\mathbb{C}^n$.   
	\end{remark}

	\begin{remark}[Second-order cointegrating functionals] 	 \label{rempolycointeg}
		Under the summability requirement for the I(2) property in Definition \ref{defid}, we may apply the Phillips-Solo device to $\dt{\Theta}(L)$ in \eqref{eqid2}, and obtain 
		\begin{equation} \label{eqid3}
			\Theta(L) = \Theta(1) + \Delta \dt{\Theta}(1) + \Delta^2 \ddt{\Theta}(L),
		\end{equation}
		where $\ddt{\Theta}(L) = - \sum_{j=0}^\infty \ddt{\Theta}_j L^j$ and $\ddt{\Theta}_j = \sum_{k=j+1}^\infty \dt{\Theta}_k$. We then may deduce from \eqref{eqid3} that  \eqref{eqid0} with $d=2$ allows the following representation: for some $\tau_0$ and $\tau_1$, 
		\begin{equation*}
			X_t = \tau_0 + \tau_1 t + \Theta(1)  \sum_{r=1}^t\sum_{s=1}^r \varepsilon_s + \dt{\Theta}(1)\sum_{s=1}^t \varepsilon_s + \dt{\nu}_t, \quad t \geq 0,
		\end{equation*}
		where $\{\dt{\nu}_t\}_{t \geq 0}$ is stationary and  $\dt{\nu}_t = \ddt{\Theta}(L) \varepsilon_t$ for each $t$. Note that for any $f \in \mathfrak C(X) = \Ann(\ran \Theta(1))$, we have $f(X_t) = f(\dt{\Theta}(1) \sum_{s=1}^t \varepsilon_s) + f(\dt{\nu}_t)$ by assuming $f(\tau_0) =f(\tau_1) = 0$. Then it may deduced from a nearly identical argument used to prove Proposition \ref{propbasic}-\ref{propbasic1} that $	\{f(X_t)\}_{t \geq 0}$ is stationary if and only if $f\in \Ann(\ran \dt{\Theta}(1))$ also holds. To sum up, $\{f(X_t)\}_{t \geq 0}$ can be made stationary under a suitable choice of $\tau_0$ and $\tau_1$ for any $f\in \Ann(\ran \dt{\Theta}(1)) \cap \Ann(\ran \Theta(1))$. We call such $f$ as a second-order cointegrating functional, which will be of interest to us in Section \ref{srep} as an important aspect of  I(2) AR processes in $\mathcal B$.
	\end{remark}

	\section{Representation of I(1) and I(2) autoregressive processes}\label{srep}
	For fixed $p \in \mathbb{N}$, suppose that a sequence $\{X_t\}_{t\geq -p+1} \subset \mathfrak L^2(\mathcal B)$ satisfies the following AR($p$) law of motion:
	\begin{equation}
		\Phi(L) X_t = \varepsilon_t, \quad t \geq 1,\label{arlaw} 
	\end{equation}
	where  $\{\varepsilon_t\}_{ t \in \mathbb{Z}} \subset \mathfrak L^2(\mathcal B)$ is an iid sequence with positive definite covariance $C_{\varepsilon_0}$ and
	\begin{equation*}
		\Phi(z) = I - \sum_{j=1}^p \Phi_j z^j, \quad \Phi_1,\ldots,\Phi_p \in \mathcal L(\mathcal B).
	\end{equation*}	
	\noindent We let the operator pencil $\Phi:\mathbb{C}\mapsto \mathcal L(\mathcal B)$ be  called the AR polynomial. The iid condition imposed on $\{\varepsilon_t\}_{ t \in \mathbb{Z}}$ can be replaced by stationarity in the covariance structure without affecting the results to be developed in this section; see Remark \ref{reminno}. We hereafter say that $\Phi$ has a unit root if it satisfies  Assumption \ref{assulinear1} below, where the following notation is employed: $\sigma(\Phi)$ denotes the spectrum of $\Phi$ given by the set $ \{ z \in \mathbb{C} : \Phi(z) \text{ is not invertible}\}$ and $D_R$ denotes the open disk with radius \(R\) centered at $0 \in \mathbb{C}$. 
	\begin{assumption}[Unit root]\label{assulinear1} \hspace{0.1cm}
		\begin{enumerate}[\normalfont(a)]\setlength\itemsep{0em}\vspace{-0.0em}
			\item\( \sigma(\Phi) \cap D_{1+\eta} = \{1\}\) for some $\eta>0$ and $\Phi(1)\neq 0$. \label{assulinear1a} 
			\item \label{assulinear1b} $\ran\Phi(1)$ and $\ker \Phi(1)$ can be complemented, i.e., for some $[\ran \Phi(1)]_{\complement} \subset \mathcal B$ and $[\ker \Phi(1)]_{\complement}\subset \mathcal B$, \begin{equation}\mathcal B = \ran \Phi(1) \oplus [\ran \Phi(1)]_{\complement}  = \ker \Phi(1) \oplus [\ker \Phi(1)]_{\complement}. \label{direct1} \end{equation}  
		\end{enumerate}  
	\end{assumption}
	Note that we require $\ran \Phi(1)$ (resp.\ $\ker \Phi(1)$) to allow a complementary subspace $[\ran \Phi(1)]_{\complement}$ (resp.\ $[\ker \Phi(1)]_{\complement}$) in Assumption \ref{assulinear1}-\ref{assulinear1b}. In general,  $[\ran \Phi(1)]_{\complement}$ and $[\ker \Phi(1)]_{\complement}$ satisfying \eqref{direct1} do not uniquely exist, but the results to be developed in this section only require their existence. 
	If $\mathcal B$ is a Hilbert space  as a special case, any closed subspace can be complemented by its orthogonal complement, hence \eqref{direct1} holds for $[\ran \Phi(1)]_\complement=[\ran \Phi(1)]^\perp$ and $[\ker \Phi(1)]_\complement=[\ker \Phi(1)]^\perp$ as long as $\ran \Phi(1)$ is closed ($\ker \Phi(1)$ is necessarily closed since $\Phi(1) \in \mathcal L(\mathcal B)$). In the existing representation theorems developed in a Hilbert space setting,  closedness of $\ran\Phi(1)$ is implied by the employed assumptions, so \eqref{direct1} consequentially holds (see Remark  \ref{remi1a} to appear in Section \ref{sec:relation}).  In a general Banach space setting, on the other hand, \eqref{direct1} does not necessarily hold when $\ran \Phi(1)$ is closed; see the example given by \citet[Theorem 3.2.20]{megginson1998}. Nevertheless, Assumption \ref{assulinear1}-\ref{assulinear1b} may not be restrictive in general and, moreover, is a fairly weaker requirement which is strictly implied by the regularity conditions on $\Phi(z)$ employed in the existing literature; a more detailed discussion will be given in Remark \ref{remi1a}.
	
	In this setting, what we seek are (i) a necessary and sufficient condition under which the AR($p$) law of motion \eqref{arlaw} allows I(1) or I(2) solutions, and (ii) a characterization of such solutions; in the case $\mathcal B=\mathbb{R}^n$ or $\mathbb{C}^n$, these issues are dealt with in Johansen's representation theory. 	Hereafter, we conveniently say that a sequence $\{X_t\}_{t \geq -p+1}$ from  \eqref{arlaw} allows the Johansen I(1) representation if it satisfies the following: for $\tau_0$ depending on initial values of \eqref{arlaw}, a stationary sequence $\{\nu_t\}_{t \geq 0}\subset \mathfrak L^2(\mathcal B)$, and $\Upsilon_{-1} \in \mathcal L(\mathcal B)$, 
	\begin{equation}
		X_t = \tau_0 + \Upsilon_{-1} \sum_{s=1}^t \varepsilon_s  + \nu_t, \quad t \geq 0.\label{bndecom}
	\end{equation}  			
	We also say that $\{X_t\}_{t \geq -p+1}$ allows the Johansen I(2) representation if it can be represented as follows: for $\tau_0$ and $\tau_1$ depending on initial values of \eqref{arlaw}, a stationary sequence $\{\nu_t\}_{t \geq 0}\subset \mathfrak L^2(\mathcal B)$, and $\Upsilon_{-2}, \Upsilon_{-1} \in \mathcal L(\mathcal B)$, 
	\begin{equation}
		X_t = \tau_0 + \tau_1 t +   \Upsilon_{-2} \sum_{r=1}^t\sum_{s=1}^r \varepsilon_s +   \Upsilon_{-1}  \sum_{s=1}^t \varepsilon_s  + \nu_t, \quad t \geq 0. \label{bndecom3}
	\end{equation}  
	
	In the case $\mathcal B=\mathbb{R}^n$ or $\mathbb{C}^n$, \cite{Johansen1991,Johansen1996} shows that  a necessary and sufficient condition for the AR($p$) law of motion \eqref{arlaw} to allow I(1) solutions is given by  that
	\begin{equation}
		\sbot{\alpha}^\intercal \Phi^{(1)}(1) \sbot{\beta} \,\,\text{ is invertible,} \label{eqjohansen}
	\end{equation}
	where, if we let $\alpha$ and $\beta$ are full-rank $n \times r$ matrices satisfying $\ran \Phi(1)=\alpha \beta^\intercal$ for some $r<n$, then $\sbot{\alpha}$ (resp.\ $\sbot{\beta}$) is a full-rank $n \times (n-r)$ matrix whose columns are orthogonal to $\alpha$ (resp.\ $\beta$) and satisfies that $\sbot{\alpha}^\intercal\sbot{\alpha}=\sbot{\beta}^\intercal\sbot{\beta}=I_{n-r}$ (the identity matrix of dimension $n-r$). The condition given by \eqref{eqjohansen} is called the Johansen I(1) condition, under which a sequence $\{X_t\}_{t \geq -p+1}$ from \eqref{arlaw} allows the Johansen I(1) representation with a certain operator $\Upsilon_{-1}$; see e.g.\ \citet[Theorem 4.2]{Johansen1996}. A similar representation result for the I(2) case is also  given by  \citet[Theorem 4.6]{Johansen1996}. Extending these results to a general Banach space setting may not be done by a simple extension:  in general, $\mathcal B$ is not equipped with an inner product and, moreover, it can be infinite dimensional;  hence we cannot rely on some important geometrical properties and matrix algebraic results that are allowed in Euclidean space.  We thus need a novel approach that relies on neither of geometrical properties induced by an inner product nor finite dimensionality of $\mathcal B$. 

	As observed in an early contribution by \cite{Schumacher1991}, the I($d$) property of solutions to the AR($p$) law of motion, characterized by a matrix pencil $\Phi(z)$, in $\mathbb{R}^n$ is determined by the behavior of the inverse $\Phi(z)^{-1}$ around $z=1$. This is also true in our  Banach space setting, so our approach for developing representation theory for I(1) and I(2) AR processes essentially boils down to examining the inverse of the AR polynomial around $z=1$. As a way to achieve this goal, we first consider the companion form of \eqref{arlaw} characterized by a linear operator pencil $\PPhi$ to be defined later, and study the behavior of $\PPhi(z)^{-1}$ around $z=1$ based on the spectral theory of linear operator pencils given in e.g.\ \cite{Kato1995}, \cite{Gohberg2013}, \cite{albrecht2011necessary}, \cite{albrecht2019fundamental} and references therein.  We then recover the behavior of $\Phi(z)^{-1}$ around $z=1$ from that of $\widetilde{\Phi}(z)^{-1}$.

	It will be convenient to fix standard notation and terminology, based on Appendix \ref{apppre4} providing a brief introduction to operator pencils,  for the subsequent discussions. For any   operator pencil $A$ and its spectrum $\sigma(A) = \{z \in \mathbb{C} : A(z) \text{ is not invertible}\}$, we let $\rho(A)$ denote the set $\mathbb{C}\setminus\sigma(A)$, which is called the resolvent set of $A$. Now suppose that $A$ allows the Laurent series at $z=z_0$ as follows: for some $d \geq 0$,
	\begin{equation}
		A(z) = \sum_{j= -d}^\infty A_j (z-z_0)^j, \quad A_{-d} \neq 0.   \label{eqtermex}
	\end{equation} 
	If $d=0$, we say that $A(z)$ is holomorphic (or equivalently, complex-differentiable)  at $z=z_0$ and let  $A^{(j)}(z_0)$ denote the $j$-th complex derivative of $A(z)$ evaluated at $z=z_0$. In this case,  \eqref{eqtermex} becomes the Taylor series of $A(z)$ at $z=z_0$, which is specially called the the Maclaurin series of  $A(z)$ if $z_0 = 0$.   If $d\neq 0$, $A(z)$ is said to have an isolated singularity at $z=z_0$. An isolated singularity with $d < \infty$ is called a pole of order $d$. A pole of order one is said to be simple. If $d=\infty$, $A(z)$ is said to have an essential singularity at $z=z_0$. The sum of the leading terms indexed by $j=-d,\ldots,-1$ is the called the principal part and the sum of the remaining terms is called the holomorphic part.

	\subsection{Relations with the literature}\label{sec:relation}
	A few different versions of the Granger-Johansen representation theorem have been proposed in the recent literature on cointegrated functional time series taking values in a Hilbert space, such as \cite{chang2016},  \cite{Hu2016}, \cite{BSS2017}, \cite{BS2018} and \cite{Franchi2017b}. Compared to those, our versions are to be developed under a general Banach space setting without relying on a richer geometry of a Hilbert space, which can help us consider  more various functional time series as subjects of the theory of cointegration, as illustrated in Section \ref{example}. Apart from such mathematical gains, we here briefly describe how our setting is related to the assumptions employed in the aforementioned papers by assuming $\mathcal B = \mathcal H$ (recall that $\mathcal H$ denotes a separable complex Hilbert space).

	We first focus on the I(1) case. Except for the paper by \cite{chang2016} providing a quite different representation result, the AR polynomial $\Phi(z)$ in the foregoing papers satisfies the following  condition: for $z \in \mathbb{C}$,
	\begin{equation} \label{eqfred}
		\dim(\ker	\Phi(z)) < \infty \,\,\, \text{and}  \,\,\,  \dim ([\ran \Phi(z)]^\perp) < \infty. 
	\end{equation} 
	If \eqref{eqfred} holds, $\Phi(z)$ is called a Fredholm operator. Fredholmness of $\Phi(z)$ can be more generally defined in our Banach space setting by replacing the latter condition in \eqref{eqfred} with $\dim (\ran \Phi(z)/\mathcal B) <\infty$. The Fredholm property, combined with the unit root assumption (Assumption \ref{assulinear1}), produces some special behavior of $\Phi(z)^{-1}$ near $z=1$, which becomes a crucial input to the existing  theorems; see \citet[Appendix A.1]{BS2018} and \citet[Appendix B]{Franchi2017b}.  An important consequence of assuming \eqref{eqfred} is that $\Upsilon_{-1}$ in \eqref{bndecom} always becomes  a finite rank operator, hence the random walk component $\Upsilon_{-1} \sum_{s=1}^t \varepsilon_s$ in \eqref{bndecom} essentially boils down to a finite dimensional unit root process. As a result,  the attractor space (resp.\ the cointegrating space) associated with the AR($p$) law of motion is necessarily finite dimensional (resp.\ infinite dimensional). On the other hand, the version of \cite{chang2016} relies on the assumption that $\Phi(1)$ is compact, which turns out to result in the opposite case, where the cointegrating space is finite dimensional and the random walk component takes values in an infinite dimensional space  unless $\mathcal H$ is finite dimensional. Their compactness assumption  is not compatible with Fredholmness of $\Phi(z)$ in an infinite dimensional setting. We thus have two qualitatively different I(1) representation results depending on two generally incompatible regularity conditions on  $\Phi(z)$. 
	
	To the best of the author's knowledge, the existing representation theorems for I(2) AR processes in a  general Hilbert space setting were recently provided by \cite{BS2018} and \cite{Franchi2017b}, where Fredholmless of $\Phi(z)$ is an essential assumption for their representation theory. Similar to the I(1) case, the Fredholm assumption makes $\Upsilon_{-2}$ and $\Upsilon_{-1}$ in \eqref{bndecom3} become finite rank operators, hence the random walk component  $\Upsilon_{-2} \sum_{r=1}^t\sum_{s=1}^r \varepsilon_s +   \Upsilon_{-1}  \sum_{s=1}^t \varepsilon_s$  is intrinsically a finite dimensional unit root process. This requirement entailed by the Fredholm assumption  not only  compels  the attractor space associated with I(2) solutions to be finite dimensional, but also places some more restrictions on their cointegrating behavior; a more detailed discussion will be given in Section  \ref{sec:repi2}.

	As discussed above, any regularity conditions imposed on $\Phi(z)$ may compel solutions to the AR($p$) law of motion to  have some specific characteristics. It is thus desirable to develop representation theory for I(1) and I(2) AR processes under minimal regularity conditions on $\Phi(z)$. We in this paper require  weaker conditions on $\Phi(z)$ than either of Fredholmness or compactness, which naturally makes our representation theory place weaker restrictions on solutions to the AR($p$) law of motion.  As an example, the random walk component of I(1) or I(2) solutions can be either finite dimensional or infinite dimensional in our results; this is in contrast to that the component is required to be exclusively finite dimensional or infinite dimensional depending on the employed regularity condition on $\Phi(z)$ in the recent literature.  
	
	\begin{remark}\label{remi1a} 
		As discussed, in the existing literature, $\Phi(z)$ is either Fredholm or compact. Fredholm property \eqref{eqfred} implies that  $\ran \Phi(1)$ (resp.\ $\ker \Phi(1)$) allows a finite  (resp.\ an infinite) dimensional complementary subspace. This is also true in a more general situation where $\mathcal B$ is not necessarily a Hilbert space and $\Phi(z)$ is a Fredholm operator acting on $\mathcal B$; this can be shown from Remark \ref{rem33} and the fact that both $\ker\Phi(1)$ and $\mathcal B/\ran\Phi(1)$ are finite dimensional in this case. 	If $\Phi(z)$ is compact and satisfies Assumption \ref{assulinear1}-\ref{assulinear1a}, then $ \Phi(1)$ is necessarily a finite rank operator \citep[Lemma 1]{chang2016}. Thus, $\ran \Phi(1)$ (resp.\ $\ker \Phi(1)$)  allows an infinite (resp.\ a finite)  dimensional complementary subspace. Note that Fredholmness or compactness of $\Phi(z)$ places some specific dimensionality restrictions on $[\ran \Phi(1)]_\complement$ and  $[\ker \Phi(1)]_\complement$, while no such restrictions are required by Assumption \ref{assulinear1}-\ref{assulinear1b}. 
	\end{remark}		
	\begin{remark}\label{remi1b} 
		As will be shown in Proposition \ref{propjohanseni1}, the random walk component in the I(1) case always takes values in $\ker \Phi(1)$, whose dimension is finite (resp.\ infinite) if $\Phi(z)$ is Fredholm (resp.\ compact) under Assumption \ref{assulinear1}-\ref{assulinear1a}. This shows where the difference  between the existing I(1) representation results about the dimensionality of the random walk component originates from. 
	\end{remark}

	\subsection{Linearization of the AR polynomial}\label{linearization}
	Consider the product Banach space  $\mathcal B^p$ equipped with the norm $\|(x_1,\ldots,x_p)\|_{\mathcal B^p} = \sum_{j=1}^p \|x_j\|_{\mathcal B}$ for any $(x_1,\ldots,x_p) \in \mathcal B^p$. We let $I_p$ denote the identity map acting on $\mathcal B^p$. In fact, the AR($p$) law of motion \eqref{arlaw} may be understood as the following AR(1) law of motion in \(\mathcal B^p\):
	\begin{equation}\label{a1companion0}
		\PPhi(L) \widetilde{X}_t =  \widetilde{\varepsilon}_t,
	\end{equation}
	where	$\PPhi:\mathbb{C} \mapsto \mathcal L(\mathcal B^p)$ is a linear operator pencil given by 	$\PPhi(z) = I_p - z \PPhi_1$ and
	\begin{equation}
		\widetilde{X}_t=  \small \left[\begin{matrix} X_t  \\ X_{t-1}  \\  \vdots \\ X_{t-p+1}    \end{matrix} \right]\normalsize, \quad 
		\PPhi_1  = \small \left[\begin{matrix} {\Phi}_1 & {\Phi}_2 & \cdots & {\Phi}_{p-1} & {\Phi}_p \\ \idb & 0 & \cdots &0 & 0 \\ \vdots & \vdots& \ddots & \vdots &\vdots \\ 0 &0 & \cdots &  \idb & 0\end{matrix} \right]\normalsize, \quad \small \widetilde{\varepsilon}_t=  \left[\begin{matrix} \varepsilon_t  \\ 0  \\  \vdots \\ 0   \end{matrix} \right].\normalsize \label{a1companion}
	\end{equation}
	Commonly, \eqref{a1companion0} is called the companion form of \eqref{arlaw}; see e.g. \citet[p.\ 15]{Johansen1996} or \citet[p.\ 128, 161]{Bosq2000}. From a  mathematical point of view, the behavior of $\Phi(z)^{-1}$ that we want to know may be obtained from that of  $\PPhi(z)^{-1}$, which is as described in  Proposition \ref{propa1} below, where the following notation is employed: $\Pi_p:\mathcal B^p \mapsto \mathcal B$ and  $\Pi_p^\ast:\mathcal B \mapsto \mathcal B^p$ denote the maps defined by
	\begin{equation} \label{coorproj}\Pi_p(x_1,x_2,\ldots,x_p) = x_1, \quad \quad  \Pi_p^*(x_1) = (x_1,0,\ldots,0). \end{equation}
	
	\begin{proposition}\label{propa1} 	Under {Assumption} \ref{assulinear1}, the operator pencils $\PPhi$ and $\Phi$  satisfy the following.
		\begin{enumerate}[\normalfont(i)]\setlength\itemsep{0em}
			\item\label{propa11} $\sigma(\PPhi) = \sigma(\Phi)$ \,and\, $\Pi_p \PPhi(z)^{-1}\Pi_p^\ast = \Phi(z)^{-1}$. 
			\item\label{propa12} 	Under Assumption \ref{assulinear1}, if either of $\PPhi(z)^{-1}$ or $\Phi(z)^{-1}$ has a pole of order $d$ (resp.\ essential singularity) at $z=1$, then the other has a pole of order $d$ (resp.\ essential singularity) at $z=1$.  
		\end{enumerate}
	\end{proposition}
	Proposition \ref{propa1}-\ref{propa11} shows that $\PPhi(z)$ inherits the unit root property of $\Phi(z)$ given by  Assumption \ref{assulinear1}, and $\Phi(z)^{-1}$ can be recovered from $\PPhi(z)^{-1}$ using the maps given in \eqref{coorproj}. Moreover,  Proposition \ref{propa1}-\ref{propa12} implies that we can obtain a necessary and sufficient condition for $\Phi(z)^{-1}$ to have a pole of order 1 or 2 at $z=1$ by finding such a condition   for $\PPhi(z)^{-1}$. These results will become useful in the development of our representation theorems for I(1) and I(2) AR processes.


	\subsection{Representation of I(1) autoregressive processes}	\label{sec:repi1}
	In Section \ref{sec:repi1companion} we develop our representation theory for I(1) AR processes resorting to the companion form AR(1) representation \eqref{a1companion0}; this is done by studying the spectral properties of $\PPhi(z)$ under Assumption \ref{assulinear1}. We then discuss on how the results obtained via the companion form can be reformulated in terms of the behavior of the AR polynomial $\Phi(z)$ in Section \ref{sec:repi1noncompanion}. 
	
	
	\subsubsection{Representation via the companion form}	\label{sec:repi1companion}
	Resorting to the companion form \eqref{a1companion0}, we  in this section provide a necessary and sufficient condition for the AR($p$) law of motion \eqref{arlaw} to admit I(1) solutions and a characterization of such solutions.   
	
	Under Assumption \ref{assulinear1}-\ref{assulinear1a}, we know from the results given in Appendix \ref{apppre4} (especially, see \eqref{mopencil})  that $\PPhi(z)^{-1}$ can be written as  the following Laurent series: for $d \in \mathbb{N} \cup \{\infty\}$, 
	\begin{equation} 
		\PPhi(z)^{-1} =  - \sum_{j=-d}^{-1} \NN_j (z-1)^j  - \sum_{j=0}^\infty \NN_j (z-1)^j.  \label{laurent0}
	\end{equation} 	
	For notational convenience, we  let 
	\begin{equation} \label{defPP}
		\PP \coloneqq \NN_{-1} \PPhi_1. 
	\end{equation}
	The operator given above turns out to be a projection under Assumption \ref{assulinear1}  (Lemma \ref{lemprime}-\ref{lemprime2}) and has a crucial role in the subsequent discussion. 
	What we first pursue for the development of our representation theory is a necessary and sufficient condition under which the AR($p$) law of motion \eqref{arlaw} admits I(1) solutions (or equivalently, $\PPhi(z)^{-1}$ has a simple pole at $z=1$).  When $\mathcal B = \mathbb{R}^n$ or $\mathbb{C}^n$, a well known such condition is the Johansen I(1) condition given by  \eqref{eqjohansen}, and it plays an essential role in Johansen's representation theory for I(1) AR processes.  \cite{BSS2017}, who studied the same issue for the case  $\mathcal B = \mathcal H$, revisited the Johansen I(1) condition and provided its geometric reformulation given by a certain nonorthogonal direct sum of $\mathcal H^p$; see Remark \ref{rembss}. Inspired by their direct sum condition that can be applied for $\mathcal H$ of an arbitrary dimension, we propose the following condition.	\\[0.7em]
	\noindent \textbf{I(1) condition}: $\mathcal B^p = \ran \PPhi(1) \oplus \ker \PPhi(1)$. \\[0.7em] 
	Some remarks on the I(1) condition are given in order. 
	\begin{remark} \label{rembss0}
		The I(1) condition is given as the direct sum of $\mathcal B^p$ by two fixed subspaces  $\ran \widetilde{\Phi}(1)$ and $\ker \widetilde{\Phi}(1)$, and this specific direct sum condition will be shown to be necessary and sufficient for the existence of I(1) solutions. In the case where our I(1) condition holds, it is worth noting that a unique projection whose range is $\ker \PPhi(1)$ and kernel is $\ran \PPhi(1)$ is well defined \citep[Theorem 3.2.11]{megginson1998}.   
	\end{remark}
	
	\begin{remark} \label{rembss}
		In the case where $\mathcal B=\mathcal H$ of an arbitrary dimension and $\Phi_1,\ldots,\Phi_p$ are compact operators, \cite{BSS2017} showed that the nonorthogonal direct sum  $\mathcal H^p = \ran \PPhi(1) \oplus \ker  \PPhi(1)$ is a sufficient condition for the AR($p$) law of motion \eqref{arlaw} to admit I(1) solutions; however, its  necessity was not discussed in their paper. They showed that their condition becomes  equivalent to the Johansen I(1) condition if $\mathcal H = \mathbb{R}^n$ or $\mathbb{C}^n$. The reader is referred to the results given in Section 4 (and the proofs of those) of their paper.   
	\end{remark}

	Our first result in this section not only shows that the I(1) condition is a necessary and sufficient condition for $\PPhi(z)^{-1}$ to have a simple pole at $z=1$ but also characterizes the principal part of its Laurent series.

	\begin{proposition}\label{cormain}
		Suppose that {Assumption} \ref{assulinear1} holds. The following conditions are equivalent. 
		\begin{enumerate}[\normalfont(i)]\setlength\itemsep{0em}
			\item\label{cormain1} $\PPhi(z)^{-1}$ has a simple pole at $z=1$. 	
			\item\label{cormain2} $\PP$ is the projection onto $\ker \PPhi(1)$ along $\ran  \PPhi(1)$. 
			\item\label{cormain3} The I(1) condition holds.	
		\end{enumerate}
		Under any of these  conditions, the following holds: for some $\eta > 0$, 
		\begin{equation}
			(1-z)\PPhi(z)^{-1} = \PP + (1-z)H(z), \quad z \in D_{1+\eta}, \label{eqlaurentform} 
		\end{equation}
		where $H(z)$ denotes the holomorphic part of the Laurent series of $\PPhi(z)^{-1}$ around $z=1$. Moreover, each  Maclaurin series of $(1-z)\PPhi(z)^{-1}$ and $H(z)$ is convergent on $D_{1+\eta}$.
	\end{proposition}

	Examples of the use of Proposition \ref{cormain} for verifying that $\widetilde{\Phi}(z)^{-1}$ has a simple pole at $z=1$  will be given later in this section. Proposition \ref{cormain} extends the results given by \cite{BSS2017}, which are briefly reviewed in Remark \ref{rembss}, in the sense that it provides a necessary and sufficient condition for the existence of I(1) solutions without requiring either of a Hilbert space structure or compactness of $\Phi_1,\ldots,\Phi_p$. In addition, combined with Propositions \ref{propa1}, the local behavior of $\PPhi(z)^{-1}$ around $z=1$ given by \eqref{eqlaurentform} provides a characterization of solutions to the AR($p$) law of motion \eqref{arlaw}, which  leads to our first version of  the Granger-Johansen representation theorem for I(1) AR processes given below.   
	
	\begin{proposition}\label{grti1}
		Suppose that  Assumption \ref{assulinear1} holds. Under the I(1) condition, a sequence $\{X_t\}_{t \geq -p+1}$ satisfying \eqref{arlaw} allows the Johansen I(1) representation \eqref{bndecom} with 
		\begin{equation}
			\Upsilon_{-1} = \Pi_p\PP \Pi_p^\ast,  \quad\quad  	\nu_t=\Pi_pH(L)\Pi_p^\ast\varepsilon_t = \sum_{j=0}^\infty \Pi_p \PPhi_1^j (\idbp-\PP) \Pi_p^\ast \varepsilon_{t-j}, \label{bndecomi1}
		\end{equation}
		where  $\PP$ and $H(z)$ are given in Proposition \ref{cormain}, and $\Pi_p$ and $\Pi_p^\ast$ are given in \eqref{coorproj}.
		Moreover, the AR($p$) law of motion \eqref{arlaw} does not allow I(1) solutions if the I(1) condition is not satisfied.
	\end{proposition}
	
	Proposition \ref{grti1} shows that, under our I(1) condition, solutions to the AR($p$) law of motion \eqref{arlaw} can be represented as \eqref{bndecom} similar to the Beveridge-Nelson decomposition \eqref{eqi1bn} of an I(1) cointegrated linear process.  For such a solution $\{X_t\}_{t\geq0}$, we may deduce from our discussion in Section \ref{ssinteg} that $\{f(X_t)\}_{t\geq 0}$ can be made stationary under a suitable initial condition if and only if $f \in \Ann(\Pi_p\PP \Pi_p^\ast)$. Some more remarks on the results given by Proposition \ref{grti1} are in order.


	\begin{remark} \label{remcointeg}
		Proposition \ref{grti1} may be viewed as an extension of Theorem 4.1 of \cite{BSS2017}, which provides a version of the Granger-Johansen representation theorem in a Hilbert space setting resorting to the companion form representation of a given AR($p$) law of motion. 
	\end{remark} 	
	
	\begin{remark}\label{schrem3}
		In the case $\dim(\mathcal B)=\infty$, neither of the attractor nor the cointegrating space associated with \eqref{arlaw} is compelled to be  finite dimensional in our representation theorem; this is in contrast to that  one of those subspaces is necessarily finite dimensional in the existing theorems developed in a Hilbert space setting; see Section \ref{sec:relation}.
		As a simple illustration, let $p=1$ and $\Phi_1$ be an arbitrary projection. In this case, we may deduce from Propositions \ref{propbasic}, \ref{cormain}, and \ref{grti1} that the dimension of the attractor space (resp.\ the cointegrating space) associated with  \eqref{arlaw} is equal to $\dim(\ran \Phi_1)$ (resp.\ $\dim(\ker \Phi_1)$). Since $\Phi_1$ is an arbitrary projection, all the following cases are possible: (i) $\dim (\ran \Phi_1) < \infty$ and $\dim (\ker \Phi_1) = \infty$, (ii) $\dim (\ran \Phi_1) = \infty$ and $\dim (\ker \Phi_1) <  \infty$, and (iii) $\dim (\ran \Phi_1) = \infty$ and $\dim (\ker \Phi_1) = \infty$. 
	\end{remark}
	
	In this section, the I(1) condition and solutions to the AR($p$) law of motion \eqref{arlaw} are characterized in terms of the behavior of $\widetilde{\Phi}(z)$ around $z=1$.  For this reason, our results do not clearly reveal how the proposed I(1) condition and the cointegrating behavior of I(1) solutions are related to the structure of the original AR polynomial $\Phi(z)$.  This is a natural consequence resulting from that we resort to the companion form representation of  \eqref{arlaw} to obtain our main results given in this section. In the next section, we will discuss on how these results can be recast in terms of the behavior of the AR polynomial $\Phi(z)$.  By doing so, we obtain a more detailed characterization of I(1) solutions and find the connection between our representation results and  those developed in a Hilbert/Euclidean space setting. 	We close this section with a few examples illustrating the use of Proposition \ref{cormain} for verifying that $\widetilde{\Phi}(z)^{-1}$ has a simple pole at $z=1$. 
	\begin{example} \label{exi1}
		In the example given in Section \ref{example}, $\mathcal B=C[0,1]$, $\widetilde{\Phi}(z) = I - z \Phi_1$, and $\Phi_1$ is defined by $\Phi_1(x)(u) = x(1)$ for $u \in [0,1]$. Note that $\PPhi(1)\neq 0$ and $\widetilde{\Phi}(1)x(1) = 0$ for any arbitrary  $x \in C[0,1]$. This implies that every $y \in \ran \widetilde{\Phi}(1)$ must satisfy $y(1) = 0$, from which we find that $\widetilde{\Phi}(1)$ is not invertible. Moreover, it can be verified that  $\widetilde{\Phi}(z)^{-1}$ is well defined for any $z \neq 1$ and given by
		\begin{equation*}
			\widetilde{\Phi}(z)^{-1} x(u) = x(u) + {zx(1)}/{(1-z)}, \quad  u \in [0,1]. 
		\end{equation*}
		Thus, Assumption \ref{assulinear1}-\ref{assulinear1a} is satisfied. Let $\mathrm{C}_0$ be the set of continuous functions $x\in \mathcal B $ satisfying $x(1)= 0$, and let  $\mathrm{C}_1$ be the set of constant functions. Then, $\mathcal B  = \mathrm{C}_0 \oplus \mathrm{C}_1$ may be easily deduced. We will show that $\ran	\widetilde{\Phi}(1) = \mathrm{C}_0$ and $\ker\widetilde{\Phi}(1) = \mathrm{C}_1$; this implies that both Assumption \ref{assulinear1}-\ref{assulinear1b} and the I(1) condition are satisfied, and thus $\widetilde{\Phi}(z)^{-1}$ has a simple pole at $z=1$.  It was already shown that any $x \in \ran\widetilde{\Phi}(1)$ satisfies $x(1)=0$, hence $\ran	\widetilde{\Phi}(1) \subset \mathrm{C}_0$. Moreover, for any $x \in \mathrm{C}_0$ we have $x=\widetilde{\Phi}(1)x \in \ran \widetilde{\Phi}(1)$, which implies that $\ran	\widetilde{\Phi}(1) \supset \mathrm{C}_0$.  Thus, $\ran\widetilde{\Phi}(1) = \mathrm{C}_0$ holds.	To show  $\ker\widetilde{\Phi}(1) = \mathrm{C}_1$, we first observe that 
		\begin{equation*}
			x \in \ker \widetilde{\Phi}(1) \quad\Rightarrow \quad x(u) =  x(1), \quad u \in [0,1],
		\end{equation*}
		and hence $\ker \widetilde{\Phi}(1)\subset \mathrm{C}_1$ holds. The reverse inclusion ($\mathrm{C}_1\subset \ker \widetilde{\Phi}(1)$) immediately follows from the fact that $\widetilde{\Phi}(1)x = 0$ for any constant function $x \in C[0,1]$. We thus also find that  $\ker \widetilde{\Phi}(1) = \mathrm{C}_1$. 
	\end{example}

	\begin{example}\label{ex4}
		Suppose that \(\mathcal B = C[-1,1] \), $\widetilde{\Phi}(z) = \idb-z\Phi_1$, and \(\Phi_1\) is defined by
		\[\Phi_1 x(u) = x(u)/2  \,+\, x(-u)/2,\quad x \in \mathcal B,\quad u \in [-1,1]. \] 
		One can easily show that $\sigma(\widetilde{\Phi}) = \{1\}$, and $\Phi_1$ (resp.\ $\idb-\Phi_1$) is the projection onto the space of even (resp.\ odd) functions along the space of odd (resp.\ even) functions. Moreover, in this case we have \(\mathcal B = \ran \widetilde{\Phi}(1) \oplus \ker\widetilde{\Phi}(1)\) and $\ker \widetilde{\Phi}(1) = \ran \Phi_1$. Therefore, Assumption \ref{assulinear1} and the I(1) condition are satisfied, from which it is concluded that $\widetilde{\Phi}(z)^{-1}$ has a simple pole at $z=1$. In fact, we reach the same conclusion in cases where any arbitrary projection replaces $\Phi_1$ in the above.   
	\end{example}
	
	\begin{example} \label{examplei1i2}
		Let $\mathbf{c_0}$ be the space of complex sequences converging to zero equipped with the norm $\|a\| = \sup_i |a_i|$ for $a = (a_1,a_2,\ldots) \in \mathbf{c_0}$.  The space $\mathbf{c_0}$ may be viewed as a natural generalization of a finite dimensional vector space equipped with the supremum norm, and also turns out to be a separable Banach space \citep[Examples 1.2.13 and 1.12.6]{megginson1998}.  Let $\Phi_1$ be defined by
		\begin{equation}
			\Phi_1(a_1,a_2,a_3,a_4\ldots) = (a_1,a_1+a_2,\lambda a_3,\lambda^2 a_4,\ldots), \label{eqopex}
		\end{equation}
		where $\lambda \in (0,1)$.
		Then, $\widetilde{\Phi}(z) =  \idb-z\Phi_1$ maps $a=(a_1,a_2,\ldots)$ to 
		\begin{equation*}
			(\idb-z\Phi_1)a  = ( (1-z)a_1, (1-z)a_2 - za_1, (1-z\lambda)a_3, (1-z\lambda^2)a_4, \ldots ).
		\end{equation*}
		It may be  deduced that $\widetilde{\Phi}(z)$ is injective on \(\mathbf{c_0}\) for any $z \in  D_{1+\eta} \setminus \{1\}$. Furthermore, for any  sequence  \(b=(b_1,b_2,b_3\ldots) \in \mathbf{c_0}\), we can find a sequence $a=(a_1,a_2,a_3\ldots) \in \mathbf{c_0}$ satisfying \(\widetilde{\Phi}(z)a = b\) by setting
		\begin{equation*}
			a_1 = {b_1}/(1-z), \quad a_2 =b_2/(1-z)+ zb_1/(1-z)^2, \quad a_j = {b_j}/{(1-z\lambda^{j-2})}, \quad j \geq 3.
		\end{equation*}
		This shows that $\widetilde{\Phi}(z)$ is also a surjection for $z \in  D_{1+\eta} \setminus \{1\}$. Therefore, we have shown that $\widetilde{\Phi}(z)$ is invertible on $D_{1+\eta} \setminus \{1\}$. Note also that 	$\ran \widetilde{\Phi}(1)$ and $\ker  \widetilde{\Phi}(1)$ are given as follows,
		\begin{equation}
			\ran \widetilde{\Phi}(1)= \{(0,b_1,b_2,\ldots) : \,  b_j \in \mathbb{C}, \,\lim_{j\to \infty} b_j = 0  \},  \quad \quad \ker  \widetilde{\Phi}(1)= \{(0,b_1,0,0,\ldots) : \, b_1 \in \mathbb{C}\}.  \label{eqran1} 
		\end{equation}
		The above subspaces can be complemented (see Example \ref{exi2000}), hence Assumption \ref{assulinear1} is satisfied. However, \eqref{eqran1} clearly shows that $\mathcal B \neq 	\ran \widetilde{\Phi}(1) \oplus 	\ker \widetilde{\Phi}(1)$. We thus conclude that $\widetilde{\Phi}(z)^{-1}$ does not have a simple pole at $z=1$; it will be shown in Section \ref{sec:repi2} that $\PPhi(z)^{-1}$ has a pole of order 2 at $z=1$ .
	\end{example}

	\subsubsection{Further characterization of I(1) solutions}	\label{sec:repi1noncompanion}
	If $\mathcal B = \mathbb{R}^n$ or $\mathbb{C}^n$, the Johansen I(1) condition given by \eqref{eqjohansen} is necessary and sufficient for $\Phi(z)^{-1}$ to have a simple pole at $z=1$. If we let $\PR_{[\ran \Phi(1)]^\perp}$ (resp.\  $\PR_{[\ker \Phi(1)]^\perp}$) denote the orthogonal projection onto $[\ran \Phi(1)]^\perp$ (resp.\ $[\ker \Phi(1)]^\perp$), then the Johansen I(1) condition can be alternatively understood as invertibility of $\PR_{[\ran \Phi(1)]^\perp}\Phi^{(1)}(1)(I-\PR_{[\ker \Phi(1)]^\perp})$ as a map from $\ker \Phi(1)$ to $[\ran \Phi(1)]^\perp$.  A natural generalization of this condition to a general Banach space setting may be the following: for some $[\ran \Phi(1)]_\complement$ and $[\ker \Phi(1)]_\complement$ satisfying \eqref{direct1}, 
	\begin{equation} \label{eqjohan1}
		\Lambda_{1,\complement}\coloneqq \PR_{[\ran \Phi(1)]_\complement}\Phi^{(1)}(1)(I-\PR_{[\ker \Phi(1)]_\complement}):\ker \Phi(1) \mapsto [\ran \Phi(1)]_{\complement} \text{ is invertible,}
	\end{equation}
	where $\PR_{[\ran \Phi(1)]_\complement}$ and $\PR_{[\ker \Phi(1)]_\complement}$ are the projections satisfying \begin{equation}\label{projeq}
		\ran \PR_{V_\complement} = V_\complement,\,\,\,\,  \ker\PR_{V_\complement}=V, \,\,\,\,\,\,\, (V,V_\complement)=(\ran \Phi(1),[\ran \Phi(1)]_\complement) \, \text{ or } \, (\ker \Phi(1),[\ker \Phi(1)]_\complement).
	\end{equation}
	These projections are  uniquely defined linear operators for each choice of  ${[\ran \Phi(1)]_\complement}$ and ${[\ker \Phi(1)]_\complement}$; see \citet[Theorem 3.2.11]{megginson1998}.
	In the case where $\mathcal B=\mathbb{R}^n$ or $\mathbb{C}^n$, $[\ran \Phi(1)]_\complement=[\ran \Phi(1)]^\perp$, and $[\ker \Phi(1)]_\complement=[\ker \Phi(1)]^\perp$, the condition \eqref{eqjohan1} becomes equivalent to the Johansen I(1) condition, under which we know from \citet[Theorem 4.2]{Johansen1996} that $\Phi(z)^{-1}$ has a simple pole at $z=1$. 	Proposition \ref{propjohanseni1} given below  shows that this result can be extended to  a more general Banach space setting, and also provides a characterization of $\Upsilon_{-1}$ in terms of the operators defined above. 
	
	\begin{proposition}\label{propjohanseni1} 
		Suppose that Assumption \ref{assulinear1} holds. Then the following conditions are equivalent.
		\begin{enumerate}[\normalfont(i)]
			\item The I(1) condition holds.
			\item  $\Lambda_{1,\complement}:\ker \Phi(1) \mapsto [\ran \Phi(1)]_{\complement}$ is invertible for some choice of $[\ran \Phi(1)]_\complement$ and $[\ker \Phi(1)]_\complement$.
			\item $\Lambda_{1,\complement}:\ker \Phi(1) \mapsto [\ran \Phi(1)]_{\complement}$ is invertible for every possible choice of $[\ran \Phi(1)]_\complement$ and $[\ker \Phi(1)]_\complement$.
		\end{enumerate}
		Let any of the above conditions hold. Then for any  choice of $[\ran \Phi(1)]_\complement$ and $[\ker \Phi(1)]_\complement$, a sequence \(\{X_t\}_{t \geq -p+1}\) satisfying \eqref{arlaw} allows the Johansen I(1) representation \eqref{bndecom} with $\Upsilon_{-1}$ satisfying
		\begin{align} 
			&\PR_{[\ker \Phi(1)]_\complement} 	\Upsilon_{-1} =	\Upsilon_{-1} (I-\PR_{[\ran \Phi(1)]_\complement}) = 0, \quad 	\Upsilon_{-1}:[\ran \Phi(1)]_\complement \mapsto \ker \Phi(1) = \Lambda_{1,\complement}^{-1}.  \label{n1formula}
		\end{align}
	\end{proposition}
	It is interesting that a natural generalization of the Johansen I(1) condition is equivalent to our previous necessary and sufficient condition for the existence of I(1) solutions developed in a general Banach space setting. This shows that our operator-theoretic approach  is  in fact closely related to the conventional Johansen's approach. We know from the equivalence between the three conditions given in Proposition \ref{propjohanseni1} that $[\ran\Phi(1)]_\complement$  and $[\ker\Phi(1)]_\complement$ can be arbitrarily chosen among possible candidates; this, of course, implies that $[\ran \Phi(1)]_\complement$ (resp.\ $[\ker \Phi(1)]_\complement$) can always be fixed to $[\ran \Phi(1)]^\perp$ (resp.\ $[\ker \Phi(1)]^\perp$) in a Hilbert/Euclidean space setting without loss of generality. Moreover,  \eqref{n1formula} describes how the operator $\Upsilon_{-1}$ acts as a map from $\ran \Phi(1)\oplus [\ran \Phi(1)]_\complement$ to $\ker \Phi(1)\oplus [\ker \Phi(1)]_\complement$; this in fact provides a full characterization of $\Upsilon_{-1}$, given that $\mathcal B$ allows the direct sums given in \eqref{direct1}.

	In Section \ref{sec:repi1companion}, the I(1) condition and the cointegrating behavior of I(1) solutions  are characterized in terms of some operators associated with $\PPhi(z)$ given in the companion form  \eqref{a1companion0}; however, one may be interested in characterizing those using operators associated with the original AR polynomial $\Phi(z)$. The conditions given by (ii) and (iii) in Proposition \ref{propjohanseni1} are already given in such a manner, and those are equivalent reformulations of our I(1) condition. Moreover, our  characterization of $\Upsilon_{-1}$, given by \eqref{n1formula}, helps us characterize the cointegrating behavior in a desired way; see Remarks \ref{remcointeg01} and  \ref{johansencointeg} below. Some more remarks are given for comparison between our I(1) representation result given by Proposition \ref{propjohanseni1} and the existing ones developed in a Hilbert/Euclidean space setting.

	\begin{remark} \label{remcointeg01}
		From \eqref{n1formula}, we know that the attractor space of I(1) solutions is given by $\ran \Upsilon_{-1} = \ker \Phi(1)$, which is complemented by another subspace $[\ker \Phi(1)]_\complement$ under Assumption \ref{assulinear1}-\ref{assulinear1b}. We then deduce from Proposition \ref{propbasic} that a cointegrating functional $f$ satisfies $f =g\, \PR_{[\ker \Phi(1)]_\complement}$, where $g \in \mathcal B'$. This result, of course, holds for any arbitrary choice of $[\ker \Phi(1)]_\complement$; a different choice of $[\ker \Phi(1)]_\complement$ only affects the definition of $\PR_{[\ker \Phi(1)]_\complement}$ without any further changes.
	\end{remark}   
	\begin{remark}\label{johansencointeg}
		Using the results given in Proposition \ref{propjohanseni1}, we may obtain a stronger characterization of the cointegrating behavior of I(1) solutions than that given  in Section  \ref{sec:repi1companion}. From the expression of $\Upsilon_{-1}$ given in \eqref{n1formula}, we find that a nonzero element $f \in \mathcal B'$ satisfies
		\begin{equation}
			\text{\text{$\{f(X_t)\}_{t \geq 0}$ is I($0$)}\, if and only if \,$f \in \Ann(\ker \Phi(1))$;} \label{cointegchac1}
		\end{equation}
		see Appendix \ref{appenBrr} for our proof of \eqref{cointegchac1}. The above characterization not only shows that the cointegrating space is given by $\Ann(\ker \Phi(1))$, but also establishes I(0)-ness of $\{f(X_t)\}_{t\geq0}$ for any cointegrating functional $f$. In the case $\mathcal B = \mathcal H$, any $f \in  \mathcal H'$ is identified as the map $\langle \cdot,v \rangle$ for a unique element $v\in \mathcal H$, so $f \in \Ann(\ker \Phi(1))$ is given by the map $\langle \cdot,v \rangle$ for some $v \in [\ker \Phi(1)]^\perp$. (To see why, note that $\langle x,v \rangle = 0$ for all $x \in \ker \Phi(1)$ if and only if $v \in [\ker \Phi(1)]^\perp$.) Thus, \eqref{cointegchac1} reduces to the following: for any $v \in \mathcal H \setminus \{0\}$,  
		\begin{equation}\label{cointegchac2}
			\text{\text{$\{\langle X_t,v \rangle\}_{t \geq 0}$ is I($0$)}\, if and only if \,$v \in [\ker \Phi(1)]^\perp$.}
		\end{equation}
		If $\Phi(z)$ is a Fredholm operator satisfying \eqref{eqfred} and, as a result, $\ker \Phi(1)$ is finite dimensional, \eqref{cointegchac2} reduces to the characterization of the cointegrating behavior of I(1) solutions provided by \cite{BS2018} and \cite{Franchi2017b}. However, our results given by \eqref{cointegchac1} and \eqref{cointegchac2} do not require the Fredholm assumption and, as a consequence,  $\ker \Phi(1)$ is not necessarily finite dimensional; noting that $\ker \PPhi(1) = \ker \Phi(1)$  in the case $p=1$, we see this from Remark \ref{schrem3} and/or Example \ref{ex4}.  
	\end{remark}
	
	\begin{remark} \label{bsi1condition}
		If $\mathcal B=\mathcal H$, we know from Proposition \ref{propjohanseni1} that the following is a necessary and sufficient condition for the AR($p$) law of motion \eqref{arlaw} to admit I(1) solutions: 
		\begin{equation}\label{n1formula10}
			\Lambda_{1}=\PR_{[\ran \Phi(1)]^\perp}\Phi^{(1)}(1)(I-\PR_{[\ker \Phi(1)]^\perp}):\ker \Phi(1) \mapsto [\ran \Phi(1)]^\perp \text{ is invertible.}
		\end{equation} 
		If the Fredholm assumption \eqref{eqfred} holds and thus $\ker \Phi(1)$ is finite dimensional, the above invertibility condition reduces to one of the necessary and sufficient conditions for the AR($p$) law of motion \eqref{arlaw} to admit I(1) solutions given by \citet[Theorem 3.2]{BS2018}. Their invertibility condition, a special case of \eqref{n1formula10}, was given together with an equivalent nonorthogonal direct sum condition given by  $\mathcal H = \ran \Phi(1) \oplus \Phi^{(1)}(1)\ker \Phi(1)$. In the same setting, \cite{Franchi2017b} provided a necessary and sufficient condition for the existence of I($d$) solutions for any $d \geq 1$,  which is characterized as an orthogonal direct sum of $\mathcal H$ and turns out to be useful to identify the cointegrating behavior of I($d$) solutions; for $d=1$, their condition is equivalent to the nonorthogonal direct sum condition given by \cite{BS2018} \citep[Proposition 4.6]{Franchi2017b}.  Under any of these equivalent conditions, solutions to  \eqref{arlaw} allows the Johansen I(1) representation  \eqref{bndecom} with $\Upsilon_{-1}$ of finite rank. (To see why, note that \eqref{n1formula} implies $\ran \Upsilon_{-1} = \ker \Phi(1)$ and $\ker \Phi(1)$ is necessarily finite dimensional under the Fredholm assumption.) This implies that the random walk component given in \eqref{bndecom} reduces to a finite dimensional unit root process. On the other hand, if $\Phi(z)$ satisfies Assumption \ref{assulinear1} and $\Phi(1)$ is a compact operator, then $\Phi(1)$ turns out to be a finite rank operator \citep[Lemma 1]{chang2016}. In this case, \eqref{n1formula10} reduces to the condition given by \citet[Theorem 2]{chang2016} as a sufficient condition for the existence of I(1) solutions, and the random walk component of such a solution takes values in $\ran \Upsilon_{-1} = \ker \Phi(1)$, which is infinite dimensional unless $\mathcal H$ is finite dimensional. In our results for the case $\mathcal B=\mathcal H$, $\ran \Upsilon_{-1} = \ker \Phi(1)$ is not required to be either finite or infinite dimensional; we see this from Remark \ref{schrem3}. Thus, our I(1) representation result given by Proposition \ref{propjohanseni1} complements the earlier results developed in a Hilbert space setting.
	\end{remark}
	\begin{remark} \label{johansencointeg2}
		Suppose that $\mathcal B = \mathbb{R}^n$ or $\mathbb{C}^n$,  $[\ran \Phi(1)]_\complement=[\ran \Phi(1)]^\perp$, and $[\ker \Phi(1)]_\complement=[\ker \Phi(1)]^\perp$. In this case, $\Upsilon_{-1}$ may be understood as an $n\times n$ matrix. Using the notation introduced for  \eqref{eqjohansen}, the results given by  \eqref{n1formula} can be written as $\beta(\beta^\intercal \beta)^{-1}\beta^\intercal\Upsilon_{-1}= \Upsilon_{-1} \alpha(\alpha^\intercal \alpha)^{-1}\alpha^\intercal = 0$ and $\sbot{\beta}^\intercal\Upsilon_{-1} \sbot{\alpha} = (\sbot{\alpha}^\intercal \Phi^{(1)}(1) \sbot{\beta})^{-1}$. We then may deduce that $\Upsilon_{-1}$ can be written as
		\begin{equation*}
			\Upsilon_{-1} =\sbot{\beta} \left(\sbot{\alpha}^\intercal \Phi^{(1)}(1)\sbot{\beta}\right)^{-1} \sbot{\alpha}^\intercal,
		\end{equation*}
		which is equivalent to the expression of $\Upsilon_{-1}$ given by \citet[Theorem 4.2]{Johansen1996}.
	\end{remark}

	\subsection{Representation of I(2) autoregressive processes}	\label{sec:repi2}
	In this section, we suppose that the I(1) condition fails and develop our representation theory for I(2) AR processes. 	
	The failure of the I(1) condition means  either of (i) $\mathcal B^p \neq \ran \PPhi(1) \oplus \ker \PPhi(1)$ or (ii) noninvertibility of the map $\Lambda_{1,\complement}$ given in \eqref{eqjohan1} for some  $[\ran \Phi(1)]_\complement$ and $[\ker \Phi(1)]_\complement$. 	Due to  Proposition \ref{propjohanseni1},  we know that $[\ran \Phi(1)]_\complement$ and $[\ker \Phi(1)]_\complement$ may be arbitrarily chosen among possible candidates without loss of generality. Especially, $[\ran \Phi(1)]_\complement=[\ran \Phi(1)]^\perp$ and $[\ker \Phi(1)]_\complement=[\ker \Phi(1)]^\perp$ can be assumed whenever $\mathcal B$ is a Hilbert space; this helps us compare the results to be developed with the existing I(2) results given in a Hilbert space setting.  	We define  $\PR_{[\ran \Phi(1)]_\complement}$ and $\PR_{[\ker \Phi(1)]_\complement}$ as in \eqref{projeq} and let 
	\begin{equation*}
		\mathsf{R} \coloneqq \PR_{[\ran \Phi(1)]_\complement}\Phi^{(1)}(1)\ker \Phi(1), \,\, \quad \mathsf{K} \coloneqq \{x \in \ker \Phi(1) : \Phi^{(1)}(1)x \in \ran \Phi(1)\}.
	\end{equation*}
	The conditions that we require for our I(2) representation results are summarized below in Assumption \ref{assulinear2}. 
	\begin{assumption}\label{assulinear2} \hspace{0.1cm}
		\begin{enumerate}[\normalfont(a)]
			\item\label{assulinear2a} 	 Assumption \ref{assulinear1} holds and the I(1) condition fails for some $[\ran \Phi(1)]_\complement$ and $[\ker \Phi(1)]_\complement$; if $\mathcal B$ is a Hilbert space, $[\ran \Phi(1)]_\complement=[\ran \Phi(1)]^\perp$ and $[\ker \Phi(1)]_\complement=[\ker \Phi(1)]^\perp$.
			\item\label{assulinear2b}  $\mathsf{R}$ (resp.\ $\mathsf{K}$) can be complemented in $[\ran \Phi(1)]_\complement$ (resp.\ $\ker \Phi(1)$), i.e., for some $\mathsf{R}_\complement\subset \mathcal B$ and  $\mathsf{K}_\complement\subset \mathcal B$, 	
			\begin{equation}[\ran \Phi(1)]_\complement = \mathsf{R}  \oplus \mathsf{R}_\complement, \quad \quad  \ker \Phi(1) = \mathsf{K}  \oplus \mathsf{K}_\complement. \label{direct2}\end{equation}
		\end{enumerate} 
	\end{assumption}
	The complementary subspaces $\mathsf{R}_\complement$ and $\mathsf{K}_\complement$ appearing in Assumption \ref{assulinear2} do not uniquely exist in general; however, our subsequent results only require their existence. The direct sum conditions  required by Assumption \ref{assulinear2}, i.e., \eqref{direct1} and \eqref{direct2}, may not be restrictive in general and do not invalidate that our I(2) results to be developed can complement the earlier results developed in a general Hilbert space setting; the regularity condition imposed on $\Phi(z)$ for the existing  I(2) results strictly implies the required direct sums and, moreover, places some certain  restrictions on the dimensions of the complementary subspaces appearing in Assumption \ref{assulinear2}, see Remark \ref{remfredi2}.  
	

	\begin{remark}\label{remfredi2} 
	Fredholmness of $\Phi(z)$ is an essential assumption in the existing I(2) representation theorems developed in the case $\mathcal B=\mathcal H$. In such a setting, the direct sums given by \eqref{direct1} hold for some $[\ran \Phi(1)]_\complement$ and $[\ker \Phi(1)]_\complement$ (Remark \ref{remi1a}), and both $[\ran \Phi(1)]_\complement$ and $\ker \Phi(1)$ are finite dimensional. Since every finite dimensional space can be complemented (Remark \ref{rem33}), the direct sums given by \eqref{direct2} hold for some $\mathsf{R}_\complement$ and $\mathsf{K}_\complement$.  This is also true in a more general situation where $\mathcal B$ is not necessarily a Hilbert space and $\Phi(z)$ is a Fredholm operator acting on $\mathcal B$; this can be seen from Remark \ref{rem33}. Note that the Fredholm assumption requires that $[\ran \Phi(1)]_\complement$, $\ker \Phi(1)$, $\mathsf{R}_\complement$ and $\mathsf{K}_\complement$ are all finite dimensional in \eqref{direct1} and \eqref{direct2}, which means that Fredholmness of $\Phi(z)$ naturally leads us to stronger conditions than the direct sums required by Assumption \ref{assulinear2}. 
	\end{remark}
	
	As in Section \ref{sec:repi1}, we first develop our representation theory for I(2) AR processes resorting to the companion form representation \eqref{a1companion0}. We then discuss on how the results obtained via the companion form can be recast in terms of the behavior of the AR polynomial $\Phi(z)$. 	 In the subsequent discussions, we will need a notion of a generalized inverse operator in our Banach space setting.  Appendix \ref{sginverse} introduces a suitable one extending the notion of the Moore-Penrose inverse, employed in \cite{BS2018} and \cite{Franchi2017b} for their I(2) representation results developed in a Hilbert space setting.	
	
	\subsubsection{Representation via the companion form}\label{sec:repi2companion}
	As in Section \ref{sec:repi1companion}, we first resort to linearization of $\Phi:\mathbb{C}\mapsto \mathcal B$, hence consider $\PPhi:\mathbb{C}\mapsto \mathcal B^p$ given in \eqref{a1companion0} and \eqref{a1companion}. Under Assumption \ref{assulinear2}, we have the Laurent series of $\PPhi(z)^{-1}$ near $z=1$ as in \eqref{laurent0}, and also define $\PP$ as in \eqref{defPP}. Before stating our main results, we collect some preliminary results and fix notation. 
	
	Under  Assumption \ref{assulinear2}, $\ran \PPhi(1)$ and $\ker \PPhi(1)$ can be complemented (Lemma \ref{lemi2}-\ref{lemi21}), i.e., for some $[\ran \PPhi(1)]_\complement \subset \mathcal B^p$ and $[\ker \PPhi(1)]_\complement \subset \mathcal B^p$,
	\begin{equation}
		\mathcal B^p = \ran \PPhi(1) \oplus [\ran \PPhi(1)]_{\complement}  = \ker \PPhi(1) \oplus [\ker \PPhi(1)]_{\complement}.  \label{direct1a} 	\end{equation}
	The complementary subspaces $[\ran \PPhi(1)]_\complement$ and $[\ker \PPhi(1)]_\complement$ given in \eqref{direct1a} are not uniquely determined in general, but the subsequent results only require their existence and do not depend on a specific choice of them. Under the direct sums given by \eqref{direct1a} and for any choice of $[\ran \PPhi(1)]_\complement$ and $[\ker \PPhi(1)]_\complement$, we may define the  projections $\PR_{[\ran \PPhi(1)]_\complement}$ and $\PR_{[\ker \PPhi(1)]_\complement}$ satisfying 
	\begin{equation}
		\ran \PR_{V_\complement} = V_\complement,\,\,\,\,  \ker\PR_{V_\complement}=V, \,\,\,\,\,\,\, (V,V_\complement)=(\ran \PPhi(1),[\ran \PPhi(1)]_\complement) \, \text{ or } \, (\ker \PPhi(1),[\ker \PPhi(1)]_\complement), \label{direct1aproj} 
	\end{equation}
	and the generalized inverse $\PPhi(1)^g$ of $\PPhi(1)$ satisfying 
	\begin{equation}
		\PPhi(1) \PPhi(1)^g= \idbp-\PR_{[\ran \PPhi(1)]_\complement}, \quad \,\, \PPhi(1)^g\PPhi(1) = \PR_{[\ker \PPhi(1)]_\complement}.\label{direct1apro2} 
	\end{equation}
	The operators $\PR_{[\ran \PPhi(1)]_\complement}$, $\PR_{[\ker \PPhi(1)]_\complement}$, and $\PPhi(1)^g$ are uniquely defined for each choice of $[\ran \PPhi(1)]_\complement$ and $[\ker \PPhi(1)]_\complement$; see \citet[Theorem 3.2.11]{megginson1998} and Appendix \ref{sginverse}. If $\mathcal B=\mathcal H$, $[\ran \PPhi(1)]_\complement=[\ran \PPhi(1)]^\perp$, and $[\ker \PPhi(1)]_\complement=[\ker \PPhi(1)]^\perp$, then $\PR_{[\ran \PPhi(1)]_\complement}$ (resp.\ $\PR_{[\ker \PPhi(1)]_\complement}$) becomes the orthogonal projection onto ${[\ran \PPhi(1)]^\perp}$ (resp.\ ${[\ker \PPhi(1)]^\perp}$) and $\PPhi(1)^g$ becomes the Moore-Penrose inverse of $\PPhi(1)$ (see Appendix \ref{sginverse}).
	For notational convenience, we let \begin{equation}
		\mathcal K = \ran \PPhi(1) \cap \ker \PPhi(1). \label{eqwkwk}
	\end{equation} 	
	A crucial preliminary result is that $\mathcal K \neq \{0\}$ is required for $\PPhi(z)^{-1}$ to have a pole of order 2 at $z=1$ (Lemma \ref{lemi2}-\ref{lemi23}); based on this result and the notation defined above, we propose our I(2) condition as follows.\\[0.7em]
	\noindent \textbf{I(2) condition}:  $\mathcal K \neq \{0\}$ and $\mathcal B^p =(\ran \PPhi(1) + \ker \PPhi(1)) \oplus \PPhi(1)^g \mathcal K$ for some $[\ran \PPhi(1)]_\complement$ and $[\ker \PPhi(1)]_\complement$. \\[0.7em] 	
	\noindent Some remarks on the I(2) condition are given in order. 
	\begin{remark}	\label{remi2conditiona}   
		The I(2) condition is given as the direct sum of $\mathcal B^p$ by two subspaces $\ran\PPhi(1)+\ker \PPhi(1)$ and $\PPhi(1)^g\mathcal K$, where the latter depends on the choice of $[\ran \PPhi(1)]_\complement$ and $[\ker \PPhi(1)]_\complement$ since the definition of $\PPhi(1)^g$ does so. However, it turns out that if the direct sum holds for a specific choice of $[\ran \PPhi(1)]_\complement$ and $[\ker \PPhi(1)]_\complement$, then it also holds  for any arbitrary possible choice of those; see Lemma \ref{lemi2}-\ref{lemi21aaa}. Therefore, the choice of  $[\ran \PPhi(1)]_\complement$ and  $[\ker \PPhi(1)]_\complement$ does not affect the subsequent results, and also may be arbitrarily chosen among possible candidates without loss of generality; for example, in a Hilbert space setting, we may assume that $[\ran \PPhi(1)]_\complement=[\ran \PPhi(1)]^\perp$ and $[\ker \PPhi(1)]_\complement=[\ker \PPhi(1)]^\perp$. 	
	\end{remark} 	
	
	\begin{remark}	\label{remi2conditionb}
		Suppose that $\mathcal B=\mathcal H$ and $\Phi(z)$ satisfies the Fredholm assumption \eqref{eqfred}. In this case, we may assume that $[\ran \PPhi(1)]_\complement=[\ran \PPhi(1)]^\perp$ and $[\ker \PPhi(1)]_\complement=[\ker \PPhi(1)]^\perp$ in the I(2) condition given above without loss of generality (see Remark \ref{remi2conditiona}), then $\PPhi(1)^g$ becomes the Moore-Penrose inverse $\PPhi(1)^\dag$ of $\PPhi(1)$ (see Appendix \ref{sginverse}). In this setting,   \citet[Theorem 4.2 and Remark 4.7]{BS2018}  showed that $\PPhi(z)^{-1}$ has a pole of order 2 at $z=1$ if and only if $\mathcal H^p=(\ran \PPhi(1) + \ker \PPhi(1)) \oplus \PPhi(1)^\dag \mathcal K$ for $\mathcal K \neq \{0\}$.\footnote{In fact, the direct sum condition given by \citet[Theorem 4.2 and Remark 4.7]{BS2018} is slightly different. However, with a simple algebra, it can be shown that the direct sum given in the I(2) condition is equivalent to $\mathcal B = (\ran \PPhi(1) + \ker \PPhi(1)) \oplus (I_p- \PPhi(1)^g) \mathcal K$, which is exactly comparable with their direct sum condition.} However, it is worth noting that Fredholmness of $\PPhi(z)$ requires $\ker \PPhi(1)$ (and thus $\mathcal K$) to be finite dimensional while no such restriction required in our I(2) condition; see Example \ref{exi2001} to appear later.
	\end{remark} 
	Our next result shows that the proposed I(2) condition is indeed necessary and sufficient for $\PPhi(z)^{-1}$ to have a pole of order 2 at $z=1$, and provides a partial characterization of the principal part of the Laurent series; a more detailed characterization of the principal part can be obtained in terms of some operators associated with $\PPhi(z)$, but which is postponed to Appendix \ref{app:companioni2} since more preliminary results and concepts, which are to be developed in Appendix \ref{appenA}, are required.  
	\begin{proposition} \label{grti2p}
		Suppose that Assumption \ref{assulinear2} holds. Then \(\PPhi(z)^{-1}\) has a pole of order 2 if and only if the I(2) condition holds. Under the I(2) condition,  the following holds for some $\eta >0$,
		\begin{align}
			(1-z)^2\PPhi(z)^{-1} = - \NN_{-2}  + (1-z)(\NN_{-2}+\PP) + (1-z)^2H(z),\quad z \in D_{1+\eta},  \label{i2laurenta0}
		\end{align} 
		where  $\NN_{-2}$ satisfies  $\ran \NN_{-2}= \mathcal K$, $\hspace{-0.0em}H(z)\hspace{-0.0em}$ is the holomorphic part of the Laurent series of $\hspace{-0.0em}\PPhi(z)^{-1}\hspace{-0.0em}$, and each Maclaurin series of $\hspace{-0.1em}(1-z)^2\PPhi(z)^{-1}\hspace{-0.1em}$ and $\hspace{-0.1em}H(z)\hspace{-0.1em}$ converges on \(\hspace{-0.1em}D_{1+\eta}\hspace{-0.1em}\). 
	\end{proposition}

	Examples of the use of Proposition \ref{grti2p} for verifying that $\widetilde{\Phi}(z)^{-1}$ has a pole of order 2 will be given  at the end of this section. Proposition \ref{grti2p} provides a characterization of the local behavior of $\PPhi(z)^{-1}$ around $z=1$ under the I(2) condition, which, together with  Proposition \ref{propa1}, leads to our first version of the Granger-Johansen representation theorem for I(2) AR processes as follows.   
	
	\begin{proposition}\label{grt2pp}
		Suppose that Assumption \ref{assulinear2} holds. Under the I(2) condition, a sequence \(\{X_t\}_{t \geq -p+1}\) satisfying \eqref{arlaw} allows the Johansen I(2) representation \eqref{bndecom3} with 
		\begin{align}
			\Upsilon_{-2} \hspace{-0.2em} = \hspace{-0.2em}- \Pi_p \NN_{-2} \Pi_p^\ast, \quad \Upsilon_{-1}\hspace{-0.2em} = \hspace{-0.2em} \Pi_p \left(\NN_{-2} + \PP\right)\Pi_p^\ast, \quad \nu_t \hspace{-0.2em}=\hspace{-0.2em} \Pi_pH(L)\Pi_p^\ast \varepsilon_t \hspace{-0.2em}=\hspace{-0.2em} \sum_{j=0}^\infty \Pi_p \PPhi_1^j (\idbp-\PP) \Pi_p^\ast \varepsilon_{t-j},\label{bndecom3i2}
		\end{align}
		where \(\NN_{-2}\), \(\PP\) and $H(z)$ are given in Proposition \ref{grti2p}, and $\Pi_p$ and $\Pi_p^\ast$ are given in \eqref{coorproj}. Moreover, the AR($p$) law of motion \eqref{arlaw} does not allow I(2) solutions if the I(2) condition is not satisfied. 
	\end{proposition}
	Note that the representation \eqref{bndecom3} with \eqref{bndecom3i2} is similar to the Beveridge-Nelson decomposition \eqref{eqi1bn}  of  an I(2) cointegrated linear process. From our discussion in Section \ref{ssinteg}, we may deduce that $f\in \Ann(\Pi_p\NN_{-2}\Pi_p^\ast)$ (resp.\  $f\in \Ann(\Pi_p\NN_{-2}\Pi_p^\ast)\cap \Ann(\Pi_p\PP\Pi_p^\ast)$) is a cointegrating functional (resp.\ a second-order cointegrating functional). Appendix \ref{app:companioni2} provides characterizations of $\NN_{-2}$ and $\PP$ in terms of certain operators associated with $\PPhi(z)$, which complements the results given by Proposition \ref{grt2pp}.
	
	In this section, we have shown that the AR($p$) law of motion \eqref{arlaw} admits I(2) solutions if and only if the I(2) condition holds, and provided a partial characterization of such solutions. All these results are obtained resorting to the companion form representation of \eqref{arlaw}, hence it is not clearly revealed how the I(2) condition and the cointegrating behavior of I(2) solutions are related to the structure of the original AR polynomial $\Phi(z)$. This issue will be addressed in the next section by providing a more detailed characterization of I(2) solution in terms of operators associated with $\Phi(z)$.
	Before closing this section, we give examples of the use of Proposition \ref{grti2p} for verifying that $\widetilde{\Phi}(z)^{-1}$ has a pole of second order. 
	\begin{example} \label{exi2000}	
		Consider Example \ref{examplei1i2}, where we showed that $\widetilde{\Phi}(z)^{-1}$ does not have a simple pole at $z=1$, hence know that Assumption \ref{assulinear2}-\ref{assulinear2a} holds. Note that  $\ran  \widetilde{\Phi}(1)$ (resp.\ $\ker \PPhi(1)$) allows a complementary subspace $[\ran \widetilde{\Phi}(1)]_{\complement}$ (resp.\ $[\ker \PPhi(1)]_\complement$); specifically, $[\ran \widetilde{\Phi}(1)]_{\complement}$ and $[\ker \PPhi(1)]_\complement$ may be set to  	
		\begin{equation*}
			[\ran \widetilde{\Phi}(1)]_{\complement} = \{(b_1,0,0,\ldots),\, b_1 \in \mathbb{C}\}, \quad  \quad 		[\ker \widetilde{\Phi}(1)]_{\complement} =  \{(b_1,0,b_2,b_3,\ldots),\, :\,b_j \in \mathbb{C}, \,  \lim_{j\to \infty} b_j = 0 \}.
		\end{equation*}
Then it can be shown that Assumption \ref{assulinear2}-\ref{assulinear2b} is also satisfied; observe that $\ran \Phi(1) \oplus \mathsf{R} = \ran \Phi(1) + \Phi^{(1)}(1)\ker \Phi(1) = \ran \PPhi(1) + \ker \PPhi(1)=\ran \PPhi(1)$ and $\mathsf{K}=\ker \PPhi(1)$, hence  the required direct sums hold for $\mathsf{R}_\complement=[\ran \PPhi(1)]_\complement$ and $\mathsf{K}_\complement = \{0\}$.
	For any \((b_1,0,\ldots) \in[\ran  \widetilde{\Phi}(1)]_{\complement}\), we find that $ -\widetilde{\Phi}(1)(b_1,0,\ldots) = (0,b_1,0,\ldots)$.
	Pre-composing  both sides of this equation with $\PPhi(1)^g$, we obtain 
			\begin{equation} \label{eqran2}
			(-b_1,0,\ldots)=  \widetilde{\Phi}(1)^g  (0,b_1,0,\ldots),
		\end{equation}
	where the equality holds since \((b_1,0,\ldots) \in [\ker \widetilde{\Phi}(1)]_{\complement} \) and $\PPhi(1)^g\PPhi(1) =\PR_{[\ker \PPhi(1)]_\complement}$.		From \eqref{eqran1}, \eqref{eqran2} and the fact that \(\mathcal K =  \ran \widetilde{\Phi}(1)\cap \ker \widetilde{\Phi}(1)= \ker \widetilde{\Phi}(1)\), we deduce that $\widetilde{\Phi}(1)^g  \mathcal K = \{(b_1,0,\ldots),\, b_1 \in \mathbb{C}\}$, which is a complementary subspace of  \(\ran  \widetilde{\Phi}(1) + \ker \widetilde{\Phi}(1)\).  Thus, $\widetilde{\Phi}(z)^{-1}$ has a pole of order 2 at $z=1$.
	\end{example}	
	
	\begin{example}  \label{exi2001} In the setting of Examples \ref{examplei1i2} and \ref{exi2000}, we observed that $\ker \PPhi(1)$ and $\mathcal K$ are finite dimensional under the I(2) condition. However, these subspaces may not be finite dimensional in general. To see this, we will consider slight modifications of $\Phi_1$; under any of the changes in $\Phi_1$ to be given below, it can be easily shown that Assumption \ref{assulinear2} is still satisfied.  
		 We first replace \eqref{eqopex} with the following, 
		\begin{equation*}
			\Phi_1(a_1,a_2,a_3,a_4\ldots) = (a_1,a_2+a_1, a_3+a_1, a_4, a_5+a_4, a_6+a_4, a_7, a_8+a_7,a_9+a_7,\ldots). \label{eqopex2}
		\end{equation*}
		It can be easily shown that $\PPhi(z)$ is invertible for $z\neq 1$. We note that  $\ran \PPhi(1)$ and $\ker \PPhi(1)$ are given by 
		\begin{align}
			&\ran \widetilde{\Phi}(1)= \{(0,b_1,b_1,0,b_2,b_2,0,b_3,b_3,\ldots) : \,  b_j \in \mathbb{C}, \,\lim_{j\to \infty} b_j = 0  \},  \label{eqran1a}  \\ &\ker  \widetilde{\Phi}(1)= \{(0,b_1,b_2,0,b_3,b_4,0, b_5,b_6,\ldots) : \, b_j \in \mathbb{C}, \,\lim_{j\to \infty} b_j = 0  \}.  \label{eqran1aa} 
		\end{align}
		Since  $\mathcal K = \ran \PPhi(1) \cap \ker \PPhi(1)= \ran \PPhi(1)$, we know from \eqref{eqran1a} and \eqref{eqran1aa} that both $\ker \PPhi(1)$ and $\mathcal K$ are infinite dimensional.  With some algebra, it may be deduced that $\ran \PPhi(1)$ (resp.\ $\ker \PPhi(1)$) allows a complementary subspace $[\ran \PPhi(1)]_\complement$ (resp.\ $[\ker \PPhi(1)]_\complement$), and especially  $[\ker \PPhi(1)]_\complement$  may be set to 
		\begin{equation}
			[\ker \PPhi(1)]_{\complement} = \{(b_1,0,b_2,b_3,0,b_4,b_5,0,b_6,\ldots) : \,  b_j \in \mathbb{C}, \,\lim_{j\to \infty} b_j = 0\}.  \label{eqran1aa2} 
		\end{equation}
		Moreover, note that for $(b_1,0,b_2,b_3,0,b_4,b_5,0,b_6,\ldots) \in [\ker \PPhi(1)]_{\complement}$,
		\begin{equation}
			-\PPhi(1)(b_1,0,b_2,b_3,0,b_4,b_5\ldots) = (0,b_1,b_1,0,b_3,b_3,\ldots). \label{exeq0112}
		\end{equation}
		Since $(b_1,0,b_2,b_3,0,b_4,b_5\ldots) \in [\ker \PPhi(1)]_{\complement}$ and $\PPhi(1)^g\PPhi(1) = \PR_{[\ker \PPhi(1)]_\complement}$,  we find the following equality by pre-composing  both sides of \eqref{exeq0112} with $\PPhi(1)^g$,
		\begin{align}
			&-(b_1,0,b_2,b_3,0,b_4,b_5\ldots) = \PPhi(1)^g(0,b_1,b_1,0,b_3,b_3,\ldots).  \label{eqran1aaa} 
		\end{align}
		From \eqref{eqran1a}, \eqref{eqran1aa} and \eqref{eqran1aaa}, we deduce that $\ran \PPhi(1) +\ker \PPhi(1) = \ker \PPhi(1)$ and $\PPhi(1)^g \mathcal K  = [\ker \PPhi(1)]_\complement$, hence $\widetilde{\Phi}(z)^{-1}$ has a pole of order 2 at $z=1$.
		
		Now we replace \eqref{eqopex} with the following, 
		\begin{equation*}
			\Phi_1(a_1,a_2,a_3,a_4\ldots) = (a_1,a_2+a_1, a_3, a_4, a_5,\ldots). \label{eqopex2a}
		\end{equation*}
		Then,  we find that 
		\begin{align}
			&\ran \widetilde{\Phi}(1)= \{(0,b_1,0,0,\ldots) : \,  b_1 \in \mathbb{C}\}, \quad \ker  \widetilde{\Phi}(1)= \{(0,b_1,b_2,b_3,\ldots) : \, b_j \in \mathbb{C}, \,\lim_{j\to \infty} b_j = 0 \},  \label{eqran1aa2a} 
		\end{align}
		and thus $\mathcal K = \ran \widetilde{\Phi}(1)$. From similar arguments, it can be shown that our I(2) condition holds. In this case, as may be deduced from \eqref{eqran1aa2a}, $\ker \PPhi(1)$ is infinite dimensional but $\mathcal K$ is finite dimensional.  
	\end{example}	
	\subsubsection{Further characterization of I(2) solutions}	\label{sec:repi2noncompanion}
	Suppose that $\mathcal B=\mathbb{R}^n$  or $\mathbb{C}^n$ and the Johansen I(1) condition fails. With the notation introduced for \eqref{eqjohansen}, we let $\varpi$ and $\varrho$ be full-rank $(n-r) \times s$ matrices satisfying  $\sbot{\alpha}^\intercal \Phi^{(1)}(1) \sbot{\beta} = \varpi\varrho^\intercal$, where $s<n-r$. Let $\sbot{\varpi}$ (resp.\ $\sbot{\varrho}$) be a full-rank $n \times (n-r-s)$ matrix whose columns are orthogonal to those of $(\alpha, \sbot{\alpha} {\varpi})$ (resp.\ $(\beta, \sbot{\beta} {\varrho})$). Without loss of generality, we may assume that $\sbot{\varpi}^\intercal\sbot{\varpi}=\sbot{\varrho}^\intercal\sbot{\varrho}=I_{n-r-s}$ (the identity matrix of dimension $n-r-s$). \cite{johansen1992representation,Johansen1996} provides a necessary and sufficient condition for $\Phi(z)^{-1}$ to have a pole of order 2 at $z=1$, which is given by that 
	\begin{equation}
		\sbot{\varpi}^\intercal \left(\frac{1}{2}\Phi^{(2)}(1) - \Phi^{(1)}(1)\Phi(1)^\dag \Phi^{(1)}(1)\right) \sbot{\varrho} \,\, \text{ is invertible,} \label{eqjohansen2}
	\end{equation}
	where $\Phi(1)^\dag$ is the Moore-Penrose inverse of $\Phi(1)$.  
	To provide a natural generalization of the Johansen I(2) condition that can be applied to our Banach space setting, as we did in Section \ref{sec:repi1noncompanion} for the I(1) case, we first review some preliminary results that hold under Assumption \ref{assulinear2} and fix notation. 	
	
Under Assumption \ref{assulinear2}, we let $\PR_{[\ran \Phi(1)]_\complement}$ and  $\PR_{[\ker \Phi(1)]_\complement}$ be defined as in \eqref{projeq}, and also let  $\Phi(1)^g$ denote the generalized inverse satisfying the following properties:
\begin{equation}
	\Phi(1) \Phi(1)^g= \idbp-\PR_{[\ran \Phi(1)]_\complement}, \quad \,\, \Phi(1)^g\Phi(1) = \PR_{[\ker \Phi(1)]_\complement}.\label{direct1apro2a} 
\end{equation}
The operators  $\PR_{[\ran \Phi(1)]_\complement}$, $\PR_{[\ker \Phi(1)]_\complement}$ and  $\Phi(1)^g$ are all uniquely defined for any fixed choice of $[\ran \Phi(1)]_\complement$ and $[\ker \Phi(1)]_\complement$.	If $\mathcal B=\mathcal H$ and thus $[\ran \Phi(1)]_\complement=[\ran \Phi(1)]^\perp$ and $[\ker \Phi(1)]_\complement=[\ker\Phi(1)]^\perp$, then $\Phi(1)^g$ is equivalent to the Moore-Penrose inverse $\Phi(1)^\dag$ of $\Phi(1)$; see Appendix \ref{sginverse}.
	Moreover, we note that the direct sum conditions given in Assumption \ref{assulinear2} can be combined and equivalently formulated as follows, 
	\begin{equation}\label{newdirec}
		\mathcal B = \ran \Phi(1) \oplus \mathsf{R} \oplus \mathsf{R}_\complement = [\ker \Phi(1)]_\complement \oplus \mathsf{K} \oplus \mathsf{K}_\complement,
	\end{equation} 
	where $[\ran \Phi(1)]_\complement = \mathsf{R} \oplus \mathsf{R}_\complement$ and $\ker \Phi(1)=\mathsf{K} \oplus \mathsf{K}_\complement$. Then for any possible choice of $\mathsf{R}_\complement$ and $\mathsf{K}_\complement$, we may  define the projections $\PR_{\mathsf{R}_\complement}$ and $\PR_{\mathsf{K}}$ satisfying
	\begin{equation}
		\ran \PR_{\mathsf{R}_\complement}=\mathsf{R}_\complement, \quad \ker \PR_{\mathsf{R}_\complement}=\ran \Phi(1) \oplus \mathsf{R}, \quad 	\ran \PR_{\mathsf{K}}=\mathsf{K}, \quad \ker \PR_{\mathsf{K}}=[\ker \Phi(1)]_\complement \oplus \mathsf{K}_\complement; \label{projeceq2}
	\end{equation}	  
see \citet[Theorem 3.2.11]{megginson1998}.
	Suppose that $\mathcal B = \mathbb{R}^n$ or $\mathbb{C}^n$,  $\mathsf{R}_\complement = \mathsf{R}^\perp$ (the orthogonal complement to $\ran \Phi(1)\oplus \mathsf{R}$), and $\mathsf{K}_\complement = \mathsf{K}^\perp$ (the orthogonal complement to $[\ker \Phi(1)]^\perp \oplus \mathsf{K}$) hold. Then the Johansen I(2) condition given by \eqref{eqjohansen2} can be alternatively understood as invertibility of $\PR_{\mathsf{R}^\perp} (\frac{1}{2}\Phi^{(2)}- \Phi^{(1)}(1)\Phi(1)^\dag \Phi^{(1)}(1))\PR_{\mathsf{K}}$ as a map from $\mathsf{K}$ to $\mathsf{R}^\perp$.  Given this observation, we know that the Johansen I(2) condition is a special case of the following more general condition stated without the notion of an orthogonal complement in a general Banach space setting: for some choice of  $\mathsf{R}_\complement$ and  $\mathsf{K}_\complement$,    
	\begin{equation}
		\Lambda_{2,\complement} \coloneqq \PR_{\mathsf{R}_{\complement}}\left(\frac{1}{2}\Phi^{(2)}(1) - \Phi^{(1)}(1)\Phi(1)^g\Phi^{(1)}(1) \right)\PR_{\mathsf{K}}:\mathsf{K} \mapsto {\mathsf{R}_{\complement}} \text{ is invertible.} \label{eqjohan2} 
	\end{equation}
	We will show in Proposition \ref{propjohanseni2} below that the above condition is equivalent to our I(2) condition developed in a general Banach space setting, and also characterize $\Upsilon_{-2}$ and  $\Upsilon_{-1}$ in terms of some operators associated with the  AR polynomial $\Phi(z)$ using this equivalence. To simplify mathematical expressions, we employ the following notation: \vspace{-0.7em}
	\begin{align}
		&\mathrm{M}_1 = \PR_{[\ran \Phi(1)]_\complement}\Phi^{(1)}(1)(I-\PR_{[\ker \Phi(1)]_\complement}), \notag\\
		&\mathrm{M}_2 = \frac{1}{2}\Phi^{(2)}(1) - \Phi^{(1)}(1)\Phi(1)^g \Phi^{(1)}(1), \notag \\
		&\mathrm{M}_3 = \frac{1}{6} \Phi^{(3)}(1) - \Phi^{(1)}(1) \Phi(1)^g\Phi^{(1)}(1) \Phi(1)^g\Phi^{(1)}(1).\notag
	\end{align}
	As the last piece of our preliminary results for  Proposition \ref{propjohanseni2}, we note that the generalized inverse $\mathrm{M}_1^g$ of $\mathrm{M}_1$ is  uniquely defined for any choice of $\mathsf{R}_\complement$ and $\mathsf{K}_\complement$ under Assumption \ref{assulinear2}, and  it satisfies $\ran \mathrm{M}_1^g = \mathsf{K}_\complement$ and $\ker \mathrm{M}_1^g = \ran \Phi(1) \oplus \mathsf{R}_\complement$ (Lemma \ref{lemi2add}).  

	\begin{proposition}\label{propjohanseni2}
		Suppose that Assumption \ref{assulinear2} holds. Then the following conditions are equivalent.
		\begin{enumerate}[\normalfont(i)]
			\item The I(2) condition holds.
			\item $\Lambda_{2,\complement}:\mathsf{K}\mapsto \mathsf{R}_\complement$ is invertible for some choice of  $\mathsf{R}_\complement$ and $\mathsf{K}_\complement$.	
			\item $\Lambda_{2,\complement}:\mathsf{K}\mapsto \mathsf{R}_\complement$ is invertible for every possible choice of $\mathsf{R}_\complement$ and $\mathsf{K}_\complement$.		
		\end{enumerate}
		Let any of the above conditions hold. Then for any  choice of $\mathsf{R}_\complement$ and $\mathsf{K}_\complement$,  a sequence \(\{X_t\}_{t \geq -p+1}\) satisfying \eqref{arlaw} allows the Johansen I(2) representation \eqref{bndecom3} with $\Upsilon_{-2}$  satisfying 		
		\begin{equation} \label{eqn2laurent}
			(I-\PR_{\mathsf{K}})\Upsilon_{-2}= \Upsilon_{-2} (I-\PR_{\mathsf{R}_\complement}) = 0,  \quad\quad  \Upsilon_{-2}:\mathsf{R}_\complement\mapsto\mathsf{K}  =  \Lambda_{2,\complement}^{-1}, 
		\end{equation}
		and  $\Upsilon_{-1}$ satisfying	 \vspace{-00em} 	
		\begin{align*} 
			&(I-\PR_{\mathsf{K}})\Upsilon_{-1}	(I-\PR_{\mathsf{R}_\complement}) = - \mathrm{M}_1^g, \\
			&(I-\PR_{\mathsf{K}})\Upsilon_{-1}	\PR_{\mathsf{R}_\complement} =\left(\Phi(1)^g \Phi^{(1)}(1) +   \mathrm{M}_1^g \mathrm{M}_2  \right) \Upsilon_{-2}, \\
			&\PR_{\mathsf{K}}\Upsilon_{-1}		(I-\PR_{\mathsf{R}_\complement}) = \Upsilon_{-2} \left(\Phi^{(1)}(1) \Phi(1)^g	+ \mathrm{M}_2  \mathrm{M}_1^g \right), \\
			&\PR_{\mathsf{K}}\Upsilon_{-1}	\PR_{\mathsf{R}_\complement} = \Upsilon_{-2}  \left(\mathrm{M}_3 -\mathrm{M}_2\Phi(1)^g\Phi^{(1)}(1) -\Phi^{(1)}(1)\Phi(1)^g \mathrm{M}_2 - \mathrm{M}_2  \mathrm{M}_1^g \mathrm{M}_2 \right)\Upsilon_{-2}.
		\end{align*}	
	\end{proposition}
	Proposition \ref{propjohanseni2} not only shows that the condition given by \eqref{eqjohan2} is equivalent to our I(2) condition developed in a general Banach space setting, but also implies that $\mathsf{R}_\complement$ and $\mathsf{K}_\complement$ can be arbitrarily chosen among possible candidates; this, of course,  means that $\mathsf{R}_\complement$ and $\mathsf{K}_\complement$ can always be fixed to the relevant orthogonal complements in a Hilbert/Euclidean space setting without loss of generality. Moreover, Proposition \ref{propjohanseni2} characterizes  in detail how $\Upsilon_{-2}$ and $\Upsilon_{-1}$ act on $\mathcal B$ satisfying the tripartite decompositions given by \eqref{newdirec}.  Based on such characterizations,  we can further  characterize I(2) solutions; see Remarks \ref{remcointeg02}-\ref{app:polycointeg1} below. Some more remarks are given for comparison between our I(2) representation result given by Proposition \ref{propjohanseni2} and the existing I(2) results developed in a Hilbert/Euclidean space setting.

	\begin{remark} \label{remcointeg02}
		From \eqref{eqn2laurent}, we know that the attractor space associated with I(2) solutions is given by $\ran \Upsilon_{-2} = \mathsf{K}$. Then Proposition \ref{propbasic} implies that  a cointegrating functional $f$ satisfies $f =g(I-\PR_{\mathsf{K}})$, where $g \in \mathcal B'$.  This result holds for any arbitrary choice of $\mathsf{K}_\complement$; a different choice of $\mathsf{K}_\complement$ only affects the definition of $\PR_{\mathsf{K}}$ without any further changes.
	\end{remark}  
	
	\begin{remark} \label{app:polycointeg0}
		Using the results given in Proposition \ref{propjohanseni2}, we can obtain a stronger characterization of the cointegrating behavior of I(2) solutions as follows: for any nonzero element $f\in \mathcal B'$ and for some $r \in \{0,1\}$,		
		\begin{align}  
			&\{f(X_t)\}_{t \geq 0} \text{ is I($r$)} \,\, \text{if and only if} \,\,f \in \Ann(\mathsf{K}), \label{eqcointeg1} \\ 
			&\{f(X_t)\}_{t \geq 0} \text{ is I($0$)}\,\, \text{if and only if} \,\, f \in \Ann (\ker {\Phi}(1)) \cap \Ann(\Phi(1)^g \Phi^{(1)}(1)\mathsf K). \label{eqcointeg2}
		\end{align}
		Obviously, \eqref{eqcointeg1} (resp.\ \eqref{eqcointeg2}) identifies the cointegrating space (the collection of the second-order cointegrating funcitonals). A  detailed discussion including our proofs of these results is given in  Appendix \ref{app:polycointeg}. 	If $\mathcal B = \mathcal H$ and thus the complementary subspaces are all set to the relevant orthogonal complements, then $\Phi(1)^g$ is equivalent to the Moore-Penrose inverse $\Phi(1)^\dag$ (see Appendix \ref{sginverse}) and $f \in \Ann(V)$ is given by the map $\langle \cdot,v \rangle$ for some $v \in V^\perp$ (see Remark \ref{johansencointeg}). In this case, \eqref{eqcointeg1} and \eqref{eqcointeg2} can be reformulated as follows:  for any $v \in \mathcal H \setminus \{0\}$ and for some $r \in \{0,1\}$, 
		\begin{align}
			&\text{$\{\langle X_t,v\rangle\}_{t \geq 0}$ is I($r$)}\,\, \text{if and only if} \,\,v \in \mathsf{K}^\perp, \label{eqcointeg1a} \\ 
			&\text{$\{\langle X_t,v\rangle\}_{t \geq 0}$ is I($0$)}\,\, \text{if and only if} \,\, v \in [\ker \Phi(1)]^\perp \cap [\Phi(1)^\dag \Phi^{(1)}(1)\mathsf{K}]^\perp. \label{eqcointeg2a}
		\end{align}
		If $\Phi(z)$ satisfies the Fredholm assumption \eqref{eqfred} and, as a result, both $\ker \Phi(1)$ and $\mathsf{K}$ is finite dimensional, our characterization given by \eqref{eqcointeg1a} and \eqref{eqcointeg2a} reduces to that provided  by  \citet[Remarks 4.5 and 4.6]{BS2018}. It, however, should be noted that \eqref{eqcointeg1a} and \eqref{eqcointeg2a} hold without the Fredholm assumption and, thus,  $\ker \Phi(1)$ and $\mathsf{K}$ are not required to be finite dimensional; noting that $\ker \PPhi(1) = \ker \Phi(1)$ and $\mathcal K = \mathsf{K}$ hold in the case $p=1$, we see  from Example \ref{exi2001} that those can be infinite dimensional.
	\end{remark} 
	
	\begin{remark}[Polynomial cointegration for an I(2) AR process] \label{app:polycointeg1}
		For any cointegrating functional $f \in \Ann(\mathsf{K})$, $f$ may or may not satisfy $f \in \Ann(\ker \Phi(1))$ since $\mathsf{K} \subset \ker \Phi(1)$. If $f \in \Ann(\ker \Phi(1))$ and one combines levels and first differences as in $f(X_t) -f(\Phi(1)^g \Phi^{(1)}(1) \Delta X_t)$, then such a sequence is always I(0); this phenomenon does not occur if $f \notin \Ann(\ker \Phi(1))$. A  detailed discussion including our proof of this result is given in  Appendix \ref{app:polycointeg2}. The case where the sequence of $f(X_t)$ and that of $f(\Phi(1)^g \Phi^{(1)}(1) \Delta X_t)$ are both I(1) may be understood as  polynomial cointegration or multicointegration \citep{yoo1987,granger1989investigation,engsted1997granger,kheifets2019fully} in our setting; a detailed  treatment of this topic for the case $\mathcal B=\mathcal H$ is given by \citet[Section 4.2]{Franchi2017b}. In the case $\mathcal B = \mathcal H$, $\Phi(1)^g$ may be set to the Moore-Penrose inverse $\Phi(1)^\dag$ and $f\in \Ann(\mathsf{K})$ (resp.\ $f\in \Ann(\ker \Phi(1))$) is identified as the map $\langle \cdot, v \rangle$ for some  $v \in \mathsf{K}^\perp$ (resp.\ $v \in [\ker \Phi(1)]^\perp$). Then the above characterization given in a Banach space setting can be rewritten as follows: for $f \in \mathsf{K}^\perp$, the sequence of $\langle X_t,v \rangle - \langle \Phi(1)^\dag \Phi^{(1)}(1) \Delta X_t,v \rangle$ is I(0) if and only if $f \in [\ker \Phi(1)]^\perp$. If the Fredholm assumption \eqref{eqfred} holds, this result becomes consistent with the characterization of polynomial cointegration given by \citet[Remark 4.10]{Franchi2017b}. However, $\ker \Phi(1)$ and $\mathsf{K}$ are necessarily finite dimensional under the Fredholm assumption, while such restrictions are not required for our results given in this remark; see Remark \ref{app:polycointeg0} and  Example \ref{exi2001}.
	\end{remark}


	\begin{remark} \label{bsi2condition}
		Suppose that $\mathcal B = \mathcal H$ and all the complementary subspaces are set to the relevant orthogonal complements. 
		In this case, \eqref{eqjohan2} reduces to that 
		\begin{equation}
			\PR_{\mathsf{R}^\perp}\left(\frac{1}{2}\Phi^{(2)}(1) - \Phi^{(1)}(1)\Phi(1)^\dag\Phi^{(1)}(1) \right)	\PR_{\mathsf{K}} : \mathsf{K} \mapsto \mathsf{R}^\perp \quad \text{is invertible,} \label{eq001rem}
		\end{equation}
		where $\PR_{\mathsf{R}^\perp}$ (resp.\ $\PR_{\mathsf{K}}$) is the orthogonal projection onto $\mathsf{R}^\perp$ (resp.\ $\mathsf{K}$). If $\Phi(z)$ satisfies the Fredholm assumption \eqref{eqfred} and thus both $\ker \Phi(1)$ and $\mathsf{K}$ are finite dimensional, then \eqref{eq001rem} reduces to the invertibility condition for the existence of I(2) solutions given by  \citet[Theorem 4.2]{BS2018}.  Their condition is given together with an equivalent  nonorthogonal direct sum condition in $\mathcal H$; on the other hand, \cite{Franchi2017b} provided an orthogonal direct sum condition for the existence of I(2) solutions under the same setting, see Section 4 of their paper. Under any of these equivalent conditions, solutions to the AR($p$) law of motion satisfies the Johansen I(2) representation \eqref{bndecom3} for  $\Upsilon_{-2}$ and $\Upsilon_{-1}$ of finite rank; it can be deduced from our characterizations of $\Upsilon_{-2}$ and $\Upsilon_{-1}$ given in Proposition \ref{propjohanseni2} that the ranges of these operators are finite dimensional if both $\ker \Phi(1)$ and $\mathsf{K}$ are finite dimensional. Hence the random walk component of I(2) solutions is intrinsically a finite dimensional unit root process. However, in our results for the case $\mathcal B=\mathcal H$, either $\ker \Phi(1)$ or $\mathsf{K}$ is not required to be finite dimensional and the random walk component can take values in an infinite dimensional subspace; we see this by noting that  $\ran \Upsilon_{-2} = \mathsf{K}$ and $\mathsf{K}$ can be infinite dimensional as in Example \ref{exi2001}. Thus, Proposition \ref{propjohanseni2} complements the earlier I(2) representation results developed in a Hilbert space setting.
	\end{remark}
	
	\begin{remark} \label{johansencointeg2i2}
		Suppose that $\mathcal B=\mathbb{R}^n$ or $\mathbb{C}^n$ and all the complementary subspaces are set to the relevant orthogonal complements. In this case, $\Upsilon_{-2}$ and  $\Upsilon_{-1}$ may be understood as  $n\times n$ matrices.  Using the notation introduced for \eqref{eqjohansen2}, the results given in  \eqref{eqn2laurent} can be recast as $\varrho(\varrho^\intercal \varrho)^{-1}\varrho^\intercal\Upsilon_{-2}= \Upsilon_{-2} \varpi(\varpi^\intercal \varpi)^{-1}\varpi^\intercal = 0$ and $\sbot{\varrho}^\intercal\Upsilon_{-2}\sbot{\varpi} = (\sbot{\varpi}^\intercal \left(\frac{1}{2}\Phi^{(2)}(1) - \Phi^{(1)}(1)\Phi(1)^\dag \Phi^{(1)}(1)\right) \sbot{\varrho} )^{-1}$. From these results, we find that  
		\begin{equation}\label{eqrem001add2i2}
			\Upsilon_{-2} =\sbot{\varrho} \left[\sbot{\varpi}^\intercal \left(\frac{1}{2}\Phi^{(2)}(1) - \Phi^{(1)}(1)\Phi(1)^\dag \Phi^{(1)}(1)\right) \sbot{\varrho}\right]^{-1}\sbot{\varpi}^\intercal.
		\end{equation}
		Let $\alpha_1 = \sbot{\alpha} \varpi$, $\beta_1 = \sbot{\beta} \varrho$, $\overline{\alpha}_1 = {\alpha_1} ({\alpha}_1^\intercal {\alpha}_1)^{-1}$ and $\overline{\beta}_1 = {\beta}_1 ({\beta}_1^\intercal {\beta}_1)^{-1}$. Then the operator $\mathrm{M}_1$  can be written as $\mathrm{M}_1 =  \sbot{\alpha}\varpi \varrho^\intercal \sbot{\beta}^\intercal$ and its Moore-Penrose inverse $\mathrm{M}_1^\dag$ is given by $\mathrm{M}_1^\dag= \overline{\beta}_1\overline{\alpha}_1^\intercal$. By replacing $\Upsilon_{-2}$ and $\mathrm{M}_1^g$ with \eqref{eqrem001add2i2} and $\mathrm{M}_1^\dag$ respectively in our characterization of $\Upsilon_{-1}$ given in Proposition \ref{propjohanseni2}, we can obtain $\Upsilon_{-1}$ characterized as an $n \times n$ matrix. 	These expressions for $\Upsilon_{-2}$ and $\Upsilon_{-1}$ are equivalent to those in  Johansen's representation of I(2) AR processes \citetext{see e.g.,\ \citeauthor{johansen2008},  \citeyear{johansen2008}, Theorem 5}.  The case where  $\mathcal B=\mathbb{R}^n$ or $\mathbb{C}^n$ was discussed in detail as a special case of a Hilbert space  in the recent literature, so the results given in this remark were already noted in the earlier works; see e.g.,\  \citet[Remark 4.4]{BS2018}.
	\end{remark}		

	\section{Concluding remarks}\label{sconclude}
	In this paper, we provide a suitable notion of cointegration in Banach spaces and study theoretical properties of the cointegrating space. We also extend the Granger-Johansen representation theorem to a potentially infinite dimensional Banach space setting. Compared to the existing results, our representation theorems are derived under a weaker geometry of a Banach space and weaker regularity conditions on the AR polynomial. 	As a consequence, not only more general AR($p$) law of motions can be accommodated in our representation theory, but also our results do not place any dimensionality restrictions on the random walk component of I(1) or I(2) solutions to a given AR($p$) law of motion.

	To develop our representation theorems under weaker  assumptions, this paper only focuses on the I(1) and I(2) cases; on the other hand,   \cite{Franchi2017b}  recently studied the general I($d$) case for $d \geq 1$ and provided a complete characterization of the cointegrating behavior in a convenient form based on the geometry of a Hilbert space and the spectral properties of Fredholm operator pencils. Extending their results directly to our Banach space setting, where non-Fredholm AR polynomials are allowed, seems to be nontrivial. This can certainly be further explored in the future study. 
	
	Beyond representation theory, research on the development of statistical procedures for analyzing Banach-valued cointegrated  time series also needs to be pursued to complement the  existing papers on estimation, testing, and forecasting with stationary Banach-valued time series such as e.g.,  \cite{pumo1998prediction}, \cite{bosq2002estimation},  \cite{labbas2002estimation}, \cite{dehling2005estimation}, \cite{ruiz2019strongly}, and \cite{dette2017}. There may be a lot of possibilities for further studies in this direction. 

	\bibliographystyle{apalike}
	\bibliography{C:/bibtexlib/swkrefs}
	\appendix
	\section{Preliminaries} \label{app:prelim}
	\subsection{Quotient spaces}	\label{apppre1}
	Let $V$ be a subspace of an arbitrary separable complex Banach space \(\mathcal B\) equipped with norm \(\|\cdot\|_{\mathcal B}\). The cosets of $V$ are defined as the collection of the following sets: 	$x + V = \{x + v : v \in V \}$, $x \in \mathcal B$.
	The quotient space of $V$, denoted by $\mathcal B/V$, is the vector space 	whose elements are  equivalence classes of the cosets of $V$: two cosets  $x + V$ and $y+V$ are in the same equivalence class if and only if $ x-y \in V$. In the present paper, any quotient space $\mathcal B/V$ is mostly associated with a closed subspace $V$. For such $V$, the quotient map $\pi_{\mathcal B/V}$ is defined by the map $	\pi_{{\mathcal B}/V} (x) =x + V$ for $x \in \mathcal B$ and the quotient norm $\|\cdot\|_{\mathcal B/V}$ is defined as $\|x + V\|_{\mathcal B / V} = \inf_{y \in V}\|x-y\|_{\mathcal B}$ for $x+V \in \mathcal B/V$.
	$\mathcal B / V$ equipped with the quotient norm $\|\cdot\|_{\mathcal B/V}$ is a Banach space  \citep[Theorem 1.7.7]{megginson1998}.

	\subsection{Random elements in $\mathcal B$} \label{apppre3}
	We briefly introduce Banach-valued random elements, called \(\mathcal B\)-random variables. The reader is referred to \citet[Chapter 1]{Bosq2000} for a more detailed discussion on this subject.
	
	Let $(\Omega, \mathbb{F}, \mathbb{P})$ be an underlying probability triple. A $\mathcal B$-random variable is a measurable map $X: \Omega \mapsto \mathcal B$, where $\mathcal B$ is understood to be equipped with the Borel $\sigma-$field. We say that $X$ is integrable if $E\|X\|_{\mathcal B} < \infty$. If $X$ is integrable, it turns out that there exists a unique element $EX \in \mathcal B$ such that for all $f \in \mathcal B'$, $E[f(X)] = f (EX)$.	Let $\mathfrak L^2(\mathcal{B})$ be the space of $\mathcal B$-random variables $X$ with $EX = 0$ and $E\|X\|^2 < \infty$. The covariance operator $C_X$ of $X \in \mathfrak L^2(\mathcal{B})$ is a map from $\mathcal B'$ to $\mathcal B$, defined by $C_X (f) = E[f(X)X]$ for $f \in \mathcal B'$. For $X,Y \in \mathfrak L^2(\mathcal B)$, the cross-covariance operator $C_{X,Y}$ is defined by $C_{X,Y}(f)= E[f(X)Y]$.
	
	\subsection{Generalized inverse operators}\label{sginverse}
	Let \(\mathcal B\) and \(\widetilde{\mathcal B}\) be separable complex Banach spaces and let \( A \in \mathcal L(\mathcal B,  \widetilde{\mathcal B})\).  Suppose that \(\mathcal B = \ker A \oplus V\) and \(\widetilde{\mathcal B} = \ran A \oplus W\) hold for some $V \subset \mathcal B $ and $W \subset \widetilde{\mathcal B}$. Since $A_R = A : V \mapsto \ran A$ is invertible, $A_R^{-1} : \ran A \to V$ is well defined. The generalized inverse $A^g$ of $A$ is obtained by extending the domain (resp.\ codomain) of $A_R^{-1}$ to $\mathcal B$ (resp.\ $\widetilde{\mathcal B}$); specifically, $A^g$ is defined as the map
	\begin{align*} 
		\mathcal B \ni x \quad \mapsto \quad A_R^{-1} (I-\PR_{W}) \,x \in 	\widetilde{\mathcal B},
	\end{align*}
	where \(\PR_{W}\) denotes the projection onto \(W\) along \(\ran A\).
	It can be shown that the generalized inverse \(A^g\) satisfies the following properties: $AA^gA = A$, $A^gAA^g = A^g$, $A A^g = (I-\PR_{W})$, and $A^g A =  \PR_{V}$, 
	where $\PR_{V}$ denotes the  projection on  \(V\) along \(\ker A\). 
	
	Note that $V$ and $W$ satisfying  \(\mathcal B = \ker A \oplus V\) and \(\widetilde{\mathcal B} = \ran A \oplus W\) are  not uniquely determined in general, so the above definition of $A^g$ depends on the choice of $V$ and $W$; however, for any fixed choice of $V$ and $W$, \(A^g\) is uniquely defined. In the case where  $\mathcal B$ and $\widetilde{\mathcal B}$ are Hilbert spaces, $V$ (resp.\ $W$) can be set to $[\ker A]^\perp$ (resp.\ $[\ran A]^\perp$), which makes $A^g$ become equivalent to the Moore-Penrose inverse of $A$. For a more detailed discussion on generalized inverses, see \cite{Engl1981}.

	\subsection{Operator pencils}\label{apppre4}
	Let \(U\) be some open and connected subset of the complex plane \(\mathbb{C}\). An operator pencil is an operator-valued map \(A : U \to \mathcal L(\mathcal B)\). An operator pencil \(A\) is said to be holomorphic on an open and connected set \(U_0 \subset U\) if, for each \(z_0 \in U_0\), the limit $A^{(1)}(z_0) \coloneqq \lim_{z\to z_0} {A(z)-A(z_0)}/(z-z_0)$ exists in the norm of $\mathcal L(\mathcal B)$. It turns out that if an operator pencil \(A\) is holomorphic, for every $z_0\in U_0$, we may represent \(A\) on \(U_0\) in terms of a power series $A(z) = \sum_{j=0}^\infty  A_j(z-z_0)^j$ for $z \in U_0$, where \(A_0,A_1,\ldots\) is a sequence in \(\mathcal L(\mathcal B)\). If there exists $k$ such that $A_j = 0$ for all $j \geq k$, then $A$ is called a polynomial operator pencil.  If $A_j = 0$ for all $j \geq 2$, then $A$ is called a linear operator pencil.  The collection of  \(z\in U\) at which the operator \(A(z)\) is not invertible is called the spectrum of \(A\), and denoted by \(\sigma(A)\). It turns out that the spectrum is always a closed set, and if \(A\) is holomorphic on \(U\), then \(A(z)^{-1}\) is holomorphic on \(U\setminus\sigma(A)\) \citep[p.\ 56]{Markus2012}. \(U\setminus \sigma(A)\)  is called the resolvent set of \(A\), and denoted \(\rho(A)\). If $A$ is holomorphic and $z_0$ is an isolated point of $\sigma(A)$, then $A(z)^{-1}$ allows the following Laurent series in a punctured neighborhood of $z=z_0$: 
	\begin{equation}\label{mopencil}
		A(z)^{-1} = \sum_{j= - d}^\infty A_j (z-z_0)^j , \quad d \in \mathbb{N} \cup \{\infty\}, \quad  A_j \in \mathcal L(\mathcal B).
	\end{equation} 
	By Cauchy's residue theorem, we have $A_j = -\frac{1}{2\pi i} \int_{\Gamma}\frac{A(z)^{-1}}{ (z-1)^{j+1} } dz$,  where $\Gamma \subset  \rho(\PPhi)$ is a clockwise-oriented contour around \(z_0\) such that the only element of $\sigma(A)$ included inside the contour is $z_0$.


		
	\section{Mathematical Appendix} \label{appenA} 
	We provide mathematical proofs of the results given in Sections \ref{sfcointeg}-\ref{srep}. It is sometimes convenient to consider $A\in \mathcal L(\mathcal B)$ whose domain is restricted to $V\subset \mathcal B$, which is denoted by $A{\mid_V}$; that is, $A{\mid_V}=A:V\mapsto \mathcal B$.
	\subsection{Proofs of the results given in  Sections \ref{sfcointeg} and \ref{linearization}}
	
	\begin{proof}[\normalfont\textbf{Proof of Proposition \ref{propbasic}}]
		To show (i), we take $0 \neq f\in \mathcal B'$ to both sides of  \eqref{eqi1bn} and obtain $f (\Delta^{d-1}X_t) = f(\tau_0)  + f \Theta(1)  (\sum_{s=1}^t \varepsilon_s) + f (\nu_t)$, $t\geq 0$. Then $\{f (\nu_t)\}_{t \geq 0}$ is stationary since $f$ is Borel measurable and $\{\nu_t\}_{t\geq 0}$ is stationary. Because $\mathbb{E} [(f \Theta(1) \varepsilon_t)^2] = f \Theta(1) C_{\varepsilon_0} f  \Theta(1)$, the second moment of $f \Theta(1)  (\sum_{s=1}^t \varepsilon_s)$ is given by $tf \Theta(1) C_{\varepsilon_0} f  \Theta(1)$, which increases without bound as $t$ grows unless $f \Theta(1) = 0$. Therefore, for $\{f (\Delta^{d-1}X_t)\}_{ t \geq 0}$ to be stationary, $f \Theta(1) = 0$ is required. In this case, a suitable initial condition on $\tau_0$ can be obtained by letting $f(\tau_0)=0$. 	Moreover, one can show without difficulty that $f \Theta(1) = 0$ if and only if \(f \in \Ann (\ran \Theta(1))\), so $\mathfrak C(X) = \Ann(\mathfrak A(X))$. 
		
		To show (ii), suppose that $g \in \mathfrak C(X)$. Note that $x \in \mathcal B$ allows the following unique decomposition, $x = x_{\cl(\mathfrak A(X))} + x_{V}$,  where $x_{\cl(\mathfrak A(X))} \in \cl(\mathfrak A(X))$ and  $x_{V} \in V$.
		From (i) and continuity of bounded linear functionals, we find  that $\mathfrak C(X) = \Ann(\mathfrak A(X)) = \Ann(\cl(\mathfrak A(X)))$. Therefore,  $g \in \mathfrak C(X)$ implies that $g(x) = g(x_V) = g \circ \PR_V (x)$. Now suppose that $g=f \circ \PR_V$. Then clearly $g(X_t) = g(\tau_0) + g (\nu_t)$, $t\geq 0$. This can be made stationary by letting $g(\tau_0) = 0$. 	
	\end{proof}	
	
	\begin{proof}[\normalfont\textbf{Proof of Corollary \ref{propbasic3}}] 
		Under the simplifications discussed in  Section \ref{ssinteg} for the case $\mathcal B=\mathcal H$, Proposition \ref{propbasic}-\ref{propbasic2} implies that $f \in \mathfrak C(X)$ is given by the map  $\langle \cdot,  \PR_V y \rangle$ for $y \in \mathcal H$. Then the stated result follows. 
	\end{proof} 
	
	\begin{proof}[\normalfont\textbf{Proof of Proposition \ref{propa1}}]
		Note that $\PPhi(z)$ may be viewed as the following block operator matrix,	
		\begin{align}
			\PPhi(z) &=  	\scriptsize \left(\begin{array}{c:cccccc} \idb - z{\Phi}_1 & -z{\Phi}_2  & -z{\Phi}_3  & \cdots  &-z{\Phi}_p  \\ \hdashline -z\idb  & \idb & \zero   & \cdots  & \zero \\
				\zero & -z\idb &\idb&\cdots&\zero\\
				\vdots& \vdots & \ddots & \ddots & \vdots  \\
				\zero & \zero    & \zero   & - z\idb &\idb 
			\end{array}\right) \eqqcolon \small\left(\begin{matrix}\PPhi_{[11]}(z)&  \PPhi_{[12]}(z)  \\ \PPhi_{[21]}(z) &  \PPhi_{[22]}(z)\end{matrix}\right),  \label{opmatrix}
		\end{align} \normalsize
		\noindent where $ \PPhi_{[22]}(z): \mathcal B^{p-1} \mapsto \mathcal B^{p-1}$ is invertible for all $z \in \mathbb{C}$.	Define the Schur complement of $ \PPhi_{[22]}(z)$ as $ \PPhi^+_{[11]}(z) \coloneqq  \PPhi_{[11]}(z) -  \PPhi_{[12]}(z)  \PPhi_{[22]}(z)^{-1}  \PPhi_{[21]}(z)$. From a little algebra, we find that $ \PPhi^+_{[11]}(z) =  \Phi(z)$. When  $ \PPhi_{[22]}$ is invertible, $ \PPhi(z)$ is invertible if and only if $ \PPhi^+_{[11]}(z)$ is invertible \citep[Section 2.2]{Bart2007}, 
		so $\sigma(\PPhi) = \sigma(\Phi)$. Furthermore from the Schur's formula in  \citet[p.\ 29]{Bart2007}, we have
		\begin{align}
			\PPhi&(z)^{-1}  = 	\scriptsize   \left(\begin{matrix} \Phi(z)^{-1} &  - \Phi(z)^{-1} \PPhi_{[12]}(z)  \PPhi_{[22]}(z)^{-1}\\  -\PPhi_{[22]}(z)^{-1} \PPhi_{[21]}(z)\Phi(z)^{-1}&  \PPhi_{[22]}(z)^{-1} + \PPhi_{[22]}(z)^{-1}  \PPhi_{[21]}(z) \Phi(z)^{-1} \PPhi_{[12]}(z)  \PPhi_{[22]}(z)^{-1}\end{matrix}\right), \label{schurform} 
		\end{align} \normalsize
		which shows $\Phi(z)^{-1} = \Pi_p\PPhi(z)^{-1}\Pi_p^\ast$. (iii) is deduced from \eqref{schurform} and invertibility of $\PPhi_{[22]}(z)$. 
	\end{proof}
	\subsection{Mathematical appendix to the I(1) representation (Section \ref{sec:repi1})}	 \label{appendixi1}
	\subsubsection{Preliminary results} 
	We provide important preliminary results for the subsequent discussions. Hereafter it should be noted that $\PP$ and $\NN_{j}$ may be alternatively expressed as the  following contour integrals, 
	\begin{equation}
		\PP = \frac{-1}{2\pi i} \int_{\Gamma}{(\idbp-z\PPhi_1)^{-1}}\PPhi_1 { (z-1)^{-j-1} } dz, \quad \NN_{j} = \frac{-1}{2\pi i} \int_{\Gamma}{(\idbp-z\PPhi_1)^{-1}}{ (z-1)^{-j-1} } dz, \quad j \in \mathbb{Z},  \label{rprojection}
	\end{equation} 
	where  $\Gamma \subset  \rho(\PPhi)$ is a clockwise-oriented contour around one such that the only element of $\sigma(\PPhi)$ included inside the contour is one. 	From $\PPhi(z)^{-1}\PPhi(z) = I_p = \PPhi(z)\PPhi(z)^{-1}$,  the identity map may be understood as the power series defined in a punctured neighborhood of one as follows, 
	\begin{equation}
		\sum_{j=-\infty}^\infty \left(\NN_{j-1}\PPhi_1 - \NN_j(\idbp-\PPhi_1)\right)  (z-1)^j  = \idbp = \sum_{j=-\infty}^\infty \left(\PPhi_1 \NN_{j-1} - (\idbp-\PPhi_1)\NN_j\right)  (z-1)^j. \label{idendecom1}
	\end{equation}
	The above identity expansions give us the following relationships, 
	\begin{align}
		\NN_{-1}\PPhi_1 - \NN_0(\idbp-\PPhi_1) \,\,= \,\,&I_p\,\, =\,\,  \PPhi_1 \NN_{-1} - (\idbp-\PPhi_1)\NN_0, \label{idendecomresult0}\\
		\NN_{j-1}\PPhi_1 - \NN_j(\idbp-\PPhi_1) \,\,=\,\, &0\,\, = \,\, \PPhi_1 \NN_{j-1} - (\idbp-\PPhi_1)\NN_j, \quad j \neq 0.  \label{idendecomresult}
	\end{align} 	
	The following lemma collects some essential spectral properties of $\PPhi(z)$. 
	\begin{lemma} \label{lemprime} Suppose that Assumption \ref{assulinear1} holds. Then the following hold.
		\begin{enumerate}[\normalfont(i)]\setlength\itemsep{0em} 
			\item\label{lemprime1} If $\PPhi(z)^{-1}\PPhi_1$ has a pole at $z=1$ of order $\ell$, then $\NN_{-m} \PPhi(1) = 0$ for all $m \geq \ell$. 
			\item\label{lemprime2}  	$\NN_j \PPhi_1 \NN_k = ( 1 - \eta_j - \eta_k) \NN_{j+k+1}$, where $\eta_j = 1\{j \geq 0\}$. Moreover, $\NN_{-1}\PPhi_1$ and $\NN_{-1}$ are projections.
			\item\label{lemprime3}   $\NN_j \PPhi_1 = \PPhi_1 \NN_j$ for all $j \in \mathbb{Z}$. 
			\item \label{lemprime4}  $\PPhi(z)^{-1}$ has a pole of order at most $\ell$ at $z=1$ if and only if (a) $n^{-1}\|G^{\ell-1} (\idbp-G)^n\|_{\op} \to 0$ for some $\ell \in \mathbb{N}$ and (b) $\ran (G^m)$ is closed for some \(m \in \mathbb{N}\) satisfying \(m\geq \ell\), where $G = \PPhi(1)\PP$. 
		\end{enumerate}
	\end{lemma}
	\begin{proof}
		To show (i), we note that  $\idbp = (\idbp - z \PPhi_1)^{-1}(\idbp-\PPhi_1)  - (z-1)(\idbp - z \PPhi_1)^{-1}\PPhi_1$ in a punctured neighborhood of $z=1$.  It is then clear that $(\idbp - z \PPhi_1)^{-1}(\idbp-\PPhi_1)$ must have a pole of order $\ell-1$ if $(\idbp - z \PPhi_1)^{-1}\PPhi_1$ has a pole of order $\ell \geq 1$. We therefore have  $\NN_{-m} \PPhi(1) = 0$ for all $m \geq \ell$.
		
		Our proof of (ii) is similar to those in \citet[p.\ 38]{Kato1995} and \citet[p.\ 119]{amouch2015}.  Let $\Gamma, \Gamma' \subset \rho(\Phi)$ be contours enclosing $z=1$, and assume that $\Gamma'$ is outside $\Gamma$.  Using the generalized resolvent equation \citep[p.50]{Gohberg2013}, it can be shown that 
		\begin{equation}
			\NN_j \PPhi_1 \NN_{k} =\left(\frac{1}{2\pi i}\right)^2 \int_{\Gamma'} \int_{\Gamma} \frac{(\idbp-\lambda \PPhi_1)^{-1}-(\idbp-z \PPhi_1)^{-1}}{(\lambda-z)(z-1)^{j+1}(\lambda-1)^{k+1}} dz d\lambda. \label{lemeqq01} 
		\end{equation}
		From \citet[p.38]{Kato1995}, we may deduce that
		\begin{equation}\label{lemeqq03} 
			\frac{1}{2\pi i } \int_{\Gamma} \frac{(\lambda-z)^{-1}}{(z-1)^{j+1}} dz = \eta_j (\lambda-1)^{-j-1}, \quad \quad \frac{1}{2\pi i } \int_{\Gamma'} \frac{(\lambda-z)^{-1}}{(\lambda-1)^{k+1}} d\lambda = (1-\eta_k) (z-1)^{-k-1}.
		\end{equation}
		Since we may evaluate the integral in any order, the right hand side of \eqref{lemeqq01} can be written as 	\small
		\begin{align*}	
			&\underbracket[0.140ex]{\left(\frac{1}{2\pi i}\right)^2 \int_{\Gamma'} \int_{\Gamma}  \frac{(\idbp-\lambda \PPhi_1)^{-1}}{(\lambda-z)(z-1)^{j+1}(\lambda-1)^{k+1}} dz d\lambda}_{\circled{1}}  \quad -\quad \underbracket[0.140ex]{\left(\frac{1}{2\pi i}\right)^2  \int_{\Gamma} \int_{\Gamma'} \frac{(\idbp-z \PPhi_1)^{-1}}{(\lambda-z)(z-1)^{j+1}(\lambda-1)^{k+1}} d\lambda dz}_{\circled{2}}.  
		\end{align*}\normalsize  
		From \eqref{lemeqq03} and Cauchy's residue theorem, we deduce that $ \small{\circled{1}}  = -\eta_j \NN_{j+k+1}$ and $\small{\circled{2}}  = (\eta_k-1)\NN_{j+k+1}$, from which we find that $\NN_j \PPhi_1 \NN_k = ( 1 - \eta_j - \eta_k) \NN_{j+k+1}$. By putting  $j=-1$ and $k=-1$, we obtain $\NN_{-1}\PPhi_1\NN_{-1} = \NN_{-1}$, which implies that $\NN_{-1}\PPhi_1$ is a projection. We now let $U(z) = zI - \PPhi_1$, and define $\PR_{U} = \frac{1}{2\pi i }\int_{\Gamma} U(z)^{-1}dz$. Then  $\PR_{U}$  is a projection \citep[Lemma 2.1 in Chapter \RN{1}]{Gohberg2013} and $\PR_U = \NN_{-1}$ holds \citep[Remark 3.11]{BS2018}. 
		
		To show (iii), 	note that $\PPhi_1$ and $\idbp-z \PPhi_1$ commute, and this implies that $(\idbp-z \PPhi_1)^{-1} \PPhi_1 = \PPhi_1(\idbp-z \PPhi_1)^{-1}$ commute \citep[Theorem 6.5]{Kato1995}.  We thus have $\NN_j\PPhi_1 = \PPhi_1\NN_j$ for all $j\in \mathbb{Z}$.

		To show (iv), we will first verify that  the following holds for some $\eta>0$,
		\begin{equation}\label{eqprop1}
			-(\idbp - z \PPhi_1)^{-1} \PPhi_1=\sum_{j=1}^{\infty} G^j (z-1)^{-1-j}+ \NN_{-1}\PPhi_1(z-1)^{-1} + \NN_H(z), \quad z \in D_{1+\eta}\setminus\{1\},
		\end{equation}
		where $\NN_H(z)$ is the holomorphic part of the above Laurent series. We deduce that  $\PPhi_1\NN_{-2}  =  \PPhi(1)  \NN_{-1}$ from \eqref{idendecomresult}, and $\PPhi(1)  \NN_{-1} =\PPhi(1)  \NN_{-1} \PPhi_1 \NN_{-1} = \PPhi(1) \PP = G$ from (ii)-(iii). We thus find that $G = \NN_{-2} \PPhi_1 = \PPhi_1\NN_{-2}$.  It is also deduced from (ii) and (iii) that $\NN_{-k} \PPhi_1 = G^{k-1}$ and $G^{k-1} =  \PPhi(1)^{k-1} \PP$ holds for $k \geq 2$, from which we find that \eqref{eqprop1} holds. In order for  \eqref{eqprop1} to converge for $z \in D_{1+\eta}\setminus\{1\}$,  $\lim_{k \to \infty } \|G^k\|_{\op}^{1/k} = 0$ must hold \citetext{\citeauthor{Kato1995}, \citeyear{Kato1995}, pp.\ 180--181}. In this case, (a) and (b) are necessary and sufficient for $G^{\ell}$ to be zero \citetext{\citeauthor{laursen1995}, \citeyear{laursen1995}, Lemma 3 and Corollary 7}.  If $G^{\ell}=0$,  we know, from \eqref{eqprop1} and (i), that  $-(\idbp-z \PPhi_1)^{-1} \PPhi_1$ has a pole of order at most $\ell$ at $z=1$ and $\NN_{-k}  \PPhi(1) = 0$ for all  $k \geq \ell$.  Combining these results with \eqref{idendecomresult}, we have $\NN_{-k-1}\PPhi_1 =  \NN_{-k} \PPhi(1) = 0$ for all $k \geq \ell$. Since $\NN_{-k-1}  = \NN_{-k-1}  \PPhi(1) + \NN_{-k-1} \PPhi_1$, we find that  $\NN_{-k-1} = 0$ for all $k \geq \ell$, so $(\idbp-z\PPhi_1)^{-1}$ has a pole of order at most $\ell$ at $z=1$. Conversely, if $(\idbp-z\PPhi_1)^{-1}$ has a pole of order at most $\ell$ at $z=1$,  we have  $\NN_{-k-1}\PPhi_1 = G^{k} = 0$ for $k \geq \ell$. For \(G^{\ell+1}\) to be nilpotent, (a) and (b) must hold; see \citet[Lemma 3 and Corollary 7]{laursen1995}.
	\end{proof}
	
	\subsubsection{Proofs of the main results}
	\begin{proof}[\normalfont\textbf{Proof of Proposition \ref{cormain}}]
		Since (ii) $\Rightarrow$ (iii) is trivial, we will show that (i) $\Rightarrow $ (ii) and (iii) $\Rightarrow$ (i).  This completes our proof of the equivalence of (i)-(iii). Then we will verify \eqref{eqlaurentform}.
		
		We will show (i) $\Rightarrow$ (ii). Since $\PPhi(z)^{-1}$ has a simple pole at $z=1$, we may deduce from \eqref{idendecomresult} that $ \PPhi(1) \NN_{-1} = 0$. This implies that $\ran \NN_{-1} \subset \ker  \PPhi(1)$. We then find that $\ran \PP \subset \ker  \PPhi(1)$ since $\PP = \NN_{-1}\PPhi_1$. Furthermore, $\ker \PPhi(1) \subset \ran \PP$ holds. To see this, note that if $x_k \in \ker \PPhi(1)$, 
		\begin{equation}
			\PP x_k = - \frac{1}{2\pi i} \int_{\Gamma} (\idbp-z  \PPhi_1)^{-1}  \PPhi_1 x_k dz =   \frac{1}{2\pi i}   \int_{\Gamma} (z-1)^{-1} x_k dz = x_k. \label{rprange}
		\end{equation} We thus conclude that $\ran \PP = \ker  \PPhi(1)$.  Moreover,  $\NN_{-1}	 \PPhi(1) = 0$ is deduced from \eqref{idendecomresult} and  $\idbp-\PP  = \idbp-\PPhi_1\NN_{-1}$ is deduced from Lemma \ref{lemprime}-\ref{lemprime3}. Since $\NR_{-1}\ran  \PPhi(1) = \{0\}$, we have $(\idbp-\PP)\ran  \PPhi(1) = {\ran  \PPhi(1)}$, which implies that $\ran \PPhi(1) \subset \ran(\idbp-\PP)$. On the other hand, for any $x \in \ran (\idbp-\PP)$, we have $x= (\idbp-\PPhi_1\NN_{-1} )x$ since $\idbp-\PP$ is a projection and $\PP=\PPhi_1\NN_{-1}$.  We know from \eqref{idendecomresult0} that $\PPhi_1\NN_{-1} =  \PPhi(1)\NN_0 + \idbp$. Therefore, 	$x = (\idbp-\PPhi_1\NN_{-1}) x = -  \PPhi(1)\NN_0 x$, 	which implies that  $x \in \ran  \PPhi(1)$. Hence, we  have $\ran (\idbp-\PP) \subset \ran  \PPhi(1)$, from which we conclude that $\ran (\idbp-\PP) = \ran  \PPhi(1)$ since $\ran  \PPhi(1)\subset \ran (\idbp-\PP)$ was already shown. To sum up, $\ran \PP = \ker  \PPhi(1)$ and $\ran (\idbp-\PP) = \ran  \PPhi(1)$, which means that $\PP$ is the  projection onto $\ker  \PPhi(1)$ along $\ran  \PPhi(1)$.
		
		To prove  (iii) $\Rightarrow$ (i), we first show that $\PPhi(z)^{-1}$ has a pole of order at most 2 at $z=1$. Due to Lemma \ref{lemprime}-\ref{lemprime4}, it suffices to show that $n^{-1}\|G (\idbp-G)^n\|_{\op} \to 0$ and $\ran (G^2)$ is closed for $G=\PPhi(1)\PP$. Since $\PPhi_1$ and $\PP$ commute,  $G^2 =   \PPhi(1)^2 \PP = \PP  \PPhi(1)^2$.  Moreover, it can be shown that \(  \ran(\PP \PPhi(1)^2) = \ran \PP \cap \ran  \PPhi(1)^2\).  To see this, note that for \(x \in  \ran (\PP \PPhi(1)^2) \), there exists $y\in \mathcal B^p$ such that $x =  \PP\PPhi(1)^2y=\PPhi(1)^2\PP$, where the second equality results from commutativity of \(\PP\) and \( \PPhi(1)^2\). This shows that \(x \in \ran \PP \cap \ran  \PPhi(1)^2\), hence \(  \ran(\PP \PPhi(1)^2) \subset \ran \PP \cap \ran  \PPhi(1)^2\). The reverse inclusion is trivial, so we omit its proof. Since  \( \ran( \PP \PPhi(1)^2) = \ran \PP \cap \ran ( \PPhi(1)^2)\) and \(\ran \PP\) is closed,  $\ran (\PP \PPhi(1)^2)$ is closed if \(\ran ( \PPhi(1)^2)\) is closed. Under the direct sum $\mathcal B^p = \ran  \PPhi(1) \oplus \ker  \PPhi(1)$, $\PPhi(1) \mathcal B^p =  \PPhi(1)[\ran  \PPhi(1) \oplus \ker \PPhi(1)] = \PPhi(1)\ran  \PPhi(1)$ holds.
		That is, $\ran ( \PPhi(1)^2)=\ran \PPhi(1)$, which is closed.  It remains to show that  $n^{-1}\|G (\idbp-G)^n\|_{\op} \to 0$. It can be easily deduced that $(\idbp-G)^n = (\idbp-\PP) + \PPhi_1^n \PP$.	From the fact that $\PPhi_1$ and $\PP$ commute, we have
		\begin{equation}
			n^{-1}\|G(\idbp-G)^n\|_{\op} \leq n^{-1}\| \PPhi(1) \PPhi_1^n\|_{\op} \leq n^{-1}\|\PPhi_1^n{\mid_{\ran  \PPhi(1) }}\|_{\op} \| \PPhi(1)  \|_{\op}.  \label{uppbdd}
		\end{equation} 
		Under Assumption \ref{assulinear1}, we may deduce,  from nearly identical arguments used in \cite{BSS2017} to prove a similar statement, that there exists $k  \in \mathbb{N}$ such that for all $n \geq k$, $\|\PPhi^n_1{\mid_{\ran  \PPhi(1)}}\| < a^n$ for some $a \in (0,1)$. Hence, the upper bound in \eqref{uppbdd} vanishes to zero, and we conclude that $\PPhi(z)^{-1}$ has a pole of at most 2 at $z=1$ under the direct sum $\mathcal B^p = \ran  \PPhi(1) \oplus \ker  \PPhi(1)$. From \eqref{idendecomresult0} and \eqref{idendecomresult}, we have $\NN_{-2} \PPhi(1) = 0$ and $\NN_{-2}\PPhi_1 = \NN_{-1} \PPhi(1)$.   The former (resp.\ the latter) shows that $\NN_{-2}{\mid_{\ran  \PPhi(1)}} = 0$ (resp.\ $\NN_{-2}{\mid_{\ker  \PPhi(1)}} = 0$). Since $\mathcal B^p=\ran \PPhi(1) \oplus \ker  \PPhi(1)$, we conclude that $\NN_{-2} = 0$. Hence, $\PPhi(z)^{-1}$ has a simple pole at $z=1$.  
	
		It remains only for us to show that \eqref{eqlaurentform} holds.  Let $H(z)$ denote the holomorphic part of the Laurent series given in \eqref{laurent0}. Note that if Assumption \ref{assulinear1} holds and $\PPhi(z)^{-1}$ has a simple pole, the Maclaurin series of $(1-z)\PPhi(z)^{-1} = \PP + (1-z) H(z)$  is convergent on $D_{1+\eta}$. Then from Lemma 4.1 of \cite{Johansen1996} (or its obvious extension allowing power series with operator coefficients), it may be deduced that the Maclaurin series of $H(z)$ is convergent on $D_{1+\eta}$. Now we will show $\PP = \NN_{-1}$ to complete our proof. One can deduce from our proof of (i) $ \Rightarrow $  (iii) that $\ran \PP = \ran \NN_{-1}$. Moreover, Lemma \ref{lemprime}-\ref{lemprime2} implies that $\NN_{-1}$ is a projection. Therefore, $\NN_{-1}$ is clearly a projection whose range is equal to $\ran \PP = \ker (\idbp-\PPhi_1)$. Then we find that $\PP = \NN_{-1} \PPhi_1 = \PPhi_1\NN_{-1} = \NN_{-1}$,	where the second equality is from the fact that $\PPhi_1$ and $\NN_{-1}$ commute, and the last equality is from that $\PPhi_1{\mid_{\ker (\idbp-\PPhi_1)}} = \idbp{\mid_{\ker (\idbp-\PPhi_1)}}$ and $\ran \NN_{-1} = \ker (\idbp-\PPhi_1)$.
	\end{proof}

	\begin{proof}[\normalfont\textbf{Proof of Proposition \ref{grti1}}]
		From Propositions \ref{propa1} and  \ref{cormain}, we have $(1-z)\Phi(z)^{-1} = \Pi_p\NN_{-1} \Pi_p^\ast + (1-z)\Pi_p H(z)\Pi_p^\ast$.
		Applying the linear filter induced by $(1-z) \Phi(z)^{-1}$ to \eqref{arlaw}, we obtain $\Delta X_t = \Pi_p\NN_{-1} \Pi_p^\ast \varepsilon_t  + \Delta \nu_t$, where $\nu_t = \Pi_pH(L)\Pi_p^\ast \varepsilon_{t}$ and $H(z) = \sum_{j=0}^\infty H_j z^j$ with $H_j = H^{(j)}/j!$ is convergent on $D_{1+\eta}$.  Clearly, the process given by $\Pi_p\NN_{-1} \Pi_p^\ast \sum_{s=1}^t \varepsilon_s + \nu_t$ is a solution, which is completed by adding a time invariant component $\tau_0$ given as the solution to the homogeneous equation $\Delta {X}_t = 0$.
		
		We then verify the claimed expression of \(\nu_t\) in \eqref{bndecomi1}. Once we show that 
		\begin{equation}
			H_0 =  \idbp-\PP, 	\quad H_j = \PPhi_1 H_{j-1}, \quad j \geq 1, \label{eqthm02}
		\end{equation}
		\noindent then the claimed expression given in  \eqref{bndecomi1} may be easily verified.	First, it can be shown that $H(z) = H(z)(\idbp-\PP)$ holds. To see this, note that $H(z) = -\sum_{j=0}^\infty \NN_j (z-1)^{j}$. Since Lemma \ref{lemprime}-\ref{lemprime2} implies that $\NN_j(\idbp-\PP) = \NN_j - \NN_j\PPhi_1\NN_{-1} = \NN_j$ for $j \geq 0$, we find that $H(z)(\idbp-\PP) = -\sum_{j=0}^\infty \NN_j (\idbp-\PP)(z-1)^{j} = -\sum_{j=0}^\infty \NN_j (z-1)^{j}$. This also shows that  $\PPhi(z)^{-1}(\idbp-\PP)= H(z)$, and
		\begin{equation}
			\idbp-\PP = \PPhi(z)\PPhi(z)^{-1}(\idbp-\PP) =  \PPhi(z) H(z). \label{prop2H}
		\end{equation}
		We then easily deduce that  $H_0= \idbp-\PP$ from \eqref{prop2H} evaluated at $z=0$.	Furthermore,  \eqref{prop2H} can be rewritten as $	(\idbp-\PPhi_1)H(z) - (z-1)\PPhi_1H(z) = \idbp-\PP$, from which we have
		\begin{equation}
			H^{(j)}(z)-j \PPhi_1H^{(j-1)}(z) - z \PPhi_1H^{(j)}(z) = 0, \quad j \geq 1. \label{prop2H3}
		\end{equation}
		Evaluating \eqref{prop2H3} at $z=0$, we obtain $H^{(j)}(0) = j \PPhi_1H^{(j-1)}(0) = j! \PPhi_1 H_{j-1}$,
		which verifies \eqref{eqthm02}. 
		
		It remains to show that if the I(1) condition is not satisfied then the AR($p$) law of motion \eqref{arlaw} does not allow I(1) solutions. This immediately follows from Propositions \ref{propa1} and \ref{cormain}.
	\end{proof}

	\begin{proof}[\normalfont\textbf{Proof of Proposition \ref{propjohanseni1}}]
		Throughout this proof, we write the Laurent series of $\Phi(z)^{-1}$ near $z=1$ as follows: for $d \in \mathbb{N} \cup \{\infty\}$, $	\Phi(z)^{-1} = - \sum_{j=-d}^\infty \NR_j (z-1)^j$. Since it is obvious that (iii) $\Rightarrow$ (ii),  we will only show that (i) $\Rightarrow$ (iii) and (ii) $\Rightarrow$ (i). The whole proof is divided into several parts.    \\[-11.5pt]
		
		\noindent 	{1. (i) $\Rightarrow$  (iii)} : Let $[\ran \Phi(1)]_\complement$ and $[\ker \Phi(1)]_\complement$ be arbitrarily chosen among possible candidates.   We know that $d=1$, $\NR_{-1} = \Pi_p \PP \Pi_p^\ast$, and  $\PP$ is a projection onto $\ker \PPhi(1)$ under the I(1) condition; see Propositions \ref{propa1} and \ref{cormain}. Since  $(x_1,\ldots,x_p) \in \ker \PPhi(1)$ implies that $x_1=\cdots=x_p \in \ker \Phi(1)$, it may be deduced that $\ran \NR_{-1}=\ran \Pi_p \PP \Pi_p^\ast \subset \Pi_p \ker \PPhi(1)  = \ker \Phi(1)$. 	
		Moreover,	from the coefficients of $(z-1)^{-1}$ and $(z-1)^0$ in the identity expansion $\Phi(z)^{-1} \Phi(z) = I$, we know that $\NR_{-1}$ satisfies $\NR_{-1} \Phi(1) = 0$ and $\NR_{-1} \Phi^{(1)}(1) + \NR_{0}\Phi(1) = -I$. From these equations, we observe that $-\NR_{-1} \PR_{[\ran \Phi(1)]_\complement}\Phi^{(1)}(1)(I-\PR_{[\ker \Phi(1)]_\complement}) = I-\PR_{[\ker \Phi(1)]_\complement}$ holds, hence $\ker \Phi(1) \subset \ran \NR_{-1}$. We thus find that $\ran \NR_{-1} = \ker \Phi(1)$, so  $\Lambda_{1,\complement}:\ker \Phi(1) \mapsto [\ran \Phi(1)]_{\complement}$ is an injection. Also, from the coefficients of  $(z-1)^{-1}$ and $(z-1)^0$ in the expansion $\Phi(z)\Phi(z)^{-1}  = I$, we find that $\Phi(1)\NR_{-1}  = 0$ and $\Phi^{(1)}(1)\NR_{-1}  + \Phi(1) \NR_{0} = -I$, which implies that $-\PR_{[\ran \Phi(1)]_\complement}\Phi^{(1)}(1)\NR_{-1}  = \PR_{[\ran \Phi(1)]_\complement}$ holds. Since $\ran \NR_{-1} = \ker \Phi(1)$ was already shown, it is concluded that $\Lambda_{1,\complement}:\ker \Phi(1) \mapsto [\ran \Phi(1)]_{\complement}$ is also a surjection, i.e., it is a bijection. The above arguments do not depend on the choice of $[\ran \Phi(1)]_\complement$ and $[\ker \Phi(1)]_\complement$. Thus (i) $\Rightarrow$ (iii). \\[-8.5pt]
		
		\noindent	{2. (ii) $\Rightarrow$ (i)} : Suppose that the direct sums given in \eqref{direct1} hold and $\Lambda_{1,\complement}$ is invertible for some choice of $[\ran \Phi(1)]_\complement$ and $[\ker \Phi(1)]_\complement$. We define $\mathcal Q_{\complement}:\mathcal B^p \to \mathcal B^p$ as follows,
		\begin{align} 
			\mathcal Q_{\complement} & \scriptsize= \left(\begin{array}{c:cccccc} \PR_{[\ran \Phi(1)]_\complement} &\PR_{[\ran \Phi(1)]_\complement}\sum_{j=2}^{p}\Phi_j  & \PR_{[\ran \Phi(1)]_\complement}\sum_{j=3}^{p} \Phi_j  & \cdots  &\PR_{[\ran \Phi(1)]_\complement}\Phi_p  \\ \hdashline 
				\zero & \zero &\zero&\cdots&\zero\\
				\vdots& \vdots & \vdots & \ddots & \vdots  \\
				\zero & \zero    & \zero   & \cdots&\zero 
			\end{array}\right).\label{opmatrix2}
		\end{align} \normalsize 
		\noindent Then $\mathcal Q_{\complement}$ is a projection on $\mathcal B^p$ and $\mathcal Q_{\complement}\PPhi(1) = 0$. The latter may be verified by noting that \begin{align}
			&\PPhi(1) = \scriptsize\left( \begin{matrix} \Phi(1) & \PPhi_{[12]}(1) \\ 0 & \PPhi_{[22]}(1) \end{matrix} \right)\left( \begin{matrix} I &0 \\ \PPhi_{[22]}(1)^{-1} \PPhi_{[21]}(1) & I_{p-1} \end{matrix} \right)\normalsize, \quad \quad \mathcal Q_{\complement} \scriptsize \left( \begin{matrix} \Phi(1) & \PPhi_{[12]}(1) \\ 0 & \PPhi_{[22]}(1) \end{matrix} \right) = 0,\normalsize
		\end{align} 
		where $\PPhi_{[12]}(\cdot)$, $\PPhi_{[21]}(\cdot)$, $\PPhi_{[22]}(\cdot)$ are defined in \eqref{opmatrix}.
		We thus find that $\ker \mathcal Q_{\complement} \supset \ran \PPhi(1)$. It can also be  shown that $\ker \mathcal Q_{\complement} \subset  \ran \PPhi(1)$. 	To see why, we first note that $\PR_{[\ran \Phi(1)]_\complement} (x_1 + \sum_{j=2}^{p}\Phi_j x_2 +\sum_{j=3}^{p}\Phi_j x_3 + \cdots + \Phi_px_p) = 0$ holds  for any $x = (x_1,\ldots,x_p) \in \ker \mathcal Q_{\complement}$, and  $\PPhi(1)(0,x_2,\ldots,\sum_{j=2}^p x_j) = ( -\sum_{j=2}^{p}\Phi_j x_2 -\sum_{j=3}^{p}\Phi_j x_3 - \cdots - \Phi_px_p , x_2,\ldots,x_p)$. Let $y_1 = -\sum_{j=2}^{p}\Phi_j x_2 -\sum_{j=3}^{p}\Phi_j x_3 - \cdots - \Phi_px_p$. With some algebra, we obtain the following results.
		\begin{align}
			&(x,0,\ldots,0) \in \ran\PPhi(1) \text{ for $x\in \ran \Phi(1)$}, \label{newproofeq04}\\
			&y_1= \PR_{[\ran \Phi(1)]_\complement}y_1 + (I-\PR_{[\ran \Phi(1)]_\complement})y_1, \quad \PR_{[\ran \Phi(1)]_\complement}y_1 = \PR_{[\ran \Phi(1)]_\complement}x_1. \label{newproofeq06}
		\end{align}
		Using \eqref{newproofeq04} and \eqref{newproofeq06}, we find that $(\PR_{[\ran \Phi(1)]_\complement}x_1,x_2,\ldots,x_p) \in \ran \PPhi(1)$. Combining this result with \eqref{newproofeq04}, we conclude that $x=(x_1,\ldots,x_p) \in \ran \PPhi(1)$, so $\ker \mathcal Q_{\complement} \subset  \ran \PPhi(1)$. 
		
		We have shown that $\ker \mathcal Q_{\complement} = \ran \PPhi(1)$, so we know that $\mathcal Q_{\complement}$ is a projection onto a complementary subspace of $\ran \PPhi(1)$. Let $x_k = (x_{1,k},\ldots,x_{p,k}) \in \ker \PPhi(1)$, then it may be easily shown that $x_{1,k}=\cdots=x_{p,k}$ and $x_{1,k} \in \ker \Phi(1)$. With a little algebra and from the fact that $x_{1,k}=(I-\PR_{[\ker \Phi(1)]_\complement})x_{1,k}$, we obtain 
		\begin{align}
			\mathcal Q_{\complement} x_k= \left(-\PR_{[\ran \Phi(1)]_\complement} \Phi^{(1)}(1) (I-\PR_{[\ker \Phi(1)]_\complement})x_{1,k} , 0 , \ldots , 0\right). \label{newproofeq07}
		\end{align}
		Moreover, we may deduce that $\ran \mathcal Q_{\complement} \subset \{(x_1,0,\ldots,0) \in \mathcal B^p : x_1 \in [\ran \Phi(1)]_{\complement} \}$. Combining this result with \eqref{newproofeq07} and invertibility of $\Lambda_{1,\complement}:\ker \Phi(1) \mapsto [\ran \Phi(1)]_{\complement}$, we find that $\mathcal Q_{\complement,R} = \mathcal Q_{\complement}:\ker \PPhi(1) \mapsto [\ran \PPhi(1)]_{\complement}$ is invertible. This implies that $\ker \PPhi(1)$ is a complementary subspace of $\ran \PPhi(1)$, which may be deduced from Fact 4.3 of \cite{Fabian2010} and the fact that the map $D : \ran \PPhi(1) \oplus [\ran \PPhi(1)]_{\complement} \mapsto \ran \PPhi(1) \oplus \ker \PPhi(1)$ given by $	D =  \left(\begin{smallmatrix} \idbp & 0 \\ 0 & \mathcal Q_{\complement,R}^{-1} 
		\end{smallmatrix}\right)$ 	is invertible. \\[-8.5pt]

		\noindent	{3. Formula for $\Upsilon_{-1}$} : For any arbitrary choice of $[\ran \Phi(1)]_\complement$ and $[\ker \Phi(1)]_\complement$, we  found earlier that $-\NR_{-1} \PR_{[\ran \Phi(1)]_\complement}\Phi^{(1)}(1)(I-\PR_{[\ker \Phi(1)]_\complement}) = I-\PR_{[\ker \Phi(1)]_\complement}$  and the map $\Lambda_{1,\complement}: \ker \Phi(1)\mapsto [\ran \Phi(1)]_\complement$ is invertible, from which \eqref{n1formula} follows immediately. 
	\end{proof}
	
	\subsubsection{A detailed discussion on Remark \ref{johansencointeg}} \label{appenBrr}
	From \eqref{n1formula}, we find that $\ran \Upsilon_{-1} = \ker \Phi(1)$. Therefore, for a nonzero $f \in \Ann(\ker \Phi(1))$,
	\begin{equation}
		f(X_t) = f(\tau_0) + f(\nu_t), \quad  t \geq  0.
	\end{equation}
	Let $\tau_0$ be satisfy $f(\tau_0) = 0$. We know from Proposition \ref{cormain} that $\Phi(z)^{-1} =  - \NR_{-1} (z-1)^{-1} - \sum_{j=0}^\infty \NR_{j}(z-1)^j$, $H(z)$ is convergent on $D_{1+\eta}$ for $\eta>0$ (and thus the coefficients of the Maclaurin series of $H(z)$ decay exponentially in norm) and $\NR_0 = - \sum_{j=0}^\infty \Pi_pH_j\Pi_p^\ast$. We therefore only  need to show that $f\NR_0 \neq 0$ to establish I(0)-ness of $\{f(\nu_t)\}_{t\geq 0}$. Under the I(1) condition, we know that $\NR_{-1}\Phi^{(1)}(1) +\NR_{0}\Phi(1) = -I$, which implies that $\ran \NR_{-1} + \ran \NR_{0} = \mathcal B$. In this case, for any $f \in \ker \Phi(1)$, $f \NR_0 = 0$ implies that  $f = 0$, which contradicts our assumption that $f \neq 0$. Therefore $f\NR_0 \neq 0$ is impossible.

	
	\subsection{Mathematical appendix to the I(2) representation (Section \ref{sec:repi2})} \label{appi2}	\subsubsection{Preliminary results}
	We first collect some preliminary results that are useful for the subsequent discussions. 
	\begin{lemma} \label{lemi2} Let everything be as in Section \ref{sec:repi2}. The following holds.
		\begin{enumerate}[\normalfont(i)] \setlength\itemsep{0em}
			\item\label{lemi21} Under Assumption \ref{assulinear1}, $\ran \PPhi(1)$ and $\ker \PPhi(1)$ can be complemented in $\mathcal B^p$. 
			\item \label{lemi21aa} Under Assumption \ref{assulinear2}, the following hold for some  $\mathcal R_{\complement} \subset \mathcal B^p$ and $\mathcal K_{\complement}  \subset \mathcal B^p$, \begin{equation}
			\mathcal B^p = \ran \PPhi(1) \oplus  \mathcal R \oplus \mathcal R_{\complement} =   [\ker \PPhi(1)]_\complement \oplus \mathcal K \oplus \mathcal K_{\complement},  \label{i2directa1} 	
		\end{equation}
	where $\mathcal R = \PR_{[\ran \PPhi(1)]_\complement}\ker \PPhi(1)$. 
			\item \label{lemi21aaa} Under Assumption \ref{assulinear2},  the I(2) condition is equivalent to the following: $\mathcal K \neq \{0\}$ and  $\mathcal B^p =(\ran \PPhi(1) + \ker \PPhi(1)) \oplus \PPhi(1)^g \mathcal K$ holds for any arbitrary choice of $[\ran \PPhi(1)]_\complement$ and $[\ker \PPhi(1)]_\complement$. 
			\item\label{lemi23} Under Assumption \ref{assulinear2}, $\mathcal K \neq \{0\}$ is necessary for $\PPhi(z)^{-1}$ to have a pole of order 2.   
			\item\label{lemi22}  Under Assumption \ref{assulinear2}, the operator \(Q= \PR_{[\ran \PPhi(1)]_\complement}(\idbp -\PR_{[\ker \PPhi(1)]_\complement})\) allows the generalized inverse $Q^g$ satisfying  $\ran Q^g=\mathcal K_\complement$ and $\ker Q^g =\ran \PPhi(1) \oplus \mathcal R_{\complement}$. 	
		\end{enumerate}
	\end{lemma}
	\begin{proof}
		Under Assumption \ref{assulinear1}, we may define the projection $\mathcal Q_\complement$ given  \eqref{opmatrix2}. Since $\ker \mathcal Q_\complement = \ran \PPhi(1)$, it is complemented by $\ran \mathcal Q_\complement$. Moreover, $\ker \PPhi(1)$ is also complemented. To see this, we note that $x=(x_1,\ldots,x_p) \in \ker \PPhi(1)$ implies that $x_1=x_2=\ldots=x_p$ and $x_1 \in \ker \Phi(1)$. Let $\mathcal T_\complement$ be the $(p\times p)$ operator matrix whose entries of the first column are all equal to $I-\PR_{[\ker \Phi(1)]_\complement}$ and all the other entries are equal to zero. Then it can be easily shown that this is a projection defined on $\mathcal B^p$ and its range is equal to $\ker \PPhi(1)$; hence $\ker \PPhi(1)$ is complemented by the kernel of this projection. This completes our proof of (i).

		To show (ii), note that $\ran \PPhi(1)+\ker \PPhi(1) = \ran \PPhi(1) \oplus \mathcal R$. We define $\mathcal S_{\complement}: \mathcal B^p \to \mathcal B^p$ as the block operator matrix obtained by replacing $\PR_{[\ran \Phi(1)]_\complement}$ with $\PR_{\mathsf{R}_\complement}$ in \eqref{opmatrix2}.
		Then it can be easily shown that $\mathcal S_{\complement}$ is a projection on $\mathcal B^p$ and 	$\mathcal S_{\complement}\PPhi(1) = 0$. 	If $x \in \ker \PPhi(1)$, then $x=(x_k,\ldots,x_k)$ for some $x_k \in \ker \Phi(1)$. For such $x$, we have $\mathcal S_{\complement} x = -\PR_{\mathsf{R}_{\complement}} \Phi^{(1)}(1)x_k = 0$,	which is because $\ker \PR_{\mathsf{R}_{\complement}}  = \ran \Phi(1) + \Phi^{(1)}\ker \Phi(1)$. Then we find that $	\ker \mathcal S_{\complement} \supset \ran \PPhi(1)+ \ker \PPhi(1)$.	Moreover, it can be also shown that $\ker \mathcal S_{\complement} \subset \ran\PPhi(1) + \ker\PPhi(1)$. To see why, we let $x = (x_1,\ldots,x_p) \in \ker \mathcal S_{\complement}$. Then $x = (x_1,\ldots,x_p)$ must satisfy $\PR_{\mathsf{R}_{\complement}} (x_1 + \sum_{j=2}^{p}\Phi_j x_2 +\sum_{j=3}^{p}\Phi_j x_3 + \cdots + \Phi_px_p) = 0$.	As in our proof of Proposition \ref{propjohanseni1}, if  we let $y_1 = -\sum_{j=2}^{p}\Phi_j x_2 -\sum_{j=3}^{p}\Phi_j x_3 - \cdots - \Phi_px_p$, then we have $\PR_{\mathsf{R}_{\complement}}y_1 = \PR_{\mathsf{R}_{\complement}}x_1$. It  can also  be shown that $\PPhi(1)(0,x_2,\ldots,\sum_{j=2}^p x_j) = (y_1,x_2,\ldots,x_p)$. We note that
		\begin{equation}
			(x,0,\ldots,0) \in \ran \PPhi(1) + \ker \PPhi(1) \text{ for $x \in \ran \Phi(1) + \Phi^{(1)}(1)\ker \Phi(1)$},\label{opmatrix2adadd}
		\end{equation}
		which is because for any arbitrary $u \in \mathcal B$ and $w \in \ker \Phi(1)$ we have $\PPhi(1)(v-w,v-2w,\ldots,v-pw) + (w,w,\ldots,w) = ( \Phi(1)v + \Phi^{(1)}(1)w,0,\ldots,0)$. Combining all these results with the fact that $\ran (I-\PR_{\mathsf{R}_{\complement}})= \ran \Phi(1) + \Phi^{(1)}(1)\ker\Phi(1)$, we may deduce that $(x_1,\ldots,x_p) \in \ran \PPhi(1)+\ker \PPhi(1)$. We have shown that $\mathcal S_{\complement}$ is a projection whose kernel is $\ran \PPhi(1) + \ker \PPhi(1)=\PPhi(1) \oplus \mathcal R$, hence $\PPhi(1) \oplus \mathcal R$ is complemented by $\ran \mathcal S_{\complement}$. The latter direct sum of \eqref{i2directa1} clearly holds for  $\mathcal K_\complement = \ker \PPhi(1) \cap [\ran \PPhi(1)]_\complement$.    
		
			To show (iii), we let $\PPhi(1)^g_1$ and $\PPhi(1)^g_2$ be the generalized inverse operators depending on two different choices of $[\ran \PPhi(1)]_\complement$ and $[\ker \PPhi(1)]_\complement$. Let $V=\ran \PPhi(1) + \ker \PPhi(1)$, $V_1=\PPhi(1)_1^g\mathcal K$, and $V_2 = \PPhi(1)_2^g\mathcal K$. We know from  \citet[Theorem 1.7.14 and Corollary 3.2.16]{megginson1998} that two complementary subspaces of $\ker \PPhi(1)$ must be isomorphic, which implies that $AV_1 = V_2$ and $V_1 = A^{-1}V_2$  for some invertible map $A \in \mathcal L(V_1,V_2)$ and its inverse $A^{-1}  \in \mathcal L(V_2,V_1)$. Define the map $D : V \oplus V_1 \mapsto V \oplus V_2$ given by $D =  \left(\begin{smallmatrix} \idbp & 0 \\ 0 & A \end{smallmatrix}\right)$, which is obviously invertible and its inverse $D^{-1} : V \oplus V_2 \mapsto V \oplus V_1$ is given by $D^{-1} =  \left(\begin{smallmatrix} \idbp & 0 \\ 0 & A^{-1} \end{smallmatrix}\right)$. Note that $D(V)=V$, $D(V_1)=V_2$,  $D^{-1}(V)=V$, and  $D^{-1}(V_2)=V_1$. We then deduce from Fact 4.3 of \cite{Fabian2010} that $\mathcal B^p = V \oplus V_1$ holds if and only if $\mathcal B^p = V \oplus V_2$ holds. This completes the proof.

		To show (iv), suppose that $\PPhi(z)^{-1}$ has a pole of order 2 at $z=1$ and $\mathcal K = \{0\}$. It can be shown that \(\ker \PPhi(1) \subset \ran \PP \subset \ker (\PPhi(1)^2)\) holds under Assumption \ref{assulinear1}, where the first inclusion follows from \eqref{rprange}. To see why the second inclusion holds, we  note that  $\NN_{-2} = \PPhi(1) \PP$ (Lemma \ref{lemprime}-\ref{lemprime4}), and then deduce from \eqref{idendecomresult} that $\PPhi(1) \NN_{-2} = 0$. From these results, we find that $\PPhi(1) \NN_{-2} = \PPhi(1)^2 \PP = 0$, which proves the second inclusion. Since  $\mathcal K = \{0\}$ implies that  $\ker \PPhi(1) = \ker (\PPhi(1)^2)$, we find that \(\ran \PP=  \ker \PPhi(1)\) and so $\NN_{-2} =\PPhi(1) \PP =0$. This contradicts our assumption that $\PPhi(z)^{-1}$ has a pole of order 2 at $z=1$, so  $\mathcal K \neq \{0\}$.
		
		To show (v), we  note that $\mathcal B^p = \ran \PPhi(1) \oplus \ran Q \oplus \mathcal R_{\complement}$ under the former direct sum \eqref{i2directa1}. Moreover, from the latter direct sum of \eqref{i2directa1}, we have $\mathcal B^p = \ker Q \oplus \mathcal K_\complement$, where $\ker Q = \mathcal K \oplus [\ker \PPhi(1)]_\complement$. Therefore, the generalized inverse $Q^g$ exists and $\ran Q^g=\mathcal K_\complement$ and $\ker Q^g =\ran \PPhi(1) \oplus \mathcal R_{\complement}$, see Appendix \ref{sginverse}.
	\end{proof}
	\begin{lemma} \label{lemi2add} Suppose that  Assumption \ref{assulinear2} holds. Then for any arbitrary choice of $\mathsf{R}_\complement$ and $\mathsf{K}_\complement$, the operator $\mathrm{M}_1=\PR_{[\ran \Phi(1)]_\complement}\Phi^{(1)}(1)(I-\PR_{[\ker \Phi(1)]_\complement})$ allows the generalized inverse $\mathrm{M}_1^g$ satisfying 
		\begin{equation} \label{eqlem01a}\mathrm{M}_1\mathrm{M}_1^g= (I-\PR_{\mathsf{R}_{\complement}})\PR_{[\ran\Phi(1)]_\complement},\quad \mathrm{M}_1^g\mathrm{M}_1= \PR_{\mathsf{K}_\complement}, \quad \ran \mathrm{M}_1^g=\mathsf{K}_\complement, \quad \ker \mathrm{M}_1^g = \ran \Phi(1) \oplus \mathsf{R}_\complement,\end{equation}
		where $\PR_{\mathsf{K}_\complement}$ is the projection onto $\mathsf{K}_\complement$ along $[\ker \Phi(1)]_\complement \oplus \mathsf{K}$. 
	\end{lemma}
	\begin{proof}
		Under the direct sums given in  \eqref{newdirec}, $\mathcal B = \ran \Phi(1) \oplus \ran \mathrm{M}_1 \oplus \mathsf{R}_\complement$ holds. Note also that $\ker \mathrm{M}_1 = [\ker \Phi(1)]_\complement \oplus \mathsf{K}$, hence we have $\mathcal B=\ker \mathrm{M}_1 \oplus \mathsf{K}_{\complement}$. (To see why, recall that $\mathcal B = \ker\Phi(1) \oplus  [\ker\Phi(1)]_\complement$  and $\ker\Phi(1) = \mathsf{K} \oplus \mathsf{K}_{\complement}$.) We then know from Appendix \ref{sginverse} that  $\mathrm{M}_1$  allows the generalized inverse $\mathrm{M}_1^g$ whose range is $\mathsf{K}_{\complement}$ and kernel is $\ran \Phi(1) \oplus \mathsf{R}_{\complement}$, from which \eqref{eqlem01a} immediately follows.   
	\end{proof}
	\subsubsection{Proofs of the main results}
	\begin{proof}[\normalfont\textbf{Proof of Proposition \ref{grti2p}}] Under the direct sums  \eqref{i2directa1} given in Lemma \ref{lemi2}-\ref{lemi21aa}, we let $\PR_{\mathcal R_{\complement}}$ denotes the projection onto $\mathcal R_\complement$ along $\ran\PPhi(1)\oplus \mathcal R$, which is well defined.  The whole proof is divided into several parts. \\[-11.5pt]
		
		\noindent {1. Necessity of the I(2) condition} : If the I(2) condition is not satisfied, then it must be the case either
		\begin{align}
			&\left(\ran\PPhi(1) + \ker \PPhi(1)\right) \cap \PPhi(1)^g \mathcal K \neq \{0\} \quad \text{or} \quad  \mathcal B^p \neq \ran\PPhi(1) + \ker \PPhi(1) +\PPhi(1)^g \mathcal K. \label{pfeq01}
		\end{align}
		It will be shown that \eqref{pfeq01} is false if  $\hspace{-0.0em}\PPhi(z)^{-1}\hspace{-0.0em}$ has a pole of order 2 at $\hspace{-0.1em}z=1\hspace{-0.1em}$. 	From \eqref{idendecomresult0}-\eqref{idendecomresult}, we have
		\begin{align}
			\NN_{-2} \PPhi(1) = &0 = \PPhi(1) \NN_{-2}, \label{i2eq01} \\
			\NN_{-2}\PPhi_1 - \NN_{-1}\PPhi(1)= &0 = \PPhi_1 \NN_{-2} - \PPhi(1)  \NN_{-1}, \label{i2eq02}\\
			\NN_{-1}\PPhi_1 - \NN_{0}\PPhi(1) = &\idbp  = \PPhi_1\NN_{-1} -\PPhi(1)\NN_{0}.  \label{i2eq03}
		\end{align}
		From \eqref{i2eq01}, we observe that \vspace{-0.7em}
		\begin{align}
			&\NN_{-2} = \PPhi_1 \NN_{-2} = \NN_{-2}\PPhi_1,  \label{n2com}\\ 
			&\NN_{-2} = \NN_{-2}  \PR_{[\ran\PPhi(1)]_\complement}. \label{i2eq055}
		\end{align}
		Restricting both sides of \eqref{i2eq02} to \(\ker\PPhi(1) \), we find that
		\begin{align}
			\NN_{-2}{\mid_{\ker\PPhi(1)}} = 0 \quad \Leftrightarrow  \quad 	\NN_{-2}  \ker\PPhi(1) = \{0\}. \label{i2eq06}
		\end{align}
		We note that $\PPhi(1)^g$ exists, and then deduce from \eqref{i2eq02} and \eqref{n2com}  that
		\begin{align}
			&\NN_{-1} (\idbp-\PR_{[\ran\PPhi(1)]_\complement}) =\NN_{-2}\PPhi(1)^g, \quad \PR_{[\ker\PPhi(1)]_\complement}\NN_{-1}  =  \PPhi(1)^g\NN_{-2}. \label{i2eq08} 
		\end{align}
		\eqref{i2eq055} implies that post-composing both sides of the latter equation in \eqref{i2eq08} with $\PR_{[\ran\PPhi(1)]_\complement}$ changes nothing, hence
		\begin{equation}
			\NN_{-1} \PR_{[\ran\PPhi(1)]_\complement} =  \PPhi(1)^g\NN_{-2} +  (\idbp-\PR_{[\ker\PPhi(1)]_\complement})\NN_{-1} \PR_{[\ran\PPhi(1)]_\complement}. \label{i2eq09}
		\end{equation}
		From the former equation in \eqref{i2eq08} and \eqref{i2eq09}, we find that
		\begin{align}
			\NN_{-1} =& \NN_{-2} \PPhi(1)^g+ \PPhi(1)^g\NN_{-2} + (\idbp-\PR_{[\ker\PPhi(1)]_\complement})\NN_{-1}   \PR_{[\ran\PPhi(1)]_\complement}. \label{i2eq10}
		\end{align}
		Restricting both sides of \eqref{i2eq10} to \(\mathcal K\), we obtain $\NN_{-1}{\mid_{\mathcal K}} = \NN_{-2} \PPhi(1)^g{\mid_{\mathcal K}}$, which is due to \eqref{i2eq06}. 	Further,   \(\NN_{-1}{\mid_{\mathcal K}} = \idbp{\mid_{\mathcal K}}\) is deduced by restricting both sides of \eqref{i2eq03} to $\mathcal K$, so we find that $\NN_{-1}{\mid_{\mathcal K}} = \NN_{-2}\PPhi(1)^g{\mid_{\mathcal K}} = \idbp{\mid_{\mathcal K}}$. 	Given \(\PR_{\mathcal R_{\complement}}\), the former direct sum given in \eqref{i2directa1}, equations \eqref{i2eq055}-\eqref{i2eq06}, and the fact that  $\ran \PPhi(1) \oplus \mathcal R = \ran \PPhi(1) + \ker \PPhi(1)$,  we conclude that $\NN_{-2}  = \NN_{-2} \PR_{\mathcal R_{\complement}}$. Therefore,
		\begin{align}
			\NN_{-2} \PR_{\mathcal R_{\complement}} \PPhi(1)^g{\mid_{\mathcal K}} = \idbp{\mid_{\mathcal K}}. \label{i2eq20}
		\end{align}
		This proves injectivity of $\PR_{\mathcal R_{\complement}} \PPhi(1)^g : \mathcal K \mapsto \mathcal R_{\complement}$, hence the former condition in \eqref{pfeq01} is impossible. 

		Pre-composing both sides of the latter equation in \eqref{i2eq02} with $\PR_{\mathcal R_{\complement}}\PPhi(1)^g$ and  using \eqref{n2com}, we have
		\begin{equation}
			\PR_{\mathcal R_{\complement}}\PPhi(1)^g\NN_{-2} = \PR_{\mathcal R_{\complement}}\PR_{[\ker \PPhi(1)]_\complement} \NN_{-1} = \PR_{\mathcal R_{\complement}}\NN_{-1}, \label{i2eq04a}
		\end{equation}
		where the last equality is established from the fact that the kernel of $\PR_{\mathcal R_{\complement}}$ is $ \ran \PPhi(1) + \ker \PPhi(1)$ and thus $\PR_{\mathcal R_{\complement}}(\idbp-\PR_{[\ker \PPhi(1)]_\complement})=0$ . Similarly, pre-composing both sides of the latter equation of \eqref{i2eq03} with $\PR_{\mathcal R_{\complement}}$,
		\begin{align}
			\PR_{\mathcal R_{\complement}} \PPhi_1\NN_{-1} = \PR_{\mathcal R_{\complement}} + \PR_{\mathcal R_{\complement}}\PPhi(1)\NN_{0} = \PR_{\mathcal R_{\complement}}.\label{i2eq05a}
		\end{align} 
		Note that $\PR_{\mathcal R_{\complement}} \NN_{-1} = \PR_{\mathcal R_{\complement}}\PPhi(1) \NN_{-1} +  \PR_{\mathcal R_{\complement}} \PPhi_1 \NN_{-1}$, where the first term is $0$. We therefore deduce from \eqref{i2eq05a} that $\PR_{\mathcal R_{\complement}} \NN_{-1}= \PR_{\mathcal R_{\complement}}$. Combining this  with \eqref{i2eq04a}, we obtain $\PR_{\mathcal R_{\complement}}\PPhi(1)^g\NN_{-2} = \PR_{\mathcal R_{\complement}}$, which implies sujectivity of $\PR_{\mathcal R_{\complement}} \PPhi(1)^g : \mathcal K \mapsto \mathcal R_{\complement}$. Hence, the latter condition in \eqref{pfeq01} cannot hold. \\[-8.5pt] 
		
		\noindent {2. Sufficiency of the I(2) condition} : 
		Suppose that $\PPhi(z)^{-1}$ has a pole of order $m \geq 3$ at $z=1$. We deduce that $\NN_{-m} = \NN_{-m}\PPhi_1 = \PPhi_1 \NN_{-m}$ from \eqref{idendecomresult}, and also $\NN_{-m} = G^{m-1} = \PPhi(1)^{m-1} \PP$ from Lemma \ref{lemprime}-\ref{lemprime4}. If the I(2) condition holds, any $x \in \mathcal B^p$ can be written as $x=x_{\ran} + x_{\ker} + \PPhi(1)^g x_{\mathcal K}$, where $x_{\ran}  \in \ran \PPhi(1)$, $x_{\ker} \in \ker\PPhi(1)$ and $ x_{\mathcal K} \in \mathcal K$. From the definitions of $\ran \PPhi(1)$ and $\mathcal K$, we also  know that there exist  $y_{1}\in \mathcal B^p$ and $y_2 \in \mathcal B^p$ satisfying $x_{\ran} =  \PPhi(1)y_1$, $x_{\mathcal K} =  \PPhi(1)y_2$, and $\PPhi(1)^2y_2 = 0$. We thus find that
		\begin{equation}
			\PPhi(1)^{m-1}\PP x =  \PP\PPhi(1)^{m-1} (x_{\ran} + x_{\ker} +\PPhi(1)^g x_{\mathcal K}) =\PP\PPhi(1)^{m} y_1 +  \PP\PPhi(1)^{m-1}(x_{\ker} +y_{2}) = 0, \label{eqi201}
		\end{equation} 
		where the first equality comes from commutativity of $\PP$ and $\PPhi(1)$. It is then deduced from \eqref{eqi201} that $\PPhi(1)^{m-1}\PP \mathcal B = \{0\}$, so  $\PPhi(1)^{m-1}\PP = 0$. That is, $\NN_{-m} = 0$ is concluded, which, however, contradicts our assumption that $\PPhi(z)^{-1}$ has a pole of order $m \geq 3$. In addition, $\mathcal K \neq \{0\}$ implies that $\mathcal B \neq \ran\PPhi(1) \oplus \ker \PPhi(1)$, which excludes the existence of a simple pole (see Proposition \ref{cormain}).  \\[-8.5pt]

		\noindent {3. Formula for $\NN_{-1}$} :	
		Note that $\NN_{-1} = \PPhi(1) \NN_{-1} + \PPhi_1 \NN_{-1}$. We also know from \eqref{n2com} that  $\NN_{-2}\PPhi_{1}^2 = \NN_{-2}\PPhi_{1}= \NN_{-2}$. Combining this with \eqref{i2eq02}, it is deduced that  $\NN_{-1} = \PPhi(1)  \NN_{-1} + \PPhi_1 \NN_{-1} =  \NN_{-2} + \PP$. \\[-8.5pt]
		
		\noindent {4. $\ran \NN_{-2}=\mathcal K$} : The result immediately follows from \eqref{i2eq20} and invertibility of  $\PR_{\mathcal R_{\complement}} \PPhi(1)^g: \mathcal K \to \mathcal R_{\complement}$.\\[-8.5pt] 
		
		\noindent 	{5. Holomorphicity of $(1-z)^2\PPhi(z)^{-1}$ and $H(z)$ on $D_{1+\eta}$} : 
		We know that the Maclaurin series of $(1-z)^2\PPhi(z)^{-1}$ is convergent on $D_{1+\eta}$. Then from  an obvious extension of Lemma 4.1 of \cite{Johansen1996}, it may be deduced that the Maulaurin series of $H(z)$ is convergent on $D_{1+\eta}$. 	
		\end{proof}
		
		\begin{proof}[\normalfont\textbf{Proof of Proposition \ref{grt2pp}}]
			From Propositions \ref{propa1} and  \ref{grti2p}, we have $(1-z)^2 \Phi(z)^{-1} = - \Pi_p\NN_{-2}\Pi_p^\ast + \Pi_p\mathcal (\NN_{-2}  + \PP)\Pi_p^\ast(1-z) + \Pi_p H(z)\Pi_p(1-z)^2$. Applying the linear filter induced by $(1-z)^2 \Phi(z)^{-1}$ to \eqref{arlaw}, we obtain $\Delta^2 X_t = -\Pi_p\NN_{-2}\Pi_p^\ast  \varepsilon_t + \Pi_p\mathcal (\NN_{-2}  + \PP)\Pi_p^\ast (\varepsilon_t - \varepsilon_{t-1}) + (\Delta \nu_t - \Delta \nu_{t-1})$. We then may deduce that solutions to \eqref{arlaw} satisfies \eqref{bndecom3} for some $\tau_0$ and $\tau_1$. The claimed  expression of \(\nu_t\) can be verified from nearly identical arguments used in our proof of Proposition \ref{grti1}.  Moreover, we may deduce from Propositions \ref{propa1} and \ref{grti2p} that  \eqref{arlaw} does not allow I(2) solutions if the I(2) condition is not satisfied. 	
		\end{proof}

		\begin{proof}[\normalfont\textbf{Proof of Proposition \ref{propjohanseni2}}]
			We write the Laurent series of $\Phi(z)^{-1}$ around $z=1$ as follows: for $d \in \mathbb{N} \cup \{\infty\}$,  $\Phi(z)^{-1} = - \sum_{j=-d}^\infty \NR_j (z-1)^j$. Since it is obvious that (iii) $\Rightarrow$ (ii),  we will only show that (i) $\Rightarrow$ (iii) and (ii) $\Rightarrow$ (i). The whole proof is divided into several parts.  \\[-11.5pt]
			
			\noindent 	{1. (i) $\Rightarrow$ (iii)} : Let $\mathsf{R}_\complement$ and  $\mathsf{K}_\complement$  be arbitrarily chosen among possible candidates.
			If $(x_1,\ldots,x_p) \in \mathcal K$, we know from (4.40) of \cite{BSS2017} that $x_1=	\cdots = x_p$, and there exists $y_1 \in \mathcal B$ such that $	\Phi(1)y_1 = -\Phi^{(1)}(1)x_1$. This implies that $x_1 \in \mathsf{K}$. 	Since $\ran \NN_{-2} = \mathcal K$, $\ran \NR_{-2} =\ran \Pi_p \NN_{-2} \Pi_p^\ast \subset \mathsf{K}$ holds.  	Under the I(2) condition, we know that $\Phi(z)^{-1}$ has a pole of order 2 at $z=1$ and may deduce the following from the coefficients in the identity expansion $\Phi(z)^{-1} \Phi(z) = I= \Phi(z) \Phi(z)^{-1}$, \vspace{-0.7em}
			\begin{align}
				\NR_{-2} \Phi(1) = &0 = \Phi(1)\NR_{-2},  \label{eq01rema}\\
				\NR_{-2} \Phi^{(1)}(1) + \NR_{-1} \Phi(1) =&0= \Phi^{(1)}(1)\NR_{-2}  + \Phi(1)\NR_{-1},   \label{eq01remb}\\
				\NR_{-2} {\Phi^{(2)}(1)}/{2} + \NR_{-1} \Phi^{(1)}(1) + \NR_0 \Phi(1) = -&I = {\Phi^{(2)}(1)}\NR_{-2}/{2}  + \Phi^{(1)}(1) \NR_{-1} +\Phi(1) \NR_0.  \label{eq01remc}
			\end{align}  
			From some algebra similar to that in proof of Theorem 4.2 of \cite{BS2018}, we find that $\NR_{-2}(I-\PR_{\mathsf{R}_\complement})=0$ and $-\NR_{-2} \PR_{\mathsf{R}_{\complement}} \mathrm{M}_2\PR_{\mathsf{K}} = \PR_{\mathsf{K}}$. This implies that $\mathsf{K} \subset \ran \NR_{-2}$ (hence $\ran \NR_{-2} = \mathsf{K}$ has been established) and thus $\Lambda_{2,\complement}=\PR_{\mathsf{R}_{\complement}} \mathrm{M}_2\PR_{\mathsf{K}}: \mathsf{K} \mapsto \mathsf{R}_{\complement}$ is an injection.

			Now from the latter equation in \eqref{eq01remb} and the properties of $\Phi(1)^g$, we may deduce that 
			\begin{align}\label{eq01rembadd2}
				\Phi^{(1)}(1) \NR_{-1}  = - \Phi^{(1)}(1)\Phi(1)^g\Phi^{(1)}(1)\NR_{-2}+\Phi^{(1)}(1) (I-\PR_{[\ker \Phi(1)]_\complement})\NR_{-1}. 
			\end{align}
			From the latter equation in \eqref{eq01remc} and \eqref{eq01rembadd2}, we obtain  $\PR_{\mathsf{R}_{\complement}}\mathrm{M}_2\NR_{-2} + \PR_{\mathsf{R}_{\complement}}\Phi^{(1)}(1)(I-\PR_{[\ker \Phi(1)]_\complement})\NR_{-1}  = - \PR_{\mathsf{R}_{\complement}}$.
			Note that $\PR_{\mathsf{R}_{\complement}}\Phi^{(1)}(1)\ker \Phi(1) = \{0\}$, which implies that $\PR_{\mathsf{R}_{\complement}}\Phi^{(1)}(1)(I-\PR_{[\ker \Phi(1)]_\complement})\NR_{-1} = 0$. We thus find that $	\PR_{\mathsf{R}_{\complement}}\mathrm{M}_2\NR_{-2}  = -\PR_{\mathsf{R}_{\complement}}$. Since $\ran \NR_{-2} = \mathsf{K}$ as shown above, this implies that $\Lambda_{2,\complement}=\PR_{\mathsf{R}_{\complement}} \mathrm{M}_2\PR_{\mathsf{K}}: \mathsf{K} \mapsto \mathsf{R}_{\complement}$ is also a surjection, i.e.,\ it is a bijection. The above arguments do not depend on a specific choice of $\mathsf{R}_\complement$ and $\mathsf{K}_\complement$, and thus (i) $\Rightarrow$ (iii).	
			\\[-8.5pt]
			
			\noindent	{2.(ii)  $\Rightarrow$ (i)} : Suppose that Assumption \ref{assulinear2} hold and $\Lambda_{2,\complement}$ is invertible for some choice of $\mathsf{R}_\complement$ and $\mathsf{K}_\complement$. It suffices to show that the I(2) condition holds for some choice of $\PPhi(1)^g$ (Lemma \ref{lemi21}-\ref{lemi21aaa}). Let \vspace{-0.0em}
			\begin{align} 
				\PPhi(1)^g = \scriptsize \left(\begin{matrix}  \Phi(1)^g & -\Phi(1)^g \PPhi_{[12]}\PPhi_{[22]}^{-1} \\ -\PPhi_{[22]}^{-1}\PPhi_{[21]}\Phi(1)^g & \PPhi_{[22]}^{-1}+\PPhi_{[22]}^{-1}\PPhi_{[21]}\Phi(1)^g \PPhi_{[12]}\PPhi_{[22]}^{-1}   \end{matrix}\right). \label{opmatrix3}
			\end{align}\normalfont
			From the factorization formula (2.3) of \cite{Bart2007}, it may be easily shown that $\PPhi(1)\PPhi(1)^g = \idbp-\mathcal Q_\complement$ for $\mathcal Q_\complement$, the bounded projection whose kernel is equal to $\ran \PPhi(1)$, given in \eqref{opmatrix2}. Moreover, $\PPhi(1)^g\PPhi(1)=\idbp-\mathcal T_\complement$ for $\mathcal T_\complement$, the bounded projection whose range is equal to $\ker \PPhi(1)$ given in our proof of Lemma \ref{lemi2}-\ref{lemi21}. Therefore, $\PPhi(1)^g$ given in \eqref{opmatrix3} is the generalized inverse of $\PPhi(1)$ for some $[\ran \PPhi(1)]_\complement$ and $[\ker \PPhi(1)]_\complement$ depending on the choice of $[\ran \Phi(1)]_\complement$ and $[\ker \Phi(1)]_\complement$.   As in our proof of Lemma \ref{lemi2}-\ref{lemi21aa}, we let $\mathcal S_{\complement}: \mathcal B^p \to \mathcal B^p$ denote the block operator matrix obtained by replacing $\PR_{[\ran \Phi(1)]_\complement}$ with $\PR_{\mathsf{R}_\complement}$ in the definition of $\mathcal Q_\complement$ given by \eqref{opmatrix2}, and let $\mathcal S_{\complement,[1,2]}: \mathcal B^{p-1} \mapsto \mathcal B$ denote the upper-right block of $\mathcal S_{\complement}$.	We already showed that $\mathcal S_{\complement}$ is a projection onto a complementary subspace of $\ran \PPhi(1) + \ker \PPhi(1)$ in our proof of Lemma \ref{lemi2}-\ref{lemi21aa}. Then  using the formula (2.3) of \cite{Bart2007}, we find that $\mathcal S_{\complement} \PPhi(1)^g  = \left( \begin{smallmatrix}  C_1& C_2 \\ 0 & 0
			\end{smallmatrix} \right)$, where 
			\begin{align}
				&C_1\hspace{-0.2em} =\hspace{-0.2em} \PR_{\mathsf{R}_{\complement}}\Phi(1)^g - \mathcal S_{\complement,[12]} \PPhi_{[22]}^{-1}\PPhi_{[21]}\Phi(1)^g, \notag\\
				&C_2 \hspace{-0.2em} =\hspace{-0.2em} \PR_{\mathsf{R}_{\complement}}\Phi(1)^g \PPhi_{[12]}(1) \PPhi_{[22]}(1)^{-1}\hspace{-0.2em} +\hspace{-0.2em}  \mathcal S_{\complement,[12]}(\PPhi_{[22]}(1)^{-1} \hspace{-0.2em}+\hspace{-0.2em} \PPhi_{[22]}(1)^{-1}\PPhi_{[21]}(1)\Phi(1)^g\PPhi_{[12]}(1)\PPhi_{[22]}(1)^{-1}),\notag\end{align} and 	$\PPhi_{[12]}(\cdot)$, $\PPhi_{[21]}(\cdot)$, $\PPhi_{[22]}(\cdot)$ are defined in \eqref{opmatrix}. Consider $x_k = (x_{1,k},\ldots,x_{p,k}) \in \mathcal K$, then we find that  $x_{1,k}=\cdots=x_{p,k}$ and $x_{1,k} \in \mathsf{K}$. For such $x_k$, we obtain the following after a tedious algebra,
			\begin{equation}
				\mathcal S_{\complement} \PPhi(1)^g x_k= \left(  -\PR_{\mathsf{R}_{\complement}}\mathrm{M}_2 \PR_{\mathsf{K}} x_{1,k} , 0, \ldots , 0\right). \label{newproofeq07add}
			\end{equation}
			From the definition of $\mathcal S_{\complement}$, it may be deduced that $\ran \mathcal S_{\complement} \subset \{(x_1,0,\ldots,0) \in \mathcal B^p : x_1 \in \mathsf{R}_{\complement}\}$. Combining this result with \eqref{newproofeq07add} and invertibility of the map $\Lambda_{2,\complement}=\PR_{\mathsf{R}_{\complement}}\mathrm{M}_2\PR_{\mathsf{K}}:\mathsf{K} \mapsto \mathsf{R}_{\complement}$, we may conclude that $\mathcal{S}_{\complement}:\PPhi(1)\mathcal K \mapsto \mathcal R_{\complement}$ is invertible. This implies that $\PPhi(1)^g \mathcal K$ is a complementary subspace of $\ran \PPhi(1) + \ker \PPhi(1)$ by a similar argument that we used in our proof of Proposition \ref{propjohanseni1}.	From the above proof, we know that $\mathcal K=\{0\}$ is impossible since it implies $\mathsf{K}=0$.  \\[-8.5pt]
			
			\noindent	{3. Formula for $\Upsilon_{-2}$} :  We know $-\NR_{-2} \PR_{\mathsf{R}_{\complement}} \mathrm{M}_2\PR_{\mathsf{K}} = \PR_{\mathsf{K}}$ holds under the I(2) condition. Since the map $\Lambda_{2,\complement}: \mathsf{K} \mapsto \mathsf{R}_{\complement}$ is invertible and $\Upsilon_{-2} = -\NR_{-2}$, the desired results given by \eqref{eqn2laurent} are easily obtained. \\[-8.5pt]
			
			\noindent	{4. Formula for $\Upsilon_{-1}$} :  We first establish some preliminary results. According to the direct sums given in Assumption \ref{assulinear2} and for any arbitrary choice of the complementary subspaces therein, 
			we have
			\begin{equation}\label{idendecoma}
				I = (I-\PR_{[\ran \Phi(1)]_\complement}) + (I-\PR_{\mathsf{R}_{\complement}})\PR_{[\ran \Phi(1)]_\complement} +  \PR_{\mathsf{R}_{\complement}}.
			\end{equation}
			As shown in our previous proof, for $x_k = (x_{1,k},\ldots,x_{p,k}) \in \mathcal B^p$,   $x_k \in \ker \PPhi(1)$ (resp.\ $x_k \in \mathcal K$) implies that $x_{1,k} = \cdots=x_{p,k}$ and $x_{1,k} \in \ker \Phi(1)$ (resp.\ $x_{1,k} \in \mathsf{K}$). 
			Based on the identity \eqref{idendecoma}, we will obtain explicit expressions of $\NR_{-1} (I-\PR_{[\ran \Phi(1)]_\complement})$, $\NR_{-1}(I-\PR_{\mathsf{R}_{\complement}})\PR_{[\ran \Phi(1)]_\complement}$ and $\NR_{-1}\PR_{\mathsf{R}_{\complement}}$. In the subsequent proof, we need \eqref{eq01remb}, \eqref{eq01remc} and the following  obtained from the coefficient of $(z-1)^{1}$ in the identity expansion $\Phi(z)^{-1} \Phi(z) = I = \Phi(z)\Phi(z)^{-1}$,
			\begin{align}
				\NR_{-2}{\Phi^{(3)}(1)}/{6}\hspace{-0.2em} +\hspace{-0.2em} \NR_{-1}{\Phi^{(2)}(1)}/{2} \hspace{-0.2em}+\hspace{-0.2em} \NR_0 \Phi^{(1)}(1) = &0 = 	{\Phi^{(3)}(1)}\NR_{-2}/6\hspace{-0.2em} +\hspace{-0.2em} {\Phi^{(2)}(1)}\NR_{-1}/2 \hspace{-0.2em}+ \hspace{-0.2em}\Phi^{(1)}(1)\NR_0. \label{remremeq03}
			\end{align}
			From \eqref{eq01remb}, we have
			\begin{equation}
				\NR_{-1}(I-\PR_{[\ran \Phi(1)]_\complement})  = - \NR_{-2}\Phi^{(1)}(1)\Phi(1)^g.  \label{remremeq01a}
			\end{equation}
			From \eqref{eq01remc} and the identity $\NR_{-1}\Phi^{(1)}(1) =  \NR_{-1}(I-\PR_{[\ran \Phi(1)]_\complement})\Phi^{(1)}(1) +  \NR_{-1}\PR_{[\ran \Phi(1)]_\complement}\Phi^{(1)}(1)$, we have $[\NR_{-2}{\Phi^{(2)}}(1)/2 + \NR_{-1}(I-\PR_{[\ran \Phi(1)]_\complement})\Phi^{(1)}(1) + \NR_{-1}\PR_{[\ran \Phi(1)]_\complement}\Phi^{(1)}(1)](I-\PR_{[\ker\Phi(1)]_\complement}) = -(I-\PR_{[\ker\Phi(1)]_\complement})$. Substituting \eqref{remremeq01a} into this equation and using the definition of $\mathrm{M}_1$, we have
			\begin{equation} \label{remremeq03a}
				\NR_{-1}\mathrm{M}_1 =-(I + \NR_{-2} \mathrm{M}_2)(I-\PR_{[\ker\Phi(1)]_\complement}),  
			\end{equation}
			Post-composing  both sides of \eqref{remremeq03a} with $\mathrm{M}_1^g$ and using the fact that $\mathrm{M}_1\mathrm{M}_1^g= (I-\PR_{\mathsf{R}_{\complement}})\PR_{[\ran\Phi(1)]_\complement}$ (Lemma \ref{lemi2add}), we obtain that 
			\begin{equation} \label{remremeq04}
				\NR_{-1}(I-\PR_{\mathsf{R}_{\complement}})\PR_{[\ran\Phi(1)]_\complement} = - (I+ \NR_{-2} \mathrm{M}_2)\mathrm{M}_1^g.
			\end{equation}
			Now from \eqref{remremeq03},  we have $[\NR_{-2}{\Phi^{(3)}(1)}/{6} + \NR_{-1}\Phi^{(2)}(1)/2 + \NR_0\Phi^{(1)}(1)]\PR_{\mathsf{K}} = 0$.	We note from the definition of $\mathsf{K}$ that  $\NR_0\Phi^{(1)}(1)\PR_{\mathsf{K}} = \NR_0(I-\PR_{[\ran\Phi(1)]_\complement})\Phi^{(1)}(1)\PR_{\mathsf{K}}$, and deduce from \eqref{eq01remc} that $\NR_0 (I-\PR_{[\ran\Phi(1)]_\complement}) = -\Phi(1)^g - \NR_{-2}{\Phi^{(2)}(1)}\Phi(1)^g /2- \NR_{-1}\Phi^{(1)}(1)\Phi(1)^g$. 	Combining these results, we obtain
			\begin{equation}
				\NR_{-1}\mathrm{M}_2\PR_{\mathsf{K}} = - [\NR_{-2}\mathrm{M}_3 - \NR_{-2}\mathrm{M}_2\Phi(1)^g\Phi^{(1)}(1) - \Phi(1)^g\Phi^{(1)}(1)]\PR_{\mathsf{K}}.\label{remremeq07}
			\end{equation}
			Since $\NR_{-1} \PR_{\mathsf{R}_{\complement}} \mathrm{M}_2\PR_{\mathsf{K}} = [\NR_{-1}\mathrm{M}_2-  \NR_{-1}(I-\PR_{[\ran\Phi(1)]_\complement})\mathrm{M}_2- \NR_{-1} (I-\PR_{\mathsf{R}_{\complement}})\PR_{[\ran\Phi(1)]_\complement}\mathrm{M}_2]\PR_{\mathsf{K}}$, and  $\NR_{-1}(I-\PR_{[\ran\Phi(1)]_\complement})$ and $\NR_{-1} (I-\PR_{\mathsf{R}_{\complement}})\PR_{[\ran\Phi(1)]_\complement}$ are given in \eqref{remremeq01a} and \eqref{remremeq04}, we have
			\begin{align}
				\NR_{-1}\PR_{\mathsf{R}_{\complement}} \mathrm{M}_2\PR_{\mathsf{K}}  =& - [\NR_{-2}\mathrm{M}_3 - \NR_{-2} \mathrm{M}_2\Phi(1)^g\Phi^{(1)}(1) - \Phi(1)^g\Phi^{(1)}(1)]\PR_{\mathsf{K}}  \notag \\ & + \NR_{-2}\Phi^{(1)}(1)\Phi(1)^g\mathrm{M}_2\PR_{\mathsf{K}}+ (I+\NR_{-2}\mathrm{M}_2)\mathrm{M}_1^g\mathrm{M}_2\PR_{\mathsf{K}}. \label{eqlast}
			\end{align}
			By post-composing both sides of  \eqref{eqlast} with $\NR_{-2}$, a  formula for $\NR_{-1}\PR_{\mathsf{R}_{\complement}}$ is obtained. Combining this with \eqref{remremeq01a},  \eqref{remremeq04}, and the fact that $\Upsilon_{-2} = -\NR_{-2}$ and $\Upsilon_{-1} = \NR_{-1}$, we obtain
			\begin{align} 
				\Upsilon_{-1} =& - \mathrm{M}_1^g  + \left(\Phi(1)^g \Phi^{(1)}(1) +   \mathrm{M}_1^g \mathrm{M}_2  \right) \Upsilon_{-2} + \Upsilon_{-2} \left(\Phi^{(1)}(1) \Phi(1)^g	+ \mathrm{M}_2  \mathrm{M}_1^g \right)  \notag\\ &+ \Upsilon_{-2}  \left(\mathrm{M}_3 -\mathrm{M}_2\Phi(1)^g\Phi^{(1)}(1) -\Phi^{(1)}(1)\Phi(1)^g \mathrm{M}_2 - \mathrm{M}_2  \mathrm{M}_1^g \mathrm{M}_2 \right)\Upsilon_{-2}.\label{eqn2laurent2app}
			\end{align}	
			Using the fact that  $\ran \mathrm{M}_1^g=\mathsf{K}_\complement$ and $\Upsilon_{-2}=\PR_{\mathsf{K}}\Upsilon_{-2}\PR_{\mathsf{R}_\complement}$,  the desired results are obtained from  \eqref{eqn2laurent2app}. 			
		\end{proof}

	\subsubsection{Supplementary results to Proposition \ref{grti2p}} \label{app:companioni2}
	We here characterize the principal part of the Laurent series $\PPhi(z)^{-1}$ in more detail. Let $\PR_{[\ran \PPhi(1)]_\complement}$, $\PR_{[\ker \PPhi(1)]_\complement}$ and $\PPhi(1)^g$ be defined as in Section \ref{sec:repi2companion}. Given the direct sum conditions given in \eqref{i2directa1}, which holds under Assumption \ref{assulinear2}, we  let $\PR_{\mathcal R_{\complement}}$ (resp.\ $\PR_{\mathcal K}$) be the projection onto $\mathcal R_\complement$ (resp.\ $\mathcal K$) along $\ran \PPhi(1)\oplus \mathcal R$ (resp.\ $ [\ker \PPhi(1)]_\complement \oplus\mathcal K_\complement$). Let $\widetilde{\Lambda}$ be the map given by ${\PR}_{\mathcal R_{\complement}} \PPhi(1)^g : \mathcal K \mapsto \mathcal R_{\complement}$ and let $Q=\PR_{[\ran \PPhi(1)]_\complement}(\idbp -\PR_{[\ker \PPhi(1)]_\complement}) $. Given the direct sums \eqref{direct1a}, the generalized inverse $Q^g$ is well defined and satisfies that $\ran Q^g = \mathcal K_\complement$ and $\ker Q^g = \ran \PPhi(1)\oplus \mathcal R_\complement$ (Lemma \ref{lemi2}-\ref{lemi22}). We will show that  $\widetilde{\Lambda}: \mathcal K \mapsto \mathcal R_{\complement}$ is invertible and  \(\NN_{-2}\) satisfies
	\begin{equation}
		(\idbp-\PR_{\mathcal K})	\NN_{-2}=\NN_{-2}(\idbp-\PR_{\mathcal R_{\complement}})=0,\quad \quad 	\NN_{-2}:\mathcal R_\complement \mapsto \mathcal K = \widetilde{\Lambda}^{-1},\label{opn2}
	\end{equation}
	and $\PP =   (\idbp-\Gamma_r)Q^g(\idbp-\Gamma_\ell) + \Gamma_\ell(\idbp-\Gamma_r) +  \Gamma_r$, where \(\Gamma_\ell = \PPhi(1)^g \NN_{-2} \) and  \(\Gamma_r = \NN_{-2}\PPhi(1)^g \).
	\\[-8.5pt]
	
	\noindent {1. Formula for $\NN_{-2}$} : Under the I(2) condition,  we know from our proof of Proposition \ref{grti2p} that \eqref{i2eq20} holds, $\NN_{-2}=\NN_{-2}\PR_{\mathcal R_{\complement}}$ and  $\widetilde{\Lambda} =\PR_{\mathcal R_{\complement}} \PPhi(1)^g: \mathcal K \to \mathcal R_{\complement}$ is bijective, from which \eqref{opn2} immediately follows.  \\[-8.5pt]
	
	\noindent {2. Formula for $\PP$} :
	We will verify the claimed formula for $\PP$. From the former direct sum in \eqref{i2directa1}, we have 
	\begin{equation}
		\idbp = (\idbp-\PR_{[\ran\PPhi(1)]_\complement}) + (\idbp-\PR_{\mathcal R_{\complement}})\PR_{[\ran\PPhi(1)]_\complement} + \PR_{\mathcal R_{\complement}}. \label{idendecomnew}
	\end{equation}
	Thus \(\PP =\PP(\idbp-\PR_{[\ran\PPhi(1)]_\complement})+\PP (\idbp-\PR_{\mathcal R_{\complement}})\PR_{[\ran\PPhi(1)]_\complement}+\PP \PR_{\mathcal R_{\complement}}\), we will obtain an expression for each summand.  
	
	Pre-composing both sides of the former equation in \eqref{i2eq02} with $\PPhi_1$, we obtain $\PPhi_1\NN_{-2}\PPhi_1 = \PPhi_1 \NN_{-1} \PPhi(1)$. We then deduce the following from \eqref{n2com} and the fact that $\PPhi_1 \NN_{-1} = \PP$,
	\begin{align} \label{i2eq077}
		\PP (\idbp-\PR_{[\ran\PPhi(1)]_\complement}) = \NN_{-2} \PPhi(1)^g.
	\end{align}
	Using the identity \(\idbp = (\idbp-\PR_{[\ran\PPhi(1)]_\complement}) + \PR_{[\ran\PPhi(1)]_\complement} \) and post-composing both sides of the former equation in \eqref{i2eq03} with $\idbp-\PR_{[\ker\PPhi(1)]_\complement}$, we find that $	\PP Q =  (\idbp - \PP(\idbp-\PR_{[\ran\PPhi(1)]_\complement}))(\idbp-\PR_{[\ker\PPhi(1)]_\complement})$.  Then from \eqref{i2eq077}, we obtain $\PP Q = (\idbp - \NN_{-2} \PPhi(1)^g)(\idbp-\PR_{[\ker\PPhi(1)]_\complement})$. It may be deduced from Lemma \ref{lemi2}-\ref{lemi22} that $QQ^g = (\idbp-\PR_{\mathcal R_{\complement}})\PR_{[\ran\PPhi(1)]_\complement}$. Thus, post-composing both sides with $Q^g$, we find that 
	\begin{align}
		\PP &(\idbp-\PR_{\mathcal R_{\complement}})\PR_{[\ran\PPhi(1)]_\complement}=  \left(\idbp - \NN_{-2} \PPhi(1)^g\right)Q^g. \label{n1eq02}
	\end{align}	

	Post-composing both sides of the former equation in \eqref{i2eq03} with  \(\PPhi(1)^g\PR_{\mathcal K}\), we obtain $	\PP \PPhi(1)^g\PR_{\mathcal K}  - \NN_{0}(\idbp-\PR_{[\ran\PPhi(1)]_\complement})\PR_{\mathcal K}  = \PPhi(1)^g\PR_{\mathcal K}$.	Note that \((\idbp-\PR_{[\ran\PPhi(1)]_\complement})\PR_{\mathcal K} = \PR_{\mathcal K}\). We also deduce
	from  \eqref{idendecomresult} that  $\NN_0 \PPhi_1 = \NN_{1}\PPhi(1)$,
	which implies that $\NN_0\PR_{\mathcal K} = 0$.	Combining these results, we find that $	\PP\PPhi(1)^g\PR_{\mathcal K}   = \PPhi(1)^g\PR_{\mathcal K}$.	Then from the identity decomposition \eqref{idendecomnew}, we have
	\begin{align}
		\PP\PR_{\mathcal {R_{\complement}}}\PPhi(1)^g\PR_{\mathcal K}  =  &\PPhi(1)^g\PR_{\mathcal K}  -	\PP (\idbp-\PR_{[\ran\PPhi(1)]_\complement})\PPhi(1)^g\PR_{\mathcal K} -  \PP(\idbp-\PR_{\mathcal R_\complement})\PR_{[\ran\PPhi(1)]_\complement}\PPhi(1)^g\PR_{\mathcal K}.\label{i22eq03}  
	\end{align}
	We then substitute \eqref{i2eq077} and \eqref{n1eq02} into \eqref{i22eq03}, and then  post-compose both sides with $\NN_{-2}$, noting that $\PR_{\mathcal K}\NN_{-2}=\NN_{-2}$ due to \eqref{opn2}. This gives us  an explicit formula for $\PP\PR_{\mathcal {R_{\complement}}}$. Combining this with   \eqref{i2eq077} and \eqref{n1eq02}, the claimed formula for $\PP$ can be obtained.
		\subsubsection{A detailed discussion on Remark \ref{app:polycointeg0}} \label{app:polycointeg}
		Since $\ran \Upsilon_{-2} = \mathsf{K}$, $f \in \Ann(\mathsf{K})$ is a cointegrating functional. For any nonzero $f \in \Ann( \Upsilon_{-2})$, we may deduce from the formula of $\Upsilon_{-1}$ given in Proposition \ref{propjohanseni2}  that $f \Upsilon_{-1}$ is equal to 
		\begin{equation}
			- f\mathrm{M}_1^g(I-\mathrm{M}_2\Upsilon_{-2}) + f\Phi(1)^g\Phi^{(1)}(1)\Upsilon_{-2}. \label{eqeqeq001a}
		\end{equation} 
		Using the expression of $f\Upsilon_{-1}$ given by \eqref{eqeqeq001a}, we will  show that the following holds, 
		\begin{equation}
			f\Upsilon_{-1} = 0 \quad \Leftrightarrow	\quad f \in \Ann(\Phi(1)^g\Phi^{(1)}(1) \mathsf{K}) \cap \Ann(\ker \Phi(1)).  \label{eqeqeq001b} 
		\end{equation}
		Since $\ran \Upsilon_{-2} = \mathsf{K}$, the second term in \eqref{eqeqeq001a} is zero if and only if $f \in \Ann(\Phi(1)^g\Phi^{(1)}(1)\mathsf{K})$.
		It thus only remains to show that $f \in \Ann(\mathsf{K}_{\complement})$ because $\ker \Phi(1) = \mathsf{K} \oplus \mathsf{K}_\complement$ and $f \in \Ann(\mathsf{K})$. To see this, we first show that $\ran (\mathrm{M}_1^g(I-\mathrm{M}_2\Upsilon_{-2})) = \mathsf{K}_\complement$. Since $\ran \mathrm{M}_1^g = \mathsf{K}_{\complement}$,   $\ran (\mathrm{M}_1^g(I-\mathrm{M}_2\Upsilon_{-2})) \subset \mathsf{K}_\complement$. Note that for any $V\subset\mathcal B$, $\mathrm{M}_1^g(I-\mathrm{M}_2\Upsilon_{-2})V \subset \mathrm{M}_1^g(I-\mathrm{M}_2\Upsilon_{-2})\mathcal B \subset  \mathsf{K}_\complement$ holds. Thus, if there is a subset $V$ such that $\mathrm{M}_1^g(I-\mathrm{M}_2\Upsilon_{-2})V = \mathsf{K}_\complement$, then  $\ran (\mathrm{M}_1^g(I-\mathrm{M}_2\Upsilon_{-2})) = \mathsf{K}_\complement$ holds. Let $V = (I-\PR_{\mathsf{R}_{\complement}})[\ran \Phi(1)]_{\complement}$. We then know from \eqref{eqn2laurent} that $\mathrm{M}_2\Upsilon_{-2}V = \{0\}$. Moreover,  $\ran \mathrm{M}_1^g = \mathrm{M}_1^g V$ holds since $\ker \mathrm{M}_1^g=\ran \Phi(1) \oplus \mathsf{R}_{\complement}$. From these results, we find that $\mathrm{M}_1^g(I-\mathrm{M}_2\Upsilon_{-2})V =  \mathrm{M}_1^g V = \mathsf{K}_\complement$, so $\ran (\mathrm{M}_1^g(I-\mathrm{M}_2\Upsilon_{-2})) = \mathsf{K}_{\complement}$. We thus conclude that $f \in \Ann(\mathsf{K}_{\complement})$ if and only if  the first term in \eqref{eqeqeq001a} is zero.  
		
		We have shown that \eqref{eqeqeq001b} holds for $f \in \Ann(\mathsf{K})$. We know from Proposition \ref{propjohanseni2} that (ignoring $\tau_0$ and $\tau_1$ without loss of generality) a nonzero element $f \in \Ann(\mathsf{K})$ satisfies either of the following: (i) $f\Upsilon_{-1} \neq 0$ and $f(X_t) = f\Upsilon_{-1} \left(\sum_{s=1}^t \varepsilon_s \right)+ f(\nu_t)$ or (ii) $f(X_t) = f(\nu_t)$.
		In case (i), $\{f(X_t)\}_{t\geq 0}$ is I(1) obviously. In case (ii), $f(X_t)$ is I(0) under our I(2) condition. To see this, note that we know from Proposition \ref{grti2p} that $\Phi(z)^{-1} = -\NR_{-2}(z-1)^{-2} - \NR_{-1} (z-1)^{-1} - \sum_{j=0}^\infty \NR_{j}(z-1)^j$, $H(z)$ is convergent on $D_{1+\eta}$ for $\eta>0$ (and thus the coefficients of the Maclaurin series of $H(z)$ decay exponentially in norm) and $\NR_0 = - \sum_{j=0}^\infty \Pi_pH_j\Pi_p^\ast$. As in Appendix \ref{appenBrr}, it suffices to show that $f\NR_0 \neq 0$ to establish the desired I(0)-ness.  Under the I(2) condition, we know from \eqref{eq01remc} that $\ran \NR_{-2} + \ran \NR_{-1} + \ran \NR_{0}  = I$. Since $f\NR_{-2}=f\NR_{-1}=0$ in case (ii), $f\NR_{0} = 0$ implies $f=0$, which contradicts our assumption that $f\neq 0$. We thus find that $f\NR_0 \neq 0$, so the cointegrating behavior of I(2) solutions is characterized as stated.  	
		
		\subsubsection{A detailed discussion on Remark \ref{app:polycointeg1}} \label{app:polycointeg2}
		For a nonzero $f \in \Ann(\mathsf{K})$, we  deduce from \eqref{eqeqeq001a} that $f(X_t) - f(\Phi(1)^g\Phi^{(1)}(1) \Delta X_t)$ is given by
		\begin{equation} \label{eq01rem}
			-f\left(\mathrm{M}_1^g(I-\mathrm{M}_2\Upsilon_{-2})\sum_{s=1}^t \varepsilon_s\right) + f(\dt{\nu}_t),
		\end{equation} 
		where $\dt{\nu}_t = \nu_t - \Phi(1)^g\Phi^{(1)}(1)\varepsilon_t + \Delta \nu_t$. 
		As shown in Appendix \ref{app:polycointeg},  $\ran (\mathrm{M}_1^g(I-\mathrm{M}_2\Upsilon_{-2})) = \mathsf{K}_{\complement}$. Since $\ker \Phi(1) = \mathsf{K} \oplus \mathsf{K}_\complement$ and $f \in \Ann(\mathsf{K})$, $f \notin \Ann(\ker \Phi(1))$ implies that $f \notin \Ann(\mathsf{K}_\complement)$. Therefore, the sequence given in \eqref{eq01rem} cannot be stationary for $f \notin \Ann(\ker \Phi(1))$. On the other hand, if $f \in \Ann(\ker \Phi(1))$ then the first term in \eqref{eq01rem} is zero. Hence, it only remains to prove I(0)-ness of   $\{f(\dt{\nu}_t)\}_{t\geq 0}$ for the desired result. The summability condition for the I(0) property is satisfied, which can be easily shown. We rewrite $\dt{\nu}_t$ as $\dt{\nu}_t = \sum_{j=0}^\infty \Psi_j \varepsilon_{t-j}$, and  find that $ \sum_{j=0}^\infty \Psi_j = -\NR_0 - \Phi(1)^g \Phi^{(1)}(1)\NR_{-1}$. If this operator is nonzero, then $\{f(\dt{\nu}_t)\}_{t\geq 0}$ is I(0).  Suppose by contradiction that  $f\NR_0 =  -f\Phi(1)^g \Phi^{(1)}(1)\NR_{-1}$.  Pre-composing \eqref{eq01remc} with $f\Phi(1)^g$  and using the fact that $\Phi(1)^g\Phi(1) =\PR_{[\ker \Phi(1)]_\complement}$, we find that 
		\begin{equation*}
			f\Phi(1)^g \Phi^{(2)}(1)\NR_{-2}/2 + f\Phi(1)^g \Phi^{(1)}(1)\NR_{-1} + f\PR_{[\ker \Phi(1)]_\complement} \NR_0 = -f\Phi(1)^g. 
		\end{equation*}
		In the above, $f\Phi(1)^g \Phi^{(1)}(1)\NR_{-1} = -f\NR_0$ under our supposition, and $f\PR_{[\ker \Phi(1)]_\complement}\NR_0=f\NR_{0}$ since $f(I-\PR_{[\ker \Phi(1)]_\complement}) = 0$ for any $f \in \Ann(\ker \Phi(1))$. We thus find that $f\Phi(1)^g\Phi^{(2)}(1)\NR_{-2}/2 = -f\Phi(1)^g$. From our  expression of $\NR_{-2} = - \Upsilon_{-2}$ given in \eqref{eqn2laurent}, we know that $\NR_{-2} \ran \Phi(1) = \{0\}$ and so
		\begin{equation}
			0 = f\Phi(1)^g\Phi^{(2)}(1)\NR_{-2}x = -f\Phi(1)^gx, \quad  \text{for all $x \in \ran \Phi(1)$}.\label{eq01a}
		\end{equation} 
		From the properties of a generalized inverse, we have  $\Phi(1)^g \ran \Phi(1) = [\ker \Phi(1)]_{\complement}$ and thus find that \eqref{eq01a} holds if and only if $f \in \Ann([\ker \Phi(1)]_{\complement})$. Since $f \in \Ann(\ker \Phi(1))$ and $\mathcal B = \ker \Phi(1)\oplus [\ker \Phi(1)]_{\complement}$, we conclude that $f=0$, which contradicts our assumption that $f$ is nonzero.
	
\end{document}